\documentclass[11pt,reqno]{article}

\usepackage{amssymb}
\usepackage{amsthm}
\usepackage{amsmath}
\usepackage{fullpage,mathrsfs,verbatim,graphicx,paralist}
\usepackage{picins}
\usepackage[breaklinks]{hyperref}

\renewcommand{\setminus}{\smallsetminus}
\newcommand\remove[1]{}

\newcommand{\mnote}[1]{$\ll$\textsf{#1 --Manor}$\gg$}

\renewcommand{\le}{\leqslant}
\renewcommand{\ge}{\geqslant}
\renewcommand{\leq}{\leqslant}
\renewcommand{\geq}{\geqslant}

\DeclareMathOperator{\cost}{cost}
\newcommand{\e}{\varepsilon}
\newcommand{\R}{\mathbb{R}}
\newcommand{\Z}{\mathbb{Z}}
\newcommand{\E}{\mathbb{E}}

\newcommand{\F}{\mathcal{F}}

\newcommand{\N}{\mathbb{N}}

\newcommand{\Lip}{\mathrm{Lip}}

\newcommand{\Mid}{\mathrm{Mid}}

\newcommand{\f}{\varphi}

\DeclareMathOperator{\diam}{diam} 

\DeclareMathOperator{\lca}{\mathsf{lca}}
\newcommand{\Matousek}{Matou\v{s}ek}

\newcommand{\tb}{{|\hspace{-0.9pt}|\hspace{-0.9pt}|}}

\DeclareMathOperator{\dist}{dist}

\newtheorem{theorem}{Theorem}[section]

\newtheorem{lemma}[theorem]{Lemma}

\newtheorem{proposition}[theorem]{Proposition}

\newtheorem{corollary}[theorem]{Corollary}

\newtheorem{definition}[theorem]{Definition}

\newtheorem{remark}[theorem]{Remark}

\newcommand{\eqdef}{\stackrel{\mathrm{def}}{=}}
\newcommand{\MM}{\mathscr M}
\newcommand{\NN}{\mathscr N}

\newcommand{\ANCA}{\mathcal{A}}
\newcommand{\ANCB}{\mathcal{B}}
\graphicspath{{figs/}}

\begin{document}

\title{Markov convexity and local rigidity of distorted metrics}

\author{
Manor Mendel\thanks{Supported by ISF grant 221/07, BSF grant
2006009, and a gift from Cisco research center.}\\Open University of
Israel\\{\tt manorme@openu.ac.il} \and Assaf Naor\thanks{Supported
by NSF grants CCF-0635078 and CCF-0832795, BSF grant 2006009, and
the Packard Foundation.}\\New York University\\{\tt
naor@cims.nyu.edu} }

\date{}
\maketitle

%\thispagestyle{empty}

%\begin{comment}
\begin{abstract}
It is shown that a Banach space admits an equivalent norm whose modulus of uniform convexity has power-type $p$ if and only if it is Markov $p$-convex.
%, thus obtaining a purely metrical characterization of $p$-convexity of Banach spaces.
Counterexamples are constructed to natural questions related to isomorphic uniform convexity of metric spaces, showing in particular that tree metrics fail to have the dichotomy property.
\end{abstract}
%\end{comment}

\tableofcontents

\section{Introduction}

A Banach space $(X,\|\cdot \|_X)$ is said to be finitely
representable in a Banach space $(Y,\|\cdot\|_Y)$ if there exists a
constant $D<\infty$ such that for every finite dimensional linear
subspace $F\subseteq X$ there is a linear operator $T: F\to Y$
satisfying $\|x\|_X\le \|Tx\|_Y\le D\|x\|_X$ for all $x\in F$. In
1976 Ribe~\cite{Ribe76} proved that if two Banach spaces $X$ and $Y$
are uniformly homeomorphic, i.e., there is a bijection $f:X\to Y$
such that $f$ and $f^{-1}$ are uniformly continuous, then $X$ is finitely
representable in $Y$ and vice versa. This remarkable theorem
motivated what is known today as the ``Ribe program": the search for
purely metric reformulations of basic linear concepts and invariants
from the local theory of Banach spaces. This research program was put
forth by Bourgain in 1986~\cite{Bourgain-trees}.

Since its inception, the
Ribe program attracted the work of many mathematicians, and led to
the development of several satisfactory metric theories that extend
important concepts and results of Banach space theory; see the introduction of~\cite{MN-cotype-full} for a historical discussion. So far, progress on the Ribe program has come hand-in-hand
with applications to metric geometry, group theory,
functional analysis, and computer science. The present paper
contains further progress in this direction: we obtain a metric
characterization of $p$-convexity in Banach spaces, derive some of
its metric consequences, and construct unexpected counter-examples
which indicate that further progress on the Ribe program can uncover
nonlinear phenomena that are markedly different from their Banach space counterparts.
%do not have analogues in Banach space theory.
In doing so, we answer questions posed by Lee-Naor-Peres and
Fefferman, and improve a theorem of Bates, Johnson, Lindenstrauss,
Preiss and Schechtman. These results, which will be explained in detail below,  were announced in~\cite{charlie-soda}.

%These results will be explained in detail
%below.

For $p\ge 2$, a Banach space $(X,\|\cdot\|_X)$ is said to be $p$-convex if there
exists a norm ${\tb} \cdot{\tb} $ which is equivalent to $\|\cdot \|_X$
(i.e., for some $a,b>0$, $a\|x\|_X\le {\tb} x{\tb} \le b\|x\|_X$ for all
$x\in X$), and a constant $K>0$ satisfying:
\begin{equation}\label{eq:def p-conv}
{\tb} x{\tb} ={\tb} y{\tb} =1\implies \left|\left|\left|
\frac{x+y}{2}\right|\right|\right|\le 1-K{\tb} x-y{\tb} ^p.
\end{equation}
$X$ is called superreflexive if it is $p$-convex for some $p\ge 2$
(historically, this is not the original definition of
superreflexivity\footnote{James' original definition of
superreflexivity is that a Banach space $X$ is superreflexive if
``its local structure forces reflexivity", i.e., if every Banach
space $Y$ that is finitely representable in $X$ must be reflexive.
Enflo's renorming theorem states that superreflexivity is equivalent
to having an equivalent norm $|||\cdot|||$ that is uniformly convex,
i.e., for every $\e\in (0,1)$ there exists $\delta>0$ such that if
$|||x|||=|||y|||=1$ and $|||x-y|||=\e$ then $|||x+y|||\le
2-\delta$.}, but it is equivalent to it due to a deep theorem of
Pisier~\cite{Pisier-martingales}, which builds on important work of
James~\cite{James72} and Enflo~\cite{Enf73}). For concreteness, we
recall (see, e.g., \cite{BCL}) that $L_p$ is $2$-convex for $p\in
(1,2]$ and $p$-convex for $p\in [2,\infty)$.

Ribe's theorem implies that $p$-convexity, and hence also
superreflexivity, is preserved under uniform homeomorphisms. The
first major success of the Ribe program is a famous theorem of
Bourgain~\cite{Bourgain-trees} which obtains a metrical
characterization of superreflexivity as follows.
\begin{theorem}[Bourgain's metrical characterization of superreflexivity~\cite{Bourgain-trees}]\label{thm:bourgain tree}
Let $B_n$ be the
complete unweighted binary tree of depth $n$, equipped with the
natural graph-theoretical metric. Then a Banach space $X$ is
superreflexive if and only if
\begin{equation}\label{eq:bourgain's tree}
\lim_{n\to\infty} c_X(B_n)=\infty.
\end{equation}
\end{theorem}
Here, and in what follows, given two metric spaces $(\mathscr M,d_\mathscr M)$,
 $(\mathscr N,d_\mathscr N)$, the parameter $c_\mathscr M(\mathscr N)$ denotes the smallest bi-Lipschitz distortion
with which $\mathscr N$ embeds into $\mathscr M$, i.e., the infimum
of those $D>0$ such that there exists a scaling factor $r>0$ and a
mapping $f:\mathscr N\to \mathscr M$ satisfying $rd_\mathscr
N(x,y)\le d_\mathscr M(x,y)\le Drd_\mathscr N(x,y)$ for all $x,y\in
\mathscr N$ (if no such $f$ exists then set $c_\mathscr M(\mathscr
N)=\infty$).

Bourgain's theorem characterizes superreflexivity of Banach spaces
in terms of their metric structure, but it leaves open the
characterization of $p$-convexity. The notion of $p$-convexity is
crucial for many applications in Banach space theory and metric
geometry, and it turns out that the completion of the Ribe program for
$p$-convexity requires significant additional work beyond Bourgain's
superreflexivity theorem. As a first step in this direction, Lee, Naor and
Peres~\cite{LNP-markov-convex} defined a bi-Lipschitz invariant of metric spaces called \emph{Markov convexity},
which is motivated by Ball's notion of Markov type~\cite{Ball} and
Bourgain's argument in~\cite{Bourgain-trees}.

\begin{definition}[\cite{LNP-markov-convex}]\label{def:markov-convex}

Let $\{X_t\}_{t\in \Z}$ be a Markov chain on a state space
$\Omega$. Given an integer $k\ge 0$,  we denote by $\{\widetilde
X_t(k)\}_{t\in \Z}$ the process which equals $X_t$ for time $t\le
k$, and evolves independently (with respect to the same transition
probabilities) for time $t > k$. Fix $p>0$. A metric space $(X,d_X)$
is called {Markov $p$-convex with constant $\Pi$} if for every
Markov chain $\{X_t\}_{t\in \Z}$ on a state space $\Omega$, and
every $f : \Omega \to X$,
\begin{equation}\label{eq:def-mconvex}
\sum_{k=0}^{\infty}\sum_{t\in \Z}\frac{\E\left[
d_X\left(f(X_t),f\left(\widetilde
X_t\left(t-2^{k}\right)\right)\right)^p\right]}{2^{kp}}
\le \Pi^p
\cdot \sum_{t\in \Z}\E \big[d_X(f(X_t),f(X_{t-1}))^p\big].
\end{equation}
The least constant $\Pi$ for which~\eqref{eq:def-mconvex} holds for all Markov chains is called the Markov
$p$-convexity constant of $X$, and is denoted $\Pi_p(X)$. We shall
say that $(X,d_X)$ is Markov $p$-convex if $\Pi_p(X) < \infty$.
\end{definition}

%\begin{comment}
To gain intuition for Definition~\ref{def:markov-convex},
consider the standard downward random walk starting from the root of the
 binary tree $B_n$ (with absorbing states at the leaves). For an arbitrary mapping $f$ from $B_n$ to a metric space $(X,d_X)$, the triangle inequality implies that for each $k\in \N$ we have
\begin{equation}\label{eq:explain}
\sum_{t\in \Z}\frac{\E\left[
d_X\left(f(X_t),f\left(\widetilde
X_t\left(t-2^{k}\right)\right)\right)^p\right]}{2^{kp}}\lesssim_p \sum_{t\in \Z}\E \big[d_X(f(X_t),f(X_{t-1}))^p\big],
\end{equation}
with asymptotic equality (up to constants depending only on $p$) for $k\le \frac{\log n}{2}$ when $X=B_n$ and $f$ is the identity mapping. On the other hand, if $X$ is a Markov
$p$-convex space then the sum over $k$ of the
left-hand side of~\eqref{eq:explain} is \emph{uniformly bounded}
by the right-hand side of~\eqref{eq:explain}, and therefore
Markov $p$-convex spaces cannot contain $B_n$ with distortion uniformly bounded in $n$.
%\end{comment}

We refer to~\cite{LNP-markov-convex} for more information on the
notion of Markov $p$-convexity. In particular, it is shown
in~\cite{LNP-markov-convex} that the Markov $2$-convexity constant
of an arbitrary weighted tree $T$ is, up to constant factors, the
Euclidean distortion of $T$. We refer to~\cite{LNP-markov-convex}
for $L_p$ versions of this statement and their algorithmic
applications. It was also shown in~\cite{LNP-markov-convex},
via a modification of an argument of Bourgain~\cite{Bourgain-trees},
that if a Banach space $X$ is $p$-convex then it is also Markov
$p$-convex. It was asked in~\cite{LNP-markov-convex} if the converse
is also true. Here we answer this question positively:
\begin{theorem} \label{thm:convexity-coincide}
 A Banach space is $p$-convex  if and only if
it is Markov $p$-convex.
\end{theorem}
Thus Markov $p$-convexity is equivalent to $p$-convexity in Banach
spaces, completing the Ribe program in this case. Our proof of
Theorem~\ref{thm:convexity-coincide} is based on a renorming method
of Pisier~\cite{Pisier-martingales}. It can be viewed as a
nonlinear variant of Pisier's argument, and  several
subtle changes are required in order to adapt it to a
nonlinear condition such as~\eqref{eq:def-mconvex}.

Results similar to Theorem~\ref{thm:convexity-coincide} have been
obtained for the notions of type and cotype of Banach spaces
(see~\cite{BMW,Pisier-type,Ball,NS02,MN-cotype-full,MN-type}), and
have been used to transfer some of the linear theory to the setting
of general metric spaces. This led to several applications to
problems in metric geometry. Apart from the applications of Markov
$p$-convexity that were obtained in~\cite{LNP-markov-convex}, here
we show that this invariant is preserved under Lipschitz quotients.
The notion of Lipschitz quotient was introduced by Gromov~\cite[Sec.~1.25]{Gro07}. Given two metric spaces $(X,d_X)$ and
$(Y,d_Y)$, a surjective mapping $f:X\to Y$ is called a Lipschitz
quotient if it is Lipschitz, and it is also ``Lipschitzly open" in
the sense that there exists a constant $c>0$ such that for every
$x\in X$ and $r>0$,
\begin{equation}\label{eq:def LipQ}
f\left(B_X(x,r)\right)\supseteq B_Y\left(f(x),\frac{r}{c}\right).
\end{equation}
Here we show the following result:
\begin{theorem}\label{thm:invariant}
If $(X,d_X)$ is Markov $p$-convex and $(Y,d_Y)$ is
a Lipschitz quotient of $X$, then $Y$ is also Markov $p$-convex.
\end{theorem}
In~\cite{BJLPS} Bates, Johnson, Lindenstrauss, Preiss and Schechtman
investigated in detail Lipschitz quotients of Banach spaces. Their
results imply that if $2\le p<q$ then $L_q$ is not a Lipschitz
quotient of $L_p$. Since $L_p$ is $p$-convex, it is also Markov
$p$-convex. Hence also all of its subsets are Markov $p$ convex.
But, $L_q$ is not $p$-convex, so we deduce that $L_q$ is not a
Lipschitz quotient of any subset of $L_p$. Thus our new ``invariant
approach" to the above result of~\cite{BJLPS} significantly extends it.
Note that the method of~\cite{BJLPS} is based on a differentiation
argument, and hence it crucially relies on the fact that the
Lipschitz quotient mapping is defined on all of $L_p$ and not just
on an arbitrary subset of $L_p$.

In light of Theorem~\ref{thm:convexity-coincide} it is natural to
ask if Bourgain's characterization of superreflexivity holds
for general metric spaces. Namely, is it true that for any metric
space $X$, if $\lim_{n\to \infty}c_X(B_n)=\infty$ then $X$ is Markov
$p$-convex for some $p<\infty$? This question was asked
in~\cite{LNP-markov-convex}. Here we show that the answer is
negative:
\begin{theorem}\label{thm:doubling}
There exists a metric space $(X,d_X)$ which is not Markov $p$-convex for any $p\in (0,\infty)$, yet $\lim_{n\to\infty} c_X(B_n)=\infty$. In fact, $(X,d_X)$ can be a doubling metric space, and hence $c_X(B_n)\ge 2^{\kappa n}$ for some constant $\kappa>0$.
\end{theorem}
%the Laakso graphs (defined in Section~\ref{sec:laakso})
%are not Markov $p$-convex for any $p<\infty$, even though they do
%not contain $B_n$ with distortion uniformly bounded in $n$ (the last
%assertion follows immediately from the fact that they are doubling
%with constant 16---see~\cite{LP01}).
Theorem~\ref{thm:doubling} is in sharp contrast to the previously established metric characterizations of the linear notions
of type and cotype. Specifically, it was shown by Bourgain, Milman and Wolfson~\cite{BMW} that any metric space with no nontrivial metric type
must contain the Hamming cubes $(\{0,1\}^n,\|\cdot\|_1)$ with
distortion independent of $n$. An analogous result was obtained
in~\cite{MN-cotype-full} for metric spaces with no nontrivial metric
cotype, with the Hamming cube replaced by the $\ell_\infty$ grid $(\{1,\ldots,m\}^n,\|\cdot\|_\infty$).

Our proof of Theorem~\ref{thm:doubling} is based on an analysis of the behavior of a certain Markov chain on the Laakso graphs: a sequence of combinatorial graphs whose definition is recalled in Section~\ref{sec:laakso}. As a consequence of this analysis, we obtain the following distortion lower bound:
\begin{theorem}\label{thm:laaksodist}
For any $p\ge 2$, the Laakso
graph of cardinality $n$ incurs distortion $\Omega\bigl((\log n)^{1/p}\bigr)$ in any
embedding into a $p$-convex Banach space.
\end{theorem}
Thus, in particular, for $p>2$ the $n$-point Laakso
graph incurs distortion $\Omega\bigl((\log n)^{1/p}\bigr)$ in any
embedding into $L_p$. The case of $L_p$ embeddings of the Laakso
graphs when $1<p\le 2$ was
already solved in~\cite{NewmanR02,Laakso,LN-l1-reduction,LMN} using
the uniform 2-convexity property of $L_p$. But, these proofs rely crucially on $2$-convexity and do not
extend to the case of $p$-convexity when $p>2$. Subsequent to the publication of our proof of Theorem~\ref{thm:laaksodist} in the announcement~\cite{charlie-soda}, an alternative proof of this fact was recently discovered by Johnson and Schechtman in~\cite{JS09}.

\subsection{The nonexistence of a metric dichotomy for trees}

 Bourgain's metrical characterization of superreflexivity yields the following statement:

 %~\footnote{Since Theorem~\ref{thm:bourgain dichotomy} is not stated explicitly in~\cite{Bourgain-trees}, we briefly %recall how it follows from Bourgain's argument in Section~\ref{sec:convexity-coincide}}

\begin{theorem}[Bourgain's tree dichotomy~\cite{Bourgain-trees}]\label{thm:bourgain dichotomy} For any Banach space $(X,\|\cdot\|_X)$ one of the following two dichotomic possibilities must hold true:
\begin{itemize}
\item either for all $n\in \N$ we have $c_X(B_n)=1$,
\item or there exists $\alpha=\alpha_X>0$ such that for all $n\in \N$ we have $c_X(B_n)\ge (\log n)^\alpha$.
\end{itemize}
\end{theorem}
Thus, there is a gap in the possible rates of growth of the sequence $\{c_X(B_n)\}_{n=1}^\infty$ when $X$ is a Banach space; consequently, if we were told that, say, $c_X(B_n)=O(\log\log n)$, then we would immediately deduce that actually $c_X(B_n)=1$ for all $n$.  Additional gap  results of this type are known due to the theory of nonlinear type and cotype:

\begin{theorem}[Bourgain-Milman-Wolfson cube dichotomy~\cite{BMW}]\label{thm:BMW dich}
For any metric space $(X,d_X)$ one of the following two dichotomic
possibilities must hold true:
\begin{itemize}
\item either for all $n\in \N$ we have $c_X\left(\{0,1\}^n,\|\cdot\|_1\right)=1$,
\item or there exists $\alpha=\alpha_X>0$ such that for all $n\in \N$ we have $c_X\left(\{0,1\}^n,\|\cdot\|_1\right)\ge n^\alpha$.
\end{itemize}
\end{theorem}
Theorem~\ref{thm:BMW dich} is a metric analogue of Pisier's characterization~\cite{Pis73} of Banach spaces with trivial Rademacher type. A metric analogue of the Maurey-Pisier characterization~\cite{MP76} of Banach spaces with finite Rademacher cotype yields the following dichotomy result for $\ell_\infty$ grids:

\begin{theorem}[Grid dichotomy~\cite{MN-cotype-full}]\label{thm:cotype dich}
For any metric space $(X,d_X)$ one of the following two dichotomic
possibilities must hold true:
\begin{itemize}
\item either for all $n\in \N$ we have $c_X\left(\{0,\dots,n\}^n,\|\cdot\|_\infty\right)=1$,
\item or there exists $\alpha=\alpha_X>0$ such that for all $n\in \N$ we have $c_X\left(\{0,\ldots,n\}^n,\|\cdot\|_\infty\right)\ge n^\alpha$.
\end{itemize}
\end{theorem}
We refer to the survey article~\cite{Men09} for more information on the theory of metric dichotomies.

Note that Theorem~\ref{thm:bourgain dichotomy} is stated for Banach
spaces, while Theorem~\ref{thm:BMW dich} and Theorem~\ref{thm:cotype
dich} hold for general metric spaces. One might expect that as in
the case of previous progress on Ribe's program, a metric theory of
$p$-convexity would result in a proof that Theorem~\ref{thm:bourgain
dichotomy} holds when $X$ is a general metric space. Surprisingly,
we show here that this is not true:

\begin{theorem} \label{lem:Bn->X} There exists a universal constant $C>0$ with the following property.
Assume that  $\{s(n)\}_{n=0}^\infty\subseteq [4,\infty)$ is a
nondecreasing sequence such that $\{n/s(n)\}_{n=0}^\infty$ is also
nondecreasing. Then there exists a metric space $(X,d_X)$
satisfying for all $n\ge 2$,
\begin{equation}\label{eq:sharp c_X}
s\left(\left\lfloor \frac{n}{40 s(n)}\right\rfloor
\right)\left(1-\frac{Cs(n)\log s(n)} {\log n}\right)\leq c_{X}(B_n)
\leq s(n).
\end{equation}
Thus, assuming that $s(n)=o\left(\frac{\log n}{\log\log n}\right)$, there exists a subsequence $\{n_k\}_{k=1}^\infty$ for which
\begin{equation}\label{eq:subsequence}
(1-o(1))s(n_k)\le c_X(B_{n_k})\le s(n_k).
\end{equation}
\end{theorem}
Theorem~\ref{lem:Bn->X} shows that unlike the case of Banach spaces, for general metric spaces, $c_X(B_n)$ can have an arbitrarily slow growth rate.

%The existing literature on nonlinear type and cotype implies that Hamming cubes and $\ell_\infty$ grids have the following local rigidity property:

Bourgain, Milman and Wolfson also obtained in~\cite{BMW} the following finitary version of Theorem~\ref{thm:BMW dich}:
\begin{theorem}[Local rigidity of Hamming cubes~\cite{BMW}]\label{thm:BMW nonquant} For every $\e>0$, $D>1$ and $n\in \N$ there exists $m=m(\e,D,n)\in \N$ such that
$$\lim_{n\to\infty} m(\e,D,n)=\infty,$$ and for every metric $d$ on $\{0,1\}^n$ which is bi-Lipschitz with distortion $\le D$ to the $\ell_1$ (Hamming) metric,
$$
c_{(\{0,1\}^n,d)}\left(\{0,1\}^m,\|\cdot\|_1\right)\le 1+\e.
$$
\end{theorem}
We refer to~\cite{BMW} (see also~\cite{Pisier-type}) for bounds on $m(\e,D,n)$. Informally, Theorem~\ref{thm:BMW nonquant} says that the Hamming cube $(\{0,1\}^n,\|\cdot\|_1)$ is {\em locally rigid} in the following sense: it is impossible to distort the Hamming metric on a sufficiently large hypercube without the resulting metric space containing a hardly distorted  copy of an arbitrarily large Hamming cube. Stated in this way, Theorem~\ref{thm:BMW nonquant} is a metric version of James' theorem~\cite{Jam64} that $\ell_1$ is not a distortable space. The analogue of
Theorem~\ref{thm:BMW nonquant} with the Hamming cube replaced by the $\ell_\infty$ grid $\left(\{0,\ldots,n\}^n,\|\cdot\|_\infty\right)$ is Matou\v{s}ek's BD-Ramsey theorem~\cite{Mat-BD}; see~\cite{MN-cotype-full} for quantitative results of this type in the $\ell_\infty$ case.  The following variant of Theorem~\ref{lem:Bn->X} shows that a local rigidity statement as above fails to hold true for binary trees; it can also be viewed as a negative solution of the distortion problem for the infinite binary tree (see~\cite{OS94} and~\cite[Ch. 13,\,14]{BL} for more information on the distortion problem for Banach spaces).

%(see~\cite{OS94} for the Odell-Schlumprecht solution of the distortion problem for $\ell_2$, and for more background %on this topic).

\begin{theorem}\label{thm:no local rigidity}
Let $B_\infty$ be the complete unweighted infinite binary tree. For
every $D\ge 4$ there exists a metric $d$ on $B_\infty$ that is
$D$-equivalent to the original shortest-path metric on $B_\infty$,
yet for every $\e\in (0,1)$ and $m\in \N$,
$$
c_{(B_\infty,d)}(B_m)\le D-\e\implies m\le D^{CD^2/\e}.
$$
\end{theorem}

The local rigidity problem for binary trees was studied by several
mathematicians. In particular, C. Fefferman asked (private
communication, 2005) whether $\{B_n\}_{n=1}^\infty$ have the local
rigidity property, and Theorem~\ref{thm:no local rigidity} answers
this question negatively. Fefferman also proved a partial local
rigidity result which is a non-quantitative variant of
Theorem~\ref{lem:dich-vertical-Bn} below (see also
Section~\ref{sec:Bn-dichotomy}). We are very grateful to C.
Fefferman for asking us the question that led to the
counter-examples of Theorem~\ref{lem:Bn->X} and Theorem~\ref{thm:no
local rigidity}, for sharing with us his partial positive results,
and for encouraging us to work on these questions. M. Gromov also
investigated the local rigidity problem for binary trees, and proved
(via different methods) non-quantitative partial positive results in
the spirit of Theorem~\ref{lem:dich-vertical-Bn}. We thank M. Gromov
for sharing with us his unpublished work on this topic.

%Theorem~\ref{lem:dich-vertical-Bn}  is a quantitative dichotomy
%theorem for vertically faithful embeddings of binary trees (this
%will be explained momentarily, in the overview of the proofs of
%Theorem~\ref{lem:Bn->X} and Theorem~\ref{thm:no local rigidity} that
%is in contained Section~\ref{sec:overview}). Our proof of this
%result is based on a beautiful argument of J. Matou\v{s}ek, who
%proved a special case of it (containing all the essential ideas)
%in~\cite{Mat-trees}.

The results of Theorem~\ref{lem:Bn->X} and Theorem~\ref{thm:no local
rigidity} are quite unexpected. Unfortunately, their proofs are
delicate and lengthy, and as such constitute the most involved part
of this article. In order to facilitate the understanding of these
constructions, we end the introduction with an overview of the main
geometric ideas that are used in their proofs. This is done in
Section~\ref{sec:overview} below---we recommend reading this section
first before delving into the technical details presented in
Section~\ref{sec:charlie-fails}.

\subsubsection{Overview of the proofs of Theorem~\ref{lem:Bn->X} and
Theorem~\ref{thm:no local rigidity}}\label{sec:overview}

\begin{figure}[ht]
\bigskip
\ \centering
\includegraphics[scale=0.6]{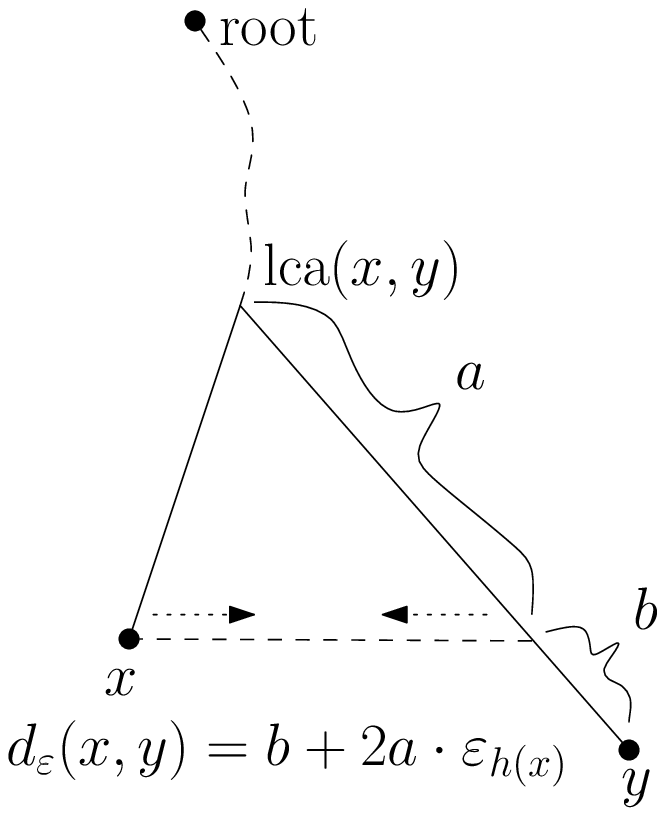}
\caption{{\em The metric $d_\e$ defined on $B_\infty$.  The arrows
indicate horizontal contraction by $\e_{h(x)}$.}}\label{fig:d eps}
\end{figure}

For $x\in B_\infty$ let $h(x)$
 be its depth, i.e., its distance from the root. Also, for $x,y\in
 B_\infty$ let $\lca(x,y)$ denote their least common ancestor.
 The tree metric on $B_\infty$ is then given by:
 $$
d_{B_\infty} (x,y) = h(x)+h(y)-2h(\lca(x,y)).
 $$
The metric space $X$ of Theorem~\ref{lem:Bn->X}  will be $B_\infty$
as a set, with a new metric defined as follows. Given a sequence
$\e=\{\e_n\}_{n=0}^\infty\subseteq (0,1]$ we define
$d_\e:B_\infty\times B_\infty\to [0,\infty)$ by
\[ d_\e(x,y)= |h(y)-h(x)|+2\e_{\min\{h(x),h(y)\}}\cdot \left[\min\{h(x),h(y)\}-h(\lca(x,y))\right].
\]
$d_\e$ does not necessarily satisfy the triangle inequality, but
under some simple conditions on the sequence $\{\e_n\}_{n=0}^\infty$
it does become a metric on $B_\infty$; see
Lemma~\ref{lem:H-tree-metric}. A pictorial description of the metric
$d_\e$ is contained in Figure~\ref{fig:d eps}. Note that when
$\e_n=1$ for all $n$, we have $d_\e=d_{B_\infty}$. Below we call the
metric spaces $(B_\infty,d_\e)$ horizontally distorted trees, or
$H$-trees, in short.

The metric space $(X,d_X)$ of Theorem~\ref{lem:Bn->X} will be
$(B_\infty,d_\e)$, where $\e_n=1/s(n)$ for all $n$. The identity
mapping of $B_n$ into the top $n$-levels of $B_\infty$ has
distortion at most $s(n)$, and therefore $c_X(B_n)\le s(n)$. The
challenge is to prove the lower bound on $c_X(B_n)$
in~\eqref{eq:sharp c_X}. Our initial approach to lower-bounding
$c_X(B_n)$ was Matou\v{s}ek's metric differentiation
proof~\cite{Mat-trees} of asymptotically sharp distortion lower
bounds for embeddings of $B_n$ into uniformly convex Banach spaces.

Following Matou\v{s}ek's terminology~\cite{Mat-trees}, for
$\delta>0$ a quadruple of points $(x,y,z,w)$ in a metric space
$(X,d_X)$ is called a $\delta$-fork if $y\in \Mid(x,z,\delta)\cap
\Mid(x,w,\delta)$, where for $a,b\in X$ the set of
$\delta$-approximate midpoints $\Mid(a,b,\delta)\subseteq X$ is
defined as the set of all $w\in X$ satisfying
$\max\{d_X(x,y),d_X(y,z)\}\le \frac{1+\delta}{2}\cdot d_X(x,z)$. The
points $z,w$ will be called below the prongs of the $\delta$-fork
$(x,y,z,w)$. Matou\v{s}ek starts with the observation that if $X$ is
a uniformly convex Banach space then in any $\delta$-fork in $X$ the
distance between the prongs must be much smaller (as $\delta\to 0$)
than $d_X(x,y)$. Matou\v{s}ek then shows that for all $D>0$, any
distortion $D$ embedding of $B_n$ into $X$ must map some $0$-fork in
$B_n$ to a $\delta$-fork in $X$, provided $n$ is large enough (as a
function of $D$ and $\delta$). This reasoning immediately implies
that $c_X(B_n)$ must be large when $X$ is a uniformly convex Banach
space, and a clever argument of Matou\v{s}ek in~\cite{Mat-trees}
turns this qualitative argument into sharp quantitative bounds.

Of course, we cannot hope to use the above argument of Matou\v{s}ek
in order to prove Theorem~\ref{lem:Bn->X}, since Bourgain's tree
dichotomy theorem (Theorem~\ref{thm:bourgain dichotomy}) does hold
true for Banach spaces. But, perhaps we can mimic this uniform
convexity argument for other target metric spaces? On the face of
it, $H$-trees are ideally suited for this purpose, since the
horizontal contractions that we introduced shrink distances between
the prongs of canonical forks  (call $(x,y,z,w)\in B_\infty$ a
canonical fork if $x$ is an ancestor of $y$ and $z,w$ are
descendants of $y$ at depth $h(x)+2(h(y)-h(x))$). It is for this
reason exactly that we defined $H$-trees.

\begin{figure}[ht]
\begin{center}
\includegraphics[scale=1]{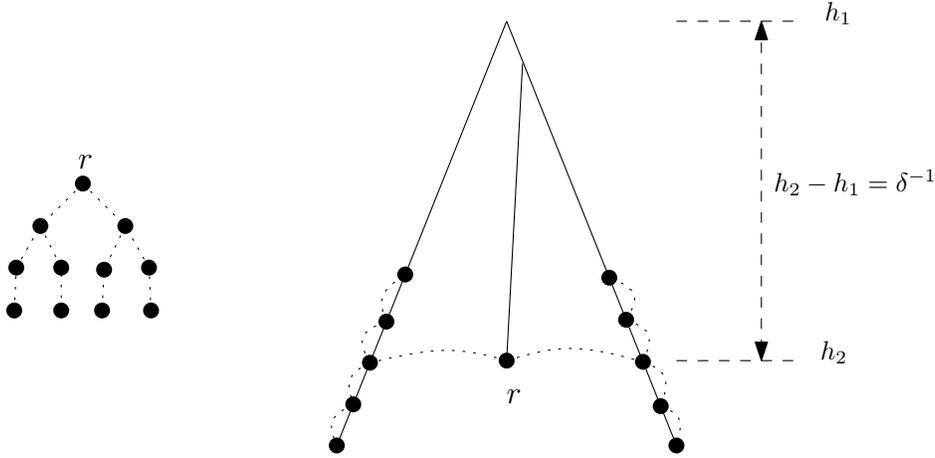}
\label{fig:B_3-4} \caption{ {\em The metric space on the right is
the $H$-tree $(B_\infty,d_\e)$, where $\e_n=\delta$ for all $n$. The
picture describes an embedding of the tree on the left ($B_3$ minus
$4$ leaves) into $(B_\infty,d_\e)$ with distortion at most $6$, yet
all ancestor/descendant distances are distorted by at most
$1+O(\delta)$.}}\label{fig:b3 minus}
\end{center}
\end{figure}

Unfortunately, the situation isn't so simple. It turns out that
$H$-trees do not behave like uniformly convex Banach spaces in terms
of the prong-contractions that they impose of $\delta$-forks.
$H$-trees can even contain larger problematic configurations that
have several undistorted $\delta$-forks; such an example is depicted
in Figure~\ref{fig:b3 minus}.

Thus, in order to prove Theorem~\ref{lem:Bn->X} it does not suffice
to use Matou\v{s}ek's argument that a bi-Lipschitz embedding of a
large enough $B_n$ must send some $0$-fork to a $\delta$-fork. But,
it turns out that this argument applies not only to forks, but also
to larger configurations.
\begin{definition}
Let $(T,d_T)$ be a tree with root $r$, and let $(X,d_X)$ be a metric
space. A mapping $f:T \to X$ is called a $D$-vertically faithful
embedding  if there exists a (scaling factor) $\lambda>0$ satisfying
for any  $x,y\in T$ such that $x$ is an ancestor of $y$,
\begin{equation}\label{eq:def vertical}
  \lambda d_T(x,y)  \le d_X(f(x),f(y)) \leq D\lambda  d_T(x,y).
  \end{equation}
\end{definition}
Recall that the distortion of a mapping $\phi:\mathscr M\to \mathscr
N$ between metric spaces $(\mathscr M,d_\mathscr M)$ and $(\mathscr
N,d_\mathscr N)$ is defined as
$$
\dist(\phi)\eqdef \left(\sup_{\substack{x,y\in \mathscr M\\ x\neq
y}}\frac{d_\mathscr N(\phi(x),\phi(y))}{d_\mathscr
M(x,y)}\right)\cdot \left(\sup_{\substack{x,y\in \mathscr M\\ x\neq
y}}\frac{d_\mathscr M(x,y)}{d_\mathscr N(\phi(x),\phi(y))}\right)\in
[1,\infty].
$$
With this terminology, we can state the following crucial result.

\begin{theorem} \label{lem:dich-vertical-Bn} There exists a
universal constant $c>0$ with the following property. Fix an integer
$t\ge 2$, $\delta,\xi\in(0,1)$, and  $D\ge 2$, and assume that $n\in
\N$ satisfies
\begin{equation}\label{eq:lower n}
 n\ge \frac{1}{\xi}D^{c(t\log t)/\delta} .
\end{equation}
Let $(X,d_X)$ be a metric space and $f:B_n \to X$   a $D$-vertically
faithful embedding. Then there exists a mapping $\phi:B_t\to B_n$
with the following properties.
\begin{itemize}
\item If $x,y\in B_t$ are such that $x$ is an ancestor of $y$, then
$\phi(x)$ is an ancestor of $\phi(y)$.
\item $\dist(\phi)\le 1+\xi$.
\item The mapping $f\circ \phi:B_t\to X$ is a $(1+\delta)$-vertically faithful embedding
of $B_t$ in $X$.
\end{itemize}
\end{theorem}

Theorem~\ref{lem:dich-vertical-Bn} is essentially due to
Matou\v{s}ek~\cite{Mat-trees}. {\Matousek} actually proved this
statement only for $t=2$, since this is all that he needed in order
to analyze forks. But, his proof extends in a straightforward way to
any $t\in \N$. Since we will use this assertion with larger $t$, for
the sake of completeness we reprove it, in a somewhat different way,
in Section~\ref{sec:Bn-dichotomy}. Note that
Theorem~\ref{lem:dich-vertical-Bn} says that $\{B_n\}_{n=1}^\infty$
do have a local rigidity property with respect to vertically
faithfully embeddings.

We solve the problem created by the existence of configurations as
those depicted in Figure~\ref{fig:b3 minus} by studying
$(1+\delta)$-vertically faithful embeddings of $B_4$, and arguing
that they must contain a large  contracted pair of points. This
claim, formalized in Lemma~\ref{lem:B5}, is proved in
Sections~\ref{sec:geometry}, \ref{sec:inembed-B4}.

We begin in Section~\ref{sec:midpoints} with studying how the metric
$P_2$ (3-point path) can be approximately embedded in
$(B_\infty,d_\e)$. We find that there are essentially only two ways
to embed it in $(B_\infty,d_\e)$, as depicted in
Figure~\ref{fig:midpoints}. We then proceed in
Section~\ref{sec:forks} to study $\delta$-forks in
$(B_\infty,d_\e)$. Since forks are formed by ``stitching" two
approximate $P_2$ metrics along a common edge (the handle), we can
limit the ``search space" using the results of
Section~\ref{sec:midpoints}. We find that there are six possible
types of different approximate forks in $(B_\infty,d_\e)$, only four
of which (depicted in Figure~\ref{fig:types}) do not have  highly
contracted prongs. Complete binary trees, and in particular $B_4$,
are composed of forks stitched together, handle to prong. In order
to study handle-to-prong stitching, we investigate in
Section~\ref{sec:3-path} how the metric $P_3$ ($4$-point path) can
be approximately embedded in $(B_\infty,d_\e)$. This is again done
by studying how two $P_2$ metrics can be stitched together, this
time bottom edge to top edge. We find that there are only three
different approximate configurations of $P_4$ in $(B_\infty,d_\e)$.

Using the machinery described above, we study in
Section~\ref{sec:inembed-B4} how the different types of forks can be
stitched together in embeddings of $B_4$ into $(B_\infty,d_\e)$,
reaching the conclusion that a large contraction is unavoidable, and
thus completing the proof of Lemma~\ref{lem:B5}. The proofs of
Theorem~\ref{lem:Bn->X} and Theorem~\ref{thm:no local rigidity} are
concluded in Section~\ref{sec:no-dich-Bn}.

\begin{comment}
\subsubsection*{Organization of the paper}

Section~\ref{sec:convexity-coincide} contains the proof of
Theorem~\ref{thm:convexity-coincide}. In
Section~\ref{sec:mconvexity-observation} we elaborate on the results
discussed in the introduction concerning the Laakso graphs,
Lipschitz quotients, and weaker notions of quotients. Lastly the
proof of Theorem~\ref{thm:Btree-non-dich} (and
Theorem~\ref{thm:non-dich-trees}) are presented in
Section~\ref{sec:charlie-fails}. This section is very long, and is
further broken down into subsections which are described at the
beginning of that section.
\end{comment}

\section{Markov $p$-convexity and $p$-convexity coincide}
\label{sec:convexity-coincide}

In this section we prove Theorem~\ref{thm:convexity-coincide}, i.e.,
that for  Banach spaces $p$-convexity and Markov $p$-convexity are
the same properties. We first show that $p$-convexity implies Markov
$p$-convexity, and in fact it implies a stronger inequality that is
stated in Proposition~\ref{prop:p-convex->sconvex} below. The
slightly weaker assertion that $p$-convexity implies Markov
$p$-convexity was first proved in~\cite{LNP-markov-convex}, based on
an argument from~\cite{Bourgain-trees}. Our argument here is different and simpler.

It was proved in~\cite{Pisier-martingales} that a Banach space $X$ is $p$-convex
if and only if it admits an equivalent norm $\|\cdot \|$ for which
there exists $K>0$ such that
 for every $a,b\in X$,
\begin{equation} \label{eq:p-convexity}
2\|a\|^p+\frac{2}{K^p}\|b\|^p\le \|a+b\|^p+\|a-b\|^p.
\end{equation}
%See also~\cite{BCL} for a proof that~\eqref{eq:p-convexity} follows from~\eqref{eq:def p-conv} (with a different %constant $K$).

\begin{proposition} \label{prop:p-convex->sconvex}
Let $\{X_t\}_{t\in \Z}$ be random variables taking values in a set $\Omega$. For every $s\in \Z$ let $\left\{\widetilde X_t(s)\right\}_{t\in \Z}$ be random variables taking values in $\Omega$, with the following property:
\begin{equation}\label{eq:pairs}
\forall\  r\le s\le t,\ (X_r,X_t)\mathrm{\  and\  } \left(X_r,\widetilde X_t(s)\right) \mathrm{\  have\  the\  same\  distribution.}
\end{equation}
\begin{comment}
\begin{enumerate}
\item Almost surely for all $t\le s$ we have $\widetilde X_t(s)=X_t$.
\item The variables
$(X_{s+1},X_{s+2},\ldots)$ and $\left(\widetilde X_{s+1}(s),\widetilde
X_{s+2}(s),\ldots\right)$ are identically distributed (but not necessarily
independent).
\end{enumerate}
\end{comment}
 Fix $p\ge 2$ and let $(X,\|\cdot\|)$ be a Banach space whose norm
satisfies~\eqref{eq:p-convexity}. Then for every $f:\Omega\to X$ we have
\begin{equation}\label{eq:stochastic}
\sum_{k=0}^{\infty}\sum_{t\in
\Z}\frac{\E\left[\left\|f(X_t)-f\left(\widetilde
X_t(t-2^{k})\right)\right\|^p\right]}{2^{kp}} \\
\le (4K)^p\sum_{t\in
\Z}\E \big[\|f(X_t)-f(X_{t-1})\|^p\big].
\end{equation}
\end{proposition}

\begin{remark}
Observe that condition~\eqref{eq:pairs} holds when $\{X_t\}_{t\in \Z}$ is a Markov chain on a state space
$\Omega$, and $\left\{\widetilde X_t(s)\right\}_{t\in \Z}$ is as in Definition~\ref{def:markov-convex}.
\end{remark}

We start by proving a useful inequality that is a simple consequence of~\eqref{eq:p-convexity}.
\begin{lemma}\label{lem:fork} Let $X$ be a Banach space whose norm
satisfies~\eqref{eq:p-convexity}.  Then for every $x,y,z,w\in X$,
\begin{equation} \label{eq:fork}
\frac{\|x-w\|^p+\|x-z\|^p}{2^{p-1}}+\frac{\|z-w\|^p}{4^{p-1}K^p}\le
\|y-w\|^p+\|z-y\|^p+2\|y-x\|^p.
\end{equation}
\end{lemma}
\begin{proof}
 For every $x,y,z,w\in X$, \eqref{eq:p-convexity} implies that
$$
\frac{\|x-w\|^p}{2^{p-1}}
+\frac{2}{K^p}\left\|y-\frac{x+w}{2}\right\|^p\le
\|y-x\|^p+\|y-w\|^p,
$$
and
$$
\frac{\|z-x\|^p}{2^{p-1}}
+\frac{2}{K^p}\left\|y-\frac{z+x}{2}\right\|^p\le
\|z-y\|^p+\|y-x\|^p.
$$
Summing these two inequalities, and applying the convexity of the
map $u\mapsto \|u\|^p$, we see that
\begin{multline*}
\|y-w\|^p+\|z-y\|^p+2\|y-x\|^p\ge
\frac{\|x-w\|^p+\|z-x\|^p}{2^{p-1}}+\frac{4}{K^p}\cdot\frac{\left\|y-\frac{x+w}{2}\right\|^p+\left\|y-\frac{z+x}{2}\right\|^p}{2}\\
\ge
\frac{\|x-w\|^p+\|z-x\|^p}{2^{p-1}}+\frac{4}{K^p}\cdot\left\|\frac{z-w}{4}\right\|^p,
\end{multline*}
implying~\eqref{eq:fork}.
\end{proof}

\begin{comment}
\begin{remark} \label{rem:fork}
{\em To get intuition on how~\eqref{eq:fork} is used below (specifically, we will return to this interpretation of~\eqref{eq:fork} in Section~\ref{sec:charlie-fails}), assume that $(x,y,z)$ and $(x,y,w)$ are both $(1+\delta)$-bi-Lipschitz to a 3 point path $\{0,1,2\}\subseteq \R$  (this is called $\delta$-fork). This means that for some (scaling factor) $\lambda>0$ we have $\lambda\le \|x-y\|,\|y-z\|,\|y-w\|\le (1+\delta)\lambda$ and $2\lambda\le \|x-z\|,\|x-w\|\le (1+\delta)2\lambda$. In this case, as a consequence of~\eqref{eq:fork}, we have
\[ 4\lambda^p+    \frac{\|z-w\|^p}{4^{p-1}K^p}\le 4(1+\delta)^p\lambda^p.\]
Hence $\frac{\|z-w\|}{\lambda} = O(K\delta^{1/p})$. Thus,  the prongs of $\delta$-forks in $p$-convex Banach spaces must be close to each other.}
\end{remark}
\end{comment}

\begin{proof}[Proof of Proposition~\ref{prop:p-convex->sconvex}]
Using Lemma~\ref{lem:fork}
we see that for every $t\in \mathbb Z$ and $k\in \N$,
\begin{multline*} \frac{\|f(X_t)-f(X_{t-2^k})\|^p+\|f(\widetilde
X_t(t-2^{k-1}))-f(
X_{t-2^{k}})\|^p}{2^{p-1}}+\frac{\|f(X_t)-f(\widetilde
X_t(t-2^{k-1}))\|^p}{4^{p-1}K^p}
\\\le
\|f(X_{t-2^{k-1}})-f(X_t)\|^p+\|f(X_{t-2^{k-1}})-f(\widetilde
X_t(t-2^{k-1}))\|^p+2\|f(X_{t-2^{k-1}})-f(X_{t-2^{k}})\|^p.
\end{multline*}
Taking expectation, and using the assumption~\eqref{eq:pairs}, we get
\begin{multline*} \frac{\E\left[\|f(X_t)-f(X_{t-2^k})\|^p\right]}{2^{p-2}}+\frac{\E\left[\|f(X_t)-f(\widetilde
X_t(t-2^{k-1}))\|^p\right]}{4^{p-1}K^p}
\\\le
2\E\left[\|f(X_{t-2^{k-1}})-f(X_t)\|^p\right]+2\E\left[\|f(X_{t-2^{k-1}})-f(X_{t-2^{k}})\|^p\right].
\end{multline*}
Dividing by $2^{(k-1)p+2}$ this becomes
\begin{multline*} \frac{\E\left[\|f(X_t)-f(X_{t-2^k})\|^p\right]}{2^{kp}}+\frac{\E\left[\|f(X_t)-f(\widetilde
X_t(t-2^{k-1}))\|^p\right]}{2^{(k+1)p}K^p}
\\\le
\frac{\E\left[\|f(X_{t-2^{k-1}})-f(X_t)\|^p\right]}{2^{(k-1)p+1}}+\frac{\E\left[\|f(X_{t-2^{k-1}})-f(X_{t-2^{k}})\|^p\right]}{2^{(k-1)p+1}}.
\end{multline*}
Summing this inequality over $k=1,\ldots,m$ and $t\in\mathbb Z$ we
get
\begin{eqnarray}\label{eq:before cancel}
&&\nonumber\!\!\!\!\!\!\!\!\!\!\!\!\!\!\sum_{k=1}^m\sum_{t\in \Z}
\frac{\E\left[\|f(X_t)-f(X_{t-2^k})\|^p\right]}{2^{kp}}+\sum_{k=1}^m\sum_{t\in\Z}\frac{\left[\E\|f(X_t)-f(\widetilde
X_t(t-2^{k-1}))\|^p\right]}{2^{(k+1)p}K^p}
\\\nonumber&\le&
\sum_{k=1}^m\sum_{t\in
\Z}\frac{\E\left[\|f(X_{t-2^{k-1}})-f(X_t)\|^p\right]}{2^{(k-1)p+1}}+
\sum_{k=1}^m\sum_{t\in \Z}\frac{\E\left[\|f(X_{t-2^{k-1}})-f(X_{t-2^{k}})\|^p\right]}{2^{(k-1)p+1}}\\
&=& \sum_{j=0}^{m-1}\sum_{s\in
\Z}\frac{\E\left[\|f(X_{s})-f(X_{s-2^{j}})\|^p\right]}{2^{jp}}.
\end{eqnarray}
\begin{comment}
\begin{eqnarray*}
&&\!\!\!\!\!\!\!\!\!\!\!\!\!\!\sum_{k=1}^m\sum_{t\in \Z}
\frac{\E\left[\|f(X_t)-f(X_{t-2^k})\|^p\right]}{2^{kp}}+\sum_{k=1}^m\sum_{t\in\Z}\frac{\left[\E\|f(X_t)-f(\widetilde
X_t(t-2^{k-1}))\|^p\right]}{2^{(k+1)p}K^p}
\\&\le&
\sum_{k=1}^m\sum_{t\in
\Z}\frac{\E\left[\|f(X_{t-2^{k-1}})-f(X_t)\|^p\right]}{2^{(k-1)p+1}}+
\sum_{k=1}^m\sum_{t\in \Z}\frac{\E\left[\|f(X_{t-2^{k-1}})-f(X_{t-2^{k}})\|^p\right]}{2^{(k-1)p+1}}\\
&=& \frac12\sum_{j=0}^{m-1}\sum_{t\in
\Z}\frac{\E\left[\|f(X_t)-f(X_{t-2^{j}})\|^p\right]}{2^{jp}}+
\frac12\sum_{j=0}^{m-1}\sum_{s\in \Z}\frac{\E\left[\|f(X_{s})-f(X_{s-2^{j}})\|^p\right]}{2^{jp}}\\
&=& \sum_{j=0}^{m-1}\sum_{s\in
\Z}\frac{\E\left[\|f(X_{s})-f(X_{s-2^{j}})\|^p\right]}{2^{jp}}.
\end{eqnarray*}
\end{comment}

It is only of interest to prove~\eqref{eq:stochastic} when $\sum_{t\in
\Z}\E \big[\|f(X_t)-f(X_{t-1})\|^p\big]<\infty$. By the triangle inequality, this implies that for every $k\in \N$ we have $\sum_{t\in
\Z}\E\left[\|f(X_{t})-f(X_{t-2^{k}})\|^p\right]<\infty$. We may therefore cancel terms in~\eqref{eq:before cancel}, arriving at the following inequality:
\begin{multline*}
\sum_{k=1}^m\sum_{t\in\Z}\frac{\E\left[\|f(X_t)-f(\widetilde
X_t(t-2^{k-1}))\|^p\right]}{2^{(k+1)p}K^p}\\\le\sum_{t\in \Z}\E\left[
\|f(X_t)-f(X_{t-1})\|^p\right]-\sum_{t\in \Z}
\frac{\E\left[\|f(X_t)-f(X_{t-2^m})\|^p\right]}{2^{mp}}\le \sum_{t\in \Z}\E\left[
\|f(X_t)-f(X_{t-1})\|^p\right].
\end{multline*}
Equivalently,
\begin{equation*}
\sum_{k=0}^{m-1}\sum_{t\in \Z}\frac{\E\left[\|f(X_t)-f(\widetilde
X_t(t-2^{k}))\|^p\right]}{2^{kp}}\le (4K)^p\sum_{t\in \Z}\E\left[
\|f(X_t)-f(X_{t-1})\|^p\right]. %\qedhere
\end{equation*}
Proposition~\ref{prop:p-convex->sconvex} now follows by letting $m\to\infty$.
\end{proof}

We next prove the more interesting direction of the equivalence of $p$-convexity and Markov $p$-convexity: a Markov
$p$-convex Banach space is also $p$-convex.

\begin{theorem}\label{thm:renorming} Let $(X,\|\cdot\|)$ be a Banach space which is Markov
$p$-convex with constant $\Pi$. Then for every $\e\in (0,1)$ there
exists a norm ${\tb} \cdot {\tb} $ on $X$ such that for all $x,y\in X$,
$$
(1-\e)\|x\|\le {\tb}  x {\tb}  \le \|x\|,
$$
and
$$
\left|\left|\left| \frac{x+y}{2}\right| \right| \right| ^p\le
\frac{{\tb}  x {\tb}  ^p+{\tb} y{\tb} ^p}{2}-\frac{1-(1-\e)^p}{4\Pi^p (p+1)}\cdot
\left|\left|\left|\frac{x-y}{2}\right| \right| \right|^p.
$$
Thus the norm ${\tb} \cdot{\tb} $ satisfies~\eqref{eq:p-convexity} with
constant $K=O\left(\frac{\Pi}{\e^{1/p}}\right)$.
\end{theorem}
\begin{proof} The fact that $X$ is Markov $p$-convex with constant
$\Pi$ implies that for every Markov chain $\{X_t\}_{t\in \Z}$ with
values in $X$, and for every $m\in \N$, we have
\begin{eqnarray}\label{eq:recall convexity}
\sum_{k=0}^m\sum_{t=1}^{2^m}\frac{\E\left[\bigl\|X_t-\widetilde
X_t(t-2^k) \bigr\|^p\right]}{2^{kp}}\le \Pi^p\sum_{t=1}^{2^m}\E\left[
\left\|X_t-X_{t-1}\right\|^p\right].
\end{eqnarray}

For $x\in X$ we shall say that a Markov chain
$\{X_t\}_{t=-\infty}^{2^m}$ is an $m$-admissible representation of
$x$ if $X_t=0$ for $t\le 0$ and $\E\left[ X_t\right] =t x$ for $t\in
\{1,\ldots,2^m\}$. Fix $\e\in (0,1)$, and denote $\eta=1-(1-\e)^p$.
For every $m\in \mathbb N$ define
\begin{equation}\label{eq:def norm}
{\tb} x{\tb} _m=\inf\left\{\left(\frac{1}{2^m}\sum_{t=1}^{2^m}\E\left[
\left\|X_t-X_{t-1}\right\|^p\right]
-\frac{\eta}{\Pi^p}\cdot \frac{1}{2^m}
\sum_{k=0}^m\sum_{t=1}^{2^m}\frac{\E\left[\bigl\|X_t-\widetilde
X_t(t-2^k) \bigr\|^p\right]}{2^{kp}}\right)^{1/p}\right\},
\end{equation}
where the infimum in~\eqref{eq:def norm} is taken over all
$m$-admissible representations of $x$. Observe that  an $m$-admissible
representation of $x$ always exists, since we can define $X_t=0$
for $t\le 0$ and $X_t=tx$ for $t\in \{1,\ldots,2^m\}$. This
example shows that ${\tb} x{\tb} _m\le \|x\|$. On the other hand, if
$\{X_t\}_{t=-\infty}^{2^m}$ is an $m$-admissible representation of
$x$ then

\begin{comment}
\begin{multline}
\sum_{t=1}^{2^m}\E
\left\|X_t-X_{t-1}\right\|^p-\frac{\eta}{\Pi^p}
\sum_{k=0}^m\sum_{t=1}^{2^m}\frac{\E\bigl\|X_t-\widetilde
X_t(t-2^k) \bigr\|^p}{2^{kp}}
\ge
(1-\eta)\sum_{t=1}^{2^m}\E
\left\|X_t-X_{t-1}\right\|^p\\
\ge (1-\e)^p\sum_{t=1}^{2^m} \left\|\E X_t-\E
X_{t-1}\right\|^p
= (1-\e)^p\sum_{t=1}^{2^m}
\|tx-(t-1)x\|^p
=2^m(1-\e)^p\|x\|^p,
\end{multline}
\end{comment}

\begin{multline}
\sum_{t=1}^{2^m}\E\left[
\left\|X_t-X_{t-1}\right\|^p\right]-\frac{\eta}{\Pi^p}
\sum_{k=0}^m\sum_{t=1}^{2^m}\frac{\E\left[\bigl\|X_t-\widetilde
X_t(t-2^k) \bigr\|^p\right]}{2^{kp}}
\stackrel{\eqref{eq:recall convexity}}{\ge}
(1-\eta)\sum_{t=1}^{2^m}\E\left[
\left\|X_t-X_{t-1}\right\|^p\right]\\
\ge (1-\e)^p\sum_{t=1}^{2^m} \left\|\E \left[X_t\right]-\E\left[
X_{t-1}\right]\right\|^p\label{eq:jensen}
= (1-\e)^p\sum_{t=1}^{2^m}
\|tx-(t-1)x\|^p
=2^m(1-\e)^p\|x\|^p,
\end{multline}
where  in the first inequality of~\eqref{eq:jensen} we used the convexity of the
function $z\mapsto \|z\|^p$. In conclusion, we see that for all
$x\in X$,
\begin{eqnarray}\label{eq:equivalence}
(1-\e)\|x\|\le {\tb} x{\tb} _m\le \|x\|.
\end{eqnarray}

Now take $x,y\in X$ and fix $\delta\in (0,1)$. Let
$\{X_t\}_{t=-\infty}^{2^m}$ be an admissible representation on $x$
and $\{Y_t\}_{t=-\infty}^{2^m}$ be an admissible representation of
$y$ which is stochastically independent of $\{X_t\}_{t=-\infty}^{2^m}$,
such that

\begin{equation}\label{eq:almost for x}
\sum_{t=1}^{2^m}\E\left[
\left\|X_t-X_{t-1}\right\|^p\right]-\frac{\eta}{\Pi^p}\sum_{k=0}^m\sum_{t=1}^{2^m}\frac{\E\left[\bigl\|X_t-\widetilde
X_t(t-2^k) \bigr\|^p\right]}{2^{kp}}  \le 2^m( {\tb} x{\tb} _m^p+\delta),
\end{equation}
and
\begin{equation}\label{eq:almost for y}
\sum_{t=1}^{2^m}\E\left[
\left\|Y_t-Y_{t-1}\right\|^p\right]-\frac{\eta}{\Pi^p}\sum_{k=0}^m\sum_{t=1}^{2^m}\frac{\E\left[\bigl\|Y_t-\widetilde
Y_t(t-2^k) \bigr\|^p\right]}{2^{kp}} \\
\le 2^m({\tb} y{\tb} _m^p+\delta).
\end{equation}

Define a Markov chain
$\{Z_t\}_{t=-\infty}^{2^{m+1}}\subseteq X$ as follows. For
$t\le-2^m$ set $Z_t=0$. With probability $\frac12$ let
$(Z_{-2^m+1},Z_{-2^m+2},\ldots,Z_{2^{m+1}})$ equal
$$
\Bigl(\underbrace{0,\ldots,0}_{2^m\
\text{times}},X_1,X_2,\ldots,X_{2^m},X_{2^m}+Y_1,X_{2^m}+Y_2,\ldots,X_{2^m}+Y_{2^m}\Bigr),
$$
and with probability $\frac12$ let
$(Z_{-2^m+1},Z_{-2^m},\ldots,Z_{2^{m+1}})$ equal
$$
\Bigl(\underbrace{0,\ldots,0}_{2^m\ \text{times}}
,Y_1,Y_2,\ldots,Y_{2^m},X_{1}+Y_{2^m},X_{2}+Y_{2^m},\ldots,X_{2^m}+Y_{2^m}\Bigr).
$$

%For brevity, $\{Z_t\}$ was defined somewhat informally. It can be
%realized by a Markov chain on the state space $\Omega\times
%\Omega'\times\{0,1\}\times \mathbb Z$ where $\Omega$ is the state
%space of $\{X_t\}$, $\Omega'$ is the state space of $\{Y_t\}$, such
%that at time $-2^m$ the Markov chain of $\{Z_t\}$ flips a fair coin
%and decides on sub-state '0' or sub-state '1'.

Hence, $Z_t=0$ for $t\le 0$, for $t\in \{1,\ldots,2^m\}$ we have $\E\left[
Z_t\right]=\frac{\E \left[X_t\right]+\E \left[Y_t\right]}{2}=t\cdot\frac{x+y}{2}$, and for $t\in
\{2^m+1,\ldots,2^{m+1}\}$ we have
\begin{equation*}
\E \left[Z_t\right]=
\frac{\E\left[X_{2^m}+Y_{t-2^m}\right]+\E\left[X_{t-2^m}+Y_{2^m}\right]}{2}
= \frac{2^mx+(t-2^m)y+(t-2^m)x+2^my}{2}=t\cdot\frac{x+y}{2}.
\end{equation*}

Thus $\{Z_t\}_{t=-\infty}^{2^{m+1}}$ is an $(m+1)$-admissible
representation of $\frac{x+y}{2}$. The definition~\eqref{eq:def
norm} implies that

\begin{equation}\label{eq:bound average}
2^{m+1} \left|\left|\left|\frac{x+y}{2}\right|\right|\right|_{m+1}^p
\le \sum_{t=1}^{2^{m+1}}\E\left[
\left\|Z_t-Z_{t-1}\right\|^p\right]
-\frac{\eta}{\Pi^p}
\sum_{k=0}^{m+1}\sum_{t=1}^{2^{m+1}}\frac{\E\left[\bigl\|Z_t -
\widetilde
Z_t(t-2^k) \bigr\|^p\right]}{2^{kp}}.
\end{equation}
Note that by definition,
\begin{equation}\label{eq:first term}
\sum_{t=1}^{2^{m+1}}\E\left[
\left\|Z_t-Z_{t-1}\right\|^p\right]
=\sum_{t=1}^{2^m}\E\left[
\left\|X_t-X_{t-1}\right\|^p\right]+\sum_{t=1}^{2^m}\E\left[
\left\|Y_t-Y_{t-1}\right\|^p\right].
\end{equation}
Moreover,
\begin{multline}\label{eq:split}
\sum_{k=0}^{m+1}\sum_{t=1}^{2^{m+1}}\frac{\E\left[\bigl\|Z_t-\widetilde
Z_t(t-2^k)
\bigr\|^p\right]}{2^{kp}}\\
=\frac{1}{2^{(m+1)p}}\sum_{t=1}^{2^{m+1}}\E\left[\bigl\|Z_t-\widetilde
Z_t(t-2^{m+1}) \bigr\|^p\right]
+
\sum_{k=0}^{m}\sum_{t=1}^{2^{m+1}}\frac{\E\left[\bigl\|Z_t-\widetilde
Z_t(t-2^k) \bigr\|^p\right]}{2^{kp}}.
\end{multline}
We bound each of the terms in~\eqref{eq:split} separately. Note
that by construction  we have for every $t\in \{1,\ldots,2^m\}$,
\begin{equation*}
Z_t-\widetilde Z_t\left(t-2^{m+1}\right)=Z_t-\widetilde
Z_t\left(1-2^{m+1}\right)
=\begin{cases} X_t-Y_t &
\mathrm{with\   probability}\ 1/4,\\
Y_t-X_t &
\mathrm{with\   probability}\ 1/4,\\
X_t-\widetilde{X}_t(1) & \mathrm{with\   probability}\ 1/4,\\
Y_t-\widetilde{Y}_t(1) & \mathrm{with\   probability}\ 1/4.\\
\end{cases}
\end{equation*}

Thus, the first term in the right hand side of~\eqref{eq:split} can be bounded from below
as follows:
\begin{multline}\label{eq:second term1}
\frac{1}{2^{(m+1)p}}\sum_{t=1}^{2^{m+1}}\E\left[\bigl\|Z_t-\widetilde
Z_t(t-2^{m+1}) \bigr\|^p\right] \ge
\frac{1}{2^{(m+1)p+1}}\sum_{t=1}^{2^{m}}\E\left[\|X_t-Y_t\|^p\right] \\
\ge \frac{1}{2^{(m+1)p+1}}\sum_{t=1}^{2^{m}}\|\E\left[
X_t\right]-\E \left[Y_t\right]\|^p
=\frac{\|x-y\|^p}{2^{(m+1)p+1}}\sum_{t=1}^{2^{m}}t^p\ge \frac{2^m\|x-y\|^p}{2^{p+1}(p+1)}.
\end{multline}

We now proceed to bound from below the second term in the right hand side of~\eqref{eq:split}. Note first that for every $k\in \{0,\ldots,m\}$ and every $t\in \{2^m+1,\ldots,2^{m+1}\}$ we have
\begin{equation*}
Z_t-\widetilde Z_t(t-2^k) =\left\{\begin{array}{ll} \left(X_{2^m}+Y_{t-2^m}\right)-\left(\widetilde X_{2^m}(t-2^k)+\widetilde Y_{t-2^m}(t-2^m-2^k)\right)&\mathrm{with\ probability\ 1/2},\\\left(Y_{2^m}+X_{t-2^m}\right)-\left(\widetilde Y_{2^m}(t-2^k)+\widetilde X_{t-2^m}(t-2^m-2^k)\right)& \mathrm{with\ probability\ 1/2}.\end{array}\right.
\end{equation*}
By Jensen's inequality, if $U,V$ are $X$-valued independent random variables with $\E[V]=0$, then $\E\left[\|U+V\|^p\right]\ge \E\left[\|U+\E[V]\|^p\right]=\E\left[\|U\|^p\right]$. Thus, since $\{X_t\}_{t=-\infty}^{2^m}$ and $\{Y_t\}_{t=-\infty}^{2^m}$ are independent,
%Using the independence of $\{X_t\}_{t=-\infty}^{2^m}$ and $\{Y_t\}_{t=-\infty}^{2^m}$, Jensen's inequality implies that
\begin{multline*}
\E\left[\bigl\|Y_{t-2^m}-\widetilde Y_{t-2^m}(t-2^m-2^k) +X_{2^m} - \widetilde X_{2^m}(t-2^k) \bigr\|^p\right]
\\\ge \E\left[\bigl\|Y_{t-2^m}-\widetilde Y_{t-2^m}(t-2^m-2^k) \bigr\|^p\right],
\end{multline*}
and
\begin{multline*}
\E\left[\bigl\|X_{t-2^m}-\widetilde X_{t-2^m}(t-2^m-2^k) +Y_{2^m} - \widetilde Y_{2^m}(t-2^k) \bigr\|^p\right]
\\\ge \E\left[\bigl\|X_{t-2^m}-\widetilde X_{t-2^m}(t-2^m-2^k) \bigr\|^p\right].
\end{multline*}
It follows that for every $k\in \{0,\ldots,m\}$ and every $t\in \{2^m+1,\ldots,2^{m+1}\}$ we have
\begin{multline}\label{eq:drop term}
\E\left[\bigl\|Z_t-\widetilde Z_t(t-2^k)\bigr\|^p\right]\\\ge\frac12\E\left[\bigl\|X_{t-2^m}-\widetilde X_{t-2^m}(t-2^m-2^k) \bigr\|^p\right]+\frac12\E\left[\bigl\|Y_{t-2^m}-\widetilde Y_{t-2^m}(t-2^m-2^k) \bigr\|^p\right].
\end{multline}
Hence,
\begin{eqnarray}\label{eq:second-term-jensen-Y}
\nonumber
&&\!\!\!\!\!\!\!\!\!\!\!\!\!\!\!\!\!
\sum_{k=0}^{m}\sum_{t=1}^{2^{m+1}}  \frac{\E\left[\bigl\|Z_t-\widetilde
Z_t(t-2^k) \bigr\|^p\right]}{2^{kp}}
 \stackrel{\eqref{eq:drop term}}{\ge}
\sum_{k=0}^{m}\sum_{t=1}^{2^m}\frac{\frac12\E\left[\bigl\|X_t-\widetilde
X_t(t-2^k) \bigr\|^p\right]+\frac12\E\left[\bigl\|Y_t-\widetilde Y_t(t-2^k)
\bigr\|^p\right]}{2^{kp}}\\ &+& \sum_{k=0}^{m}\sum_{t=2^m+1}^{2^{m+1}}\frac{\frac12\left[\E\bigl\|X_{t-2^m}-\widetilde
X_{t-2^m}(t-2^m-2^k)
\bigr\|^p\right]+\frac12\E\left[\bigl\|Y_{t-2^m}-\widetilde
Y_{t-2^m}(t-2^m-2^k) \bigr\|^p\right]}{2^{kp}} \nonumber\\
& =&
\sum_{k=0}^{m}\sum_{t=1}^{2^m}\frac{\E\left[\bigl\|X_t-\widetilde
X_t(t-2^k) \bigr\|^p\right]}{2^{kp}}+
\sum_{k=0}^{m}\sum_{t=1}^{2^m}\frac{\E\left[\bigl\|Y_t-\widetilde
Y_t(t-2^k) \bigr\|^p\right]}{2^{kp}}.
\label{eq:second term 2}
\end{eqnarray}

Combining~\eqref{eq:almost for x}, \eqref{eq:almost for y},
\eqref{eq:bound average}, \eqref{eq:first term}, \eqref{eq:split},
\eqref{eq:second term1} and~\eqref{eq:second term 2}, and letting
$\delta$ tend to $0$,  we see that
$$
2^{m+1}\left|\left|\left|\frac{x+y}{2}\right|\right|\right|_{m+1}^p\le
2^m{\tb} x{\tb} _m^p+2^m{\tb} y{\tb} _m^p
-\frac{\eta}{\Pi^p}\cdot\frac{2^m\|x-y\|^p}{2^{p+1}(p+1)},
$$
or,
\begin{equation}\label{eq:before limit}
\left|\left|\left|\frac{x+y}{2}\right|\right|\right|_{m+1}^p\le
\frac{{\tb} x{\tb} _m^p+{\tb} y{\tb} _m^p}{2}
-\frac{\eta}{4\Pi^p(p+1)}\cdot\left\|\frac{x-y}{2}\right\|^p.
\end{equation}
Define for $w\in X$,
$$
{\tb} w{\tb} =\limsup_{m\to \infty} {\tb} w{\tb} _m.
$$
Then a combination of~\eqref{eq:equivalence} and~\eqref{eq:before
limit} yields that
$$
(1-\e)\|x\|\le {\tb} x{\tb} \le \|x\|,
$$
and
\begin{multline}\label{eq:triangle}\left|\left|\left|\frac{x+y}{2}\right|\right|\right|^p\le
\frac{{\tb} x{\tb} ^p+{\tb} y{\tb} ^p}{2}-\frac{\eta}{4\Pi^p(p+1)}\cdot\left\|\frac{x-y}{2}\right\|^p\\\le \frac{{\tb} x{\tb} ^p+{\tb} y{\tb} ^p}{2}-\frac{\eta}{4\Pi^p
(p+1)}\cdot \left|\left|\left|\frac{x-y}{2}\right|\right|\right|^p.
\end{multline}
Note that~\eqref{eq:triangle} implies that the set $\{x\in X:\
{\tb} x{\tb} \le 1\}$ is convex, so that ${\tb} \cdot{\tb} $ is a norm on $X$.
This concludes the proof of Theorem~\ref{thm:renorming}.
\end{proof}

%\section{Observations on Markov convexity}
%\label{sec:mconvexity-observation}

\section{A doubling space which is not Markov $p$-convex for any $p\in (0,\infty)$}
\label{sec:laakso}

\begin{comment}
\begin{figure}[ht]
\bigskip
\ \centering
\includegraphics[scale=0.7]{lang}
\caption{The Laakso graphs.} \label{fig:lang}
\end{figure}
\end{comment}

\parpic[r]{\includegraphics[scale=0.7]{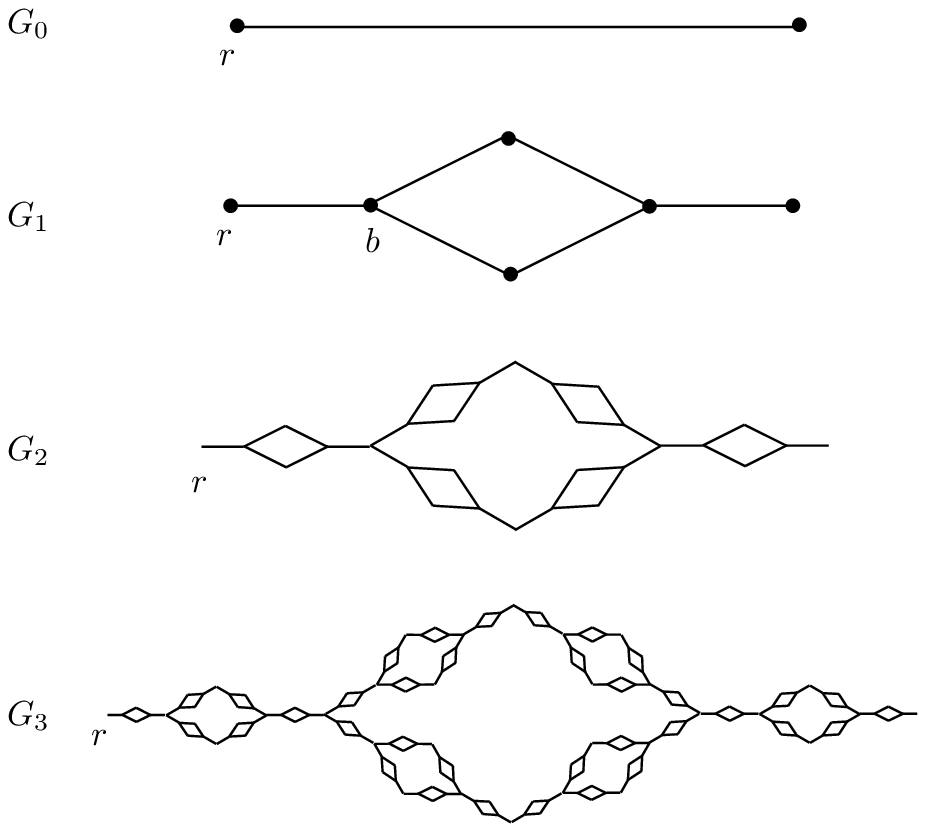}}
Consider the Laakso graphs~\cite{Laakso}, $\{G_i\}_{i=0}^\infty$, which are
defined as follows. $G_0$ is the graph consisting of one
edge of unit length. To construct $G_i$, take six copies of $G_{i-1}$ and scale
their metric by a factor of $\frac14$. We glue four of them
cyclicly by identifying pairs of endpoints, and attach at two
opposite gluing points the remaining two copies. Note that each edge of $G_i$ has length $4^{-i}$; we denoted the resulting shortest path metric on $G_i$ by $d_{G_i}$. As shown in~\cite[Thm.~2.3]{LP01}, the doubling constant of metric space $(G_i,d_{G_i})$ is at most $6$.

% See Figure~\ref{fig:lang}.

We direct $G_m$ as follows. Define the root of $G_m$ to be (an arbitrarily chosen) one of the two vertices having only one adjacent edge.
In the figure %~\ref{fig:lang}
this could be the leftmost vertex $r$. Note that in no edge the two endpoints are at the same distance from the root.
The edges of $G_m$ are then directed from the endpoint closer to the root to the endpoint further away from the root.
The resulting directed graph is acyclic.
We now define $\{X_t\}_{t=0}^{4^m}$ to be the standard random walk on the directed graph $G_m$, starting from the root. This random walk is extended  to $t\in \mathbb Z$
by stipulating that  $X_t=X_0$ for $t<0$, and  $X_t=X_{4^m}$ for $t>4^m$.

\begin{proposition} \label{prop:lb-mconvex-laakso}
For the random walk defined above,
\begin{equation} \label{eq:lang-mconvex}
\sum_{k=0}^{2m}\sum_{t\in \mathbb{Z}}
\frac{\mathbb{E}\left[d_{G_m}\bigl(X_t,\widetilde{X}_{t}(t-2^k)\bigr)^p\right]}{2^{kp}}
 \gtrsim \frac{m}{8^{p}}  \sum_{t\in \mathbb{Z}} \mathbb{E} \left[d_{G_m}(X_t, X_{t-1})^p\right].
\end{equation}
\end{proposition}
\begin{proof}For every $t\in \Z$ we have,
$$\mathbb{E} \left[d_{G_m}(X_t, X_{t-1})^p\right]=\left\{\begin{array}{ll}4^{-mp}& t\in\{0,\ldots 4^m-1\},\\
0& \mathrm{otherwise.}\end{array}\right.$$
Hence,
\begin{equation}\label{eq:RHS laakso}
 \sum_{t\in \mathbb{Z}} \mathbb{E} \left[d_{G_m}(X_t, X_{t-1})^p\right] =4^{-m(p-1)}.
 \end{equation}
%Also note that for $t\in \{0,4^m-a\}$, $d(X_t,X_{t+a})=a/4^m$.

Fix $k\in\{0,\ldots, 2m-2\}$ and write $h= \lceil k/2\rceil$. View $G_m$ as being built from $A=G_{m-h}$, where each edge of $A$ has been replaced by a copy of $G_h$. Note that for every $i\in \{0,\ldots,4^{m-h-1}+1\}$, at time $t=(4i+1)4^h$ the walk $X_t$ is at a vertex of $G_m$ which has two outgoing edges, corresponding to distinct copies of $G_h$. To see this it suffices to show that all vertices of $G_{m-h}$ that are exactly $(4i+1)$ edges away from the root, have out-degree $2$. This fact is true since $G_{m-h}$ is obtained from $G_{m-h-1}$ by replacing each edge by a copy of $G_1$, and each such copy of $G_1$ contributes one vertex of out-degree $2$, corresponding to the vertex labeled $b$ in the figure describing $G_1$.

Consider the set of times
\begin{equation*}
T_k\eqdef\left\{0,\ldots,4^m-1\right\}
\bigcap \left( \bigcup_{i=0}^{4^{m-h-1}+1} \bigl [(4i+1)4^h+4^{h-2}, (4i+1)4^h+2\cdot 4^{h-2}\bigr ] \right) .
\end{equation*}
For $t\in T_k$ find $i\in \{0,\ldots,4^{m-h-1}+1\}$ such that $t\in \bigl [(4i+1)4^h+4^{h-2}, (4i+1)4^h+2\cdot 4^{h-2}\bigr ]$. Since, by the definition of $h$,  we have $t-2^k\in\bigl[(4i+1)4^h-4^h,(4i+1)4^h\bigr)$, the walks $\{X_s\}_{s\in \Z}$ and $\{\widetilde X_s(t-2^k)\}_{s\in \Z}$ started evolving independently at some vertex lying in a copy of $G_h$ preceding a vertex $v$ of $G_m$ which has two outgoing edges, corresponding to distinct copies of $G_h$. Thus, with probability at least $\frac12$, the walks $X_t$ and $\widetilde X_t(t-2^k)$ lie on two distinct copies of $G_h$ in $G_m$, immediately following the vertex $v$, and at distance at least $4^{h-2}\cdot 4^{-m}$ and at most $2\cdot 4^{h-2}\cdot 4^{-m}$ from $v$. Hence, with probability at least $\frac12$ we have $d_{G_m}\left(X_t,\widetilde X_t(t-2^k)\right)\ge 2\cdot 4^{h-2}\cdot 4^{-m}=2^{2h-3-2m}$, and therefore,
$$
\frac{\mathbb{E}\left[d(X_t,\widetilde X_{t}(t-2^k))^p\right]}{2^{kp}} \ge \frac{\frac 12 2^{(2h-3-2m)p}}{2^{kp}}\geq 2^{-(2m+3)p-1}.
$$
We deduce that for all $k\in\{0,\ldots, 2m-2\}$,
\begin{multline}\label{eq:lower laakso sum}
\sum_{t\in \Z} \frac{\mathbb{E}\left[d(X_t,\widetilde X_{t}(t-2^k))^p\right]}{2^{kp}}\ge \sum_{t\in T_k} \frac{\mathbb{E}\left[d(X_t,\widetilde X_{t}(t-2^k))^p\right]}{2^{kp}}\ge |T_k|\cdot 2^{-(2m+3)p-1}\\\gtrsim 4^{h-2}\cdot 4^{m-h-1}\cdot 2^{-(2m+3)p-1}\gtrsim \frac{1}{8^{p}}4^{-m(p-1)}.
\end{multline}
A combination of~\eqref{eq:RHS laakso} and~\eqref{eq:lower laakso sum}  implies~\eqref{eq:lang-mconvex}.
\end{proof}

\begin{proof}[Proof of Theorem~\ref{thm:doubling}] As explained in~\cite{Laakso,LP01}, by passing to an appropriate Gromov-Hausdorff limit, there exists a doubling metric space $(X,d_X)$ that contains an isometric copy of all the Laakso graphs $\{G_m\}_{m=0}^\infty$. Proposition~\ref{prop:lb-mconvex-laakso} therefore implies that $X$ is not Markov $p$-convex for any $p\in (0,\infty)$.
\end{proof}

\begin{proof}[Proof of Theorem~\ref{thm:laaksodist}] Let $(X,d_X)$ be a Markov $p$-convex metric space, i.e, $\Pi_p(X)<\infty$. Assume that $f:G_m\to X$ satisfies
\begin{equation}\label{eq:bilip}
x,y\in G_m\implies \frac{1}{A}d_{G_m}(x,y)\le d_X(f(x),f(y))\le Bd_{G_m}(x,y).
\end{equation}
Let $\{X_t\}_{t\in \Z}$ be the random walk from Proposition~\ref{prop:lb-mconvex-laakso}. Then
\begin{eqnarray*}
\frac{m}{8^{p}A^p}  \sum_{t\in \mathbb{Z}} \mathbb{E} \left[d_{G_m}(X_t, X_{t-1})^p\right]&\stackrel{\eqref{eq:lang-mconvex}}{\lesssim}& \frac{1}{A^p}\sum_{k=0}^{2m}\sum_{t\in \mathbb{Z}}
\frac{\mathbb{E}\left[d_{G_m}\bigl(X_t,\widetilde{X}_{t}(t-2^k)\bigr)^p\right]}{2^{kp}}\\&\stackrel{\eqref{eq:bilip}}{\le}& \sum_{k=0}^{2m}\sum_{t\in \mathbb{Z}}
\frac{\mathbb{E}\left[d_{X}\bigl(f(X_t),f(\widetilde{X}_{t}(t-2^k))\bigr)^p\right]}{2^{kp}}\\&\stackrel{\eqref{eq:def-mconvex}}{\le}& \Pi_p(X)^p
\sum_{t\in \mathbb{Z}} \mathbb{E} \left[d_{X}(f(X_t), f(X_{t-1}))^p\right]\\
&\stackrel{\eqref{eq:bilip}}{\le}& \Pi_p(X)^pB^p \sum_{t\in \mathbb{Z}} \mathbb{E} \left[d_{G_m}(X_t, X_{t-1})^p\right].
\end{eqnarray*}
Thus $AB\gtrsim m^{1/p}\gtrsim (\log |G_m|)^{1/p}$.
\end{proof}

\section{Lipschitz quotients}

Say that a metric space $(Y,d_Y)$ is a $D$-Lipschitz quotient of a metric space $(X,d_X)$ if there exist $a,b>0$ with $ab\le D$ and a mapping $f:X\to Y$ such that for all $x\in X$ and $r>0$,
\begin{equation}\label{eq:def D-Q}
B_Y\left(f(x),\frac{r}{a}\right)\subseteq f\left(B_X(x,r)\right)\subseteq B_Y(f(x),br).
\end{equation}
Observe that the last inclusion in~\eqref{eq:def D-Q} is to equivalent to the fact that $f$ is $b$-Lipschitz.

The following proposition implies Theorem~\ref{thm:invariant}.

\begin{proposition} \label{prop:lipschitz-invariant}
If $(Y,D_Y)$ is a $D$-Lipschitz quotient of $(X,d_X)$ then $\Pi_p(Y) \le D\cdot \Pi_p(X)$.
\end{proposition}
\begin{proof}
Fix $f:X\to Y$ satisfying~\eqref{eq:def D-Q}. Also, fix a Markov chain $\{X_t\}_{t\in \Z}$ on a state space $\Omega$, and
a mapping $g:\Omega\to Y$.

Fix $m\in \Z$ and let $\Omega^*$ be the set of finite sequences of elements of $\Omega$ starting at time $m$, i.e., the set of sequences of the form  $(\omega_i)_{i=m}^t\in \Omega^{t-m+1}$ for all $t\ge m$. It will be convenient to consider the Markov chain  $\{X^*_t\}_{t=m}^\infty$ on $\Omega^*$ which is given by:
$$
\Pr\left[X_t^*=(\omega_m,\omega_{m+1},\ldots,\omega_t)\right]=\Pr\left[X_m=\omega_m,X_{m+1}=\omega_{m+1},\ldots,X_t=\omega_t\right].
$$
%Thus,
%\begin{equation*}
% \Pr \Bigl [X^*_t=(\omega_1,\ldots, \omega_{t-1},\omega_t)\, \big|\, X^*_{t-1}=(\omega_1,\ldots, \omega_{t-1})
%\Bigr] \; = \;
%\Pr\Bigl[X_t=\omega_t\,\big|\, X_{t-1}=\omega_{t-1}\Bigr],
%\end{equation*}
%and the rest of the transition probabilities are 0.
Also, define $g^*:\Omega^* \to Y$ by $g^*(\omega_1,\ldots,\omega_t)=g(\omega_t)$. By definition, $\{g^*(X^*_t)\}_{t=m}^\infty$ and $\{g(X_t)\}_{t=m}^\infty$ are identically distributed.

We next define a mapping $h^*:\Omega^*\to X$ such that $f\circ h^*=g^*$ and for all $(\omega_m,\ldots,\omega_t)\in \Omega^*$,
\begin{equation}\label{eq:lift}
d_X\left(h^*(\omega_m,\ldots,\omega_{t-1}),h^*(\omega_m,\ldots,\omega_{t})\right)\le a d_Y(g(\omega_{t-1}),g(\omega_t)).
\end{equation}
For $\omega^*\in \Omega^*$, we will define $h^*(\omega^*)$ by induction on the length of $\omega^*$.
If $\omega^*=(\omega_m)$, then we fix $h^*(\omega^*)$ to be an arbitrary element in $f^{-1}(g(\omega_m))$. Assume that $\omega^*=(\omega_m,\ldots, \omega_{t-1},\omega_t)$ and that $h^*(\omega_m,\ldots,\omega_{t-1})$ has been defined. Set $x=f(h^*(\omega_m,\ldots,\omega_{t-1}))=g^*(\omega_m,\ldots,\omega_{t-1})=g(\omega_{t-1})$ and $r=ad_Y(g(\omega_{t-1}),g(\omega_t))$. Since $g(\omega_t)\in B_Y\left(x,r/a\right)$, it follows from~\eqref{eq:def D-Q} there exists $y\in X$
such that $f(y)=g(\omega_t)$, and
$d_X(x,y) \le r$. We then define
$h^*((\omega_m,\ldots,\omega_{t-1},\omega_t))\eqdef y$.

Write $X_t^*=X_m^*$ for $t\le m$. By the Markov $p$-convexity of $(X,d_X)$, we have
\begin{equation} \label{eq:X-mconvex}
\sum_{k=0}^{\infty}\sum_{t\in \Z}\frac{\E \left[d_X\bigl(h^*(X^*_t),h^*(\widetilde
X^*_t(t-2^{k}))\bigr)^p\right]}{2^{kp}}
\le \Pi_p(X)^p \sum_{t\in \Z}\E\left[
d_X(h^*(X^*_t),h^*(X^*_{t-1}))^p\right].
\end{equation}
By~\eqref{eq:lift} we have for every $t\ge m+1$,
$$
d_X(h^*(X^*_t),h^*(X^*_{t-1}))\le ad_Y(g(X_t),g(X_{t-1})),
$$
while for $t\le m$ we have $d_X(h^*(X^*_t),h^*(X^*_{t-1}))=0$.
Thus,
\begin{equation}\label{eq:lift RHS}
\sum_{t\in \Z}\E\left[
d_X(h^*(X^*_t),h^*(X^*_{t-1}))^p\right]\le a^p\sum_{t\in \Z}\E\left[d_Y(g(X_t),g(X_{t-1}))^p\right].
\end{equation}
At the same time, using the fact that $f$ is $b$-Lipschitz and $f\circ h^*=g^*$, we see that if $t\ge m+2^k$,
\begin{multline*}
d_X\bigl(h^*(X^*_t),h^*(\widetilde
X^*_t(t-2^{k}))\bigr)\ge \frac{1}{b}d_Y\bigl(f(h^*(X^*_t)),f(h^*(\widetilde
X^*_t(t-2^{k})))\bigr)\\=\frac{1}{b}d_Y\bigl(g^*(X^*_t),g^*(\widetilde
X^*_t(t-2^{k}))\bigr)=\frac{1}{b}d_Y\bigl(g(X_t),g(\widetilde.
X_t(t-2^{k}))\bigr)
\end{multline*}
Thus,
\begin{equation}\label{eq:lift LHS}
\sum_{k=0}^{\infty}\sum_{t\in \Z}\frac{\E \left[d_X\bigl(h^*(X^*_t),h^*(\widetilde
X^*_t(t-2^{k}))\bigr)^p\right]}{2^{kp}}\ge \frac{1}{b^p}\sum_{k=0}^{\infty}\sum_{t=m+2^k}^\infty\frac{\E \left[d_Y\bigl(g(X_t),g(\widetilde
X_t(t-2^{k}))\bigr)^p\right]}{2^{kp}}.
\end{equation}
By combining~\eqref{eq:lift RHS} and~\eqref{eq:lift LHS} with~\eqref{eq:X-mconvex}, and letting $m$ tend to $-\infty$, we get the inequality:
$$
\sum_{k=0}^{\infty}\sum_{t\in \Z}\frac{\E \left[d_Y\bigl(g(X_t),g(\widetilde
X_t(t-2^{k}))\bigr)^p\right]}{2^{kp}}\le \left(ab\Pi_p(X)\right)^p\sum_{t\in \Z}\E\left[d_Y(g(X_t),g(X_{t-1}))^p\right].
$$
Since this inequality holds for every Markov chain $\{X_t\}_{t\in \Z}$ and every $g:\Omega\to Y$, and since $ab\le D$, we have proved that $\Pi_p(Y)\le D\Pi_p(X)$, as required.
\end{proof}

\section{A dichotomy theorem  for vertically faithful embeddings of trees}
\label{sec:Bn-dichotomy}

In this section we prove Theorem~\ref{lem:dich-vertical-Bn}. The
proof  naturally breaks into two parts. The first is the following
BD Ramsey property of paths (which can be found non-quantitatively
in~\cite{Mat-BD}, where also the BD Ramsey terminology is
explained).

A mapping $\phi:\MM\to \NN$ is called a {\em rescaled isometry} if
$\dist(\phi)=1$, or equivalently there exists $\lambda>0$ such that
$d_\NN(\phi(x),\phi(y))=\lambda d_\MM(x,y)$ for all $x,y\in \MM$.
For $n\in \N$ let $P_n$ denote the $n$-path, i.e., the set
$\{0,\ldots,n\}$ equipped with the metric inherited from the real
line.

\begin{proposition}\label{prop:path-boosting}
Fix $\delta\in (0,1)$, $D\ge 2$ and $t,n\in \N$ satisfying $n\ge
D^{(4t\log t)/\delta}$. If $f:P_n\to X$ satisfies $\dist(f)\le D$
then there exists a rescaled isometry $\phi:P_t\to P_n$ such that
$\dist(f\circ \phi)\le 1+\delta$.
\end{proposition}

Given a metric space $(X,d_X)$ and a nonconstant mapping $f:P_n\to X$, define
\begin{equation*}\label{eq:defT}
T(X,f)\eqdef \frac{d_X(f(0),f(n)}{n\max_{i\in \{1,\ldots,n\}}d_X(f(i-1),f(i))}=\frac{d_X(f(0),f(n))}{n\|f\|_{\Lip}}.
\end{equation*}
If $f$ is a constant mapping (equivalently $\max_{i\in \{1,\ldots,n\}}d_X(f(i-1),f(i))=0$) then we set $T(X,f)=0$. Note that by the triangle inequality we always have $T(X,f)\le 1$.

\begin{lemma} \label{lem:sub-mult}
For every $m,n\in\mathbb N$ and $f:P_{mn}\to X$, there exist rescaled isometries  $\phi^{(n)}:P_n\to P_{mn}$ and $\phi^{(m)}:P_m\to P_{mn}$, such that
\[T(X,f) \leq T\left(X,f\circ \phi^{(m)}\right) \cdot T\left(X,f\circ \phi^{(n)}\right) .\]
%In particular, $T_{mn}(X)\le T_m(X) \cdot T_n(X)$.
\end{lemma}
\begin{proof}
Fix $f:P_{mn}\to X$ and define $\phi^{(m)}:P_m\to P_{mn}$ by $\phi^{(m)}(i)=in$. Then,
\begin{equation}\label{eq:phi m}
d_X(f(0),f({mn}))\leq T\left(X,f\circ\phi^{(m)}\right)  m \max_{i\in \{1,\ldots,m\}} d_X(f({(i-1)n}),f({in})).
\end{equation}
Similarly, for every $i\in\{1,\ldots, m\}$ define $\phi_i^{(n)}:P_n\to P_{mn}$ by $\phi^{(n)}_i(j)=(i-1)n+j$. Then
\begin{equation}\label{eq:phi n} d_X(f({(i-1)n}),f({in})) \leq T\left(X,f\circ\phi_i^{(n)}\right) n \max_{j\in\{1,\ldots,n\}}
d_X(f({(i-1)n+j-1}), f({(i-1)n+j})).
\end{equation}
Letting $i\in \{1,\ldots,m\}$ be such that $T\left(X,f\circ\phi_i^{(n)}\right)$ is maximal, and $\phi^{(n)}=\phi^{(n)}_i$, we conclude that
\[ d_X(f(0),f({mn}))\stackrel{\eqref{eq:phi m}\wedge\eqref{eq:phi n}}{\leq} T\left(X,f\circ\phi^{(m)}\right)T\left(X,f\circ\phi^{(n)}\right)mn \max_{i\in \{1,\ldots,mn\}} d_X(f({i-1}),f(i))
.\qedhere\]
% $T(X,f)\leq T(X,f\circ\iota^{(m)}) \cdot T(X,f\circ\iota_i^{(n)})$.
\end{proof}

\begin{lemma} \label{lem:dist-lb}
For every $f:P_m \to X$ we have $\dist(f) \ge 1/ T(X,f)$.
\end{lemma}
\begin{proof} Assuming $a|i-j|\le d_X(f(i),f(j))\le b|i-j|$ for all $i,j\in P_m$, the claim is $bT(X,f)\ge a$. Indeed,
$am\leq d_X(f(0),f(m)) \leq T(X,f) m
\max_{i=\in\{1,\ldots,m\}} d_X(f({i-1}),f(i)) \leq T(X,f) bm$.
\end{proof}

\begin{lemma}\label{lem:ln-embed}
Fix $f:P_m\to X$.
If $0< \e<1/m$ and $T(X,f)\ge1-\e$, then $\dist (f) \le 1/(1-m\e)$.
%In particular, If $T_m(X)>1-\e$, then
%$c_X(P_m)\leq 1+2 m \e$.
%In particular, if $T(X)=1$, then $c_X(P_m)=1$ for every $m\in \mathbb{N}$.
\end{lemma}
\begin{proof}
%Fix $f:P_m\to X$ such that
%\begin{equation}\label{eq:2}
%d(f(0),f(m)) \geq (1-\e)m \max_i d(f({i-1}),f(i)) .
%\end{equation}
Denote  $b=\max_{i\in \{1,\ldots,n\}} d_X(f(i),f({i-1}))>0$. For every
$0\leq i<j\leq m$ we have $ d_X(f(i),f(j))\leq \sum_{\ell=i+1}^j d_X(f({\ell-1}),f(\ell))\leq b|j-i|$, and
\begin{multline*}
(1-\e)mb\le T(X,f)mb=d_X(f(0),f(m))\\\le d_X(f(0),f(i))+d_X(f(i),f(j))+d_X(f(j),f(m))\le d_X(f(i),f(j))+b(m+i-j).
\end{multline*}
Thus $d_X(f(i),f(j))\ge b(j-i-m\e)\ge (1-m\e )b|j-i|$.
\end{proof}

\begin{proof}[Proof of Proposition~\ref{prop:path-boosting}] Set $k=\lfloor \log_t n \rfloor$ and denote by $I$ the identity mapping from $P_{t^k}$ to $P_n$. By Lemma~\ref{lem:dist-lb} we have $T(X,f\circ I)\ge 1/D$. An iterative application of Lemma~\ref{lem:sub-mult} implies that there exists a rescaled isometry $\phi:P_t\to P_{t^k}$ such that
$$
T(X,f\circ I\circ \phi)\ge D^{-1/k}\ge e^{-2\log D/\log_t n}\ge e^{-\delta/(2t)}\ge 1-\frac{\delta}{2t}.
$$
By Lemma~\ref{lem:ln-embed} we therefore have $\dist(f\circ I\circ\phi)\le 1/(1-\delta/2)\le 1+\delta$.
\end{proof}

The second part of the proof of Theorem~\ref{lem:dich-vertical-Bn}
uses the following combinatorial lemma due to
Matou\v{s}ek~\cite{Mat-trees}. Denote by $T_{k,m}$ the complete
rooted tree of height $m$, in which every non-leaf vertex has $k$
children. For a rooted tree $T$, denote by $\mathrm{SP}(T)$ the
set of all unordered pairs $\{x,y\}$ of distinct vertices of $T$
such that $x$ is an ancestor of $y$.

\begin{lemma}[{\cite[Lem.~5]{Mat-trees}}] \label{lem:mat-ramsey}
Let $m,r,k\in \N$ satisfy $k\geq
r^{(m+1)^2}$. Suppose that each of the pairs from
$\mathrm{SP}(T_{k,m})$ is colored by one of $r$ colors. Then there
exists a copy $T'$ of $B_m$ in this $T_{k,m}$ such that the color
of any pair $\{x,y\}\in\mathrm{SP}(T')$ only depends on the levels
of $x$ and $y$.
\end{lemma}

\begin{proof}[Proof of Lemma~\ref{lem:dich-vertical-Bn}]
Let $f:B_n\to X$ be a $D$-vertically faithful embedding, i.e., for
some $\lambda>0$ it satisfies \begin{equation}\label{eq:f vert}
\lambda d_{B_n}(x,y)\le d_X(f(x),f(y))\le D\lambda d_{B_n}(x,y)
\end{equation}
whenever $x,y\in B_n$ are such that $x$ is an ancestor of $y$.

Let $k,\ell\in \N$ be auxiliary parameters to be determined later,
and define $m=\left\lfloor n/(k\ell)\right\rfloor $.  We first
construct a mapping $g: T_{2^k,m} \to B_n$ in a top-down manner as
follows. If $r$ is the root of $T_{2^k,m}$ then $g(r)$ is defined to
be the root of $B_n$. Having defined $g(u)$, let
$v_1,\ldots,v_{2^k}\in T_{2^k,m}$ be the children of $u$, and let
$w_1,\ldots,w_{2^k}\in B_n$ be the descendants of $g(u)$ at depth
$k$ below $g(u)$. For each $i\in \{1,\ldots,2^k\}$ let $g(v_i)$ be
an arbitrary descendant of $w_i$ at depth $h(g(u))+\ell k$. Note
that for this construction to be possible we need to have $m\ell
k\le n$, which is ensured by our choice of $m$.

By construction, if $x,y\in T_{2^k,m}$ and $x$ is an ancestor of
$y$, then $g(x)$ is an ancestor of $g(y)$ and $d_{B_n}(g(x),g(y))=
\ell k d_{T_{2^k,m}}(x,y)$. Also, if $x,y\in T_{2^k,m}$ and
$\lca(x,y)=u$, then we have $h(\lca(g(x),g(y)))\in
\{h(g(u)),h(g(u))+1,\ldots,h(g(u))+k-1\}$. This implies that
\begin{equation*}\label{eq:dist g}
((\ell-1) k+1) d_{T_{2^k,m}}(x,y)\le d_{B_n}(x,y)\le \ell k
d_{T_{2^k,m}}(x,y).
\end{equation*}
Thus, assuming $\ell\ge 2$, we have $\dist(g)\le 1+2/\ell$.
Moreover, denoting $F=f\circ g$ and using~\eqref{eq:f vert}, we see
that if $x,y\in T_{2^k,m}$ are such that $x$ is an ancestor of $y$
then
\begin{equation}\label{eq:F vert}
k\ell\lambda d_{T_{2^k,m}}(x,y)\le d_X(F(x),F(y))\le D\ell k\lambda
d_{T_{2^k,m}}(x,y).
\end{equation}

Color every pair $\{x,y\}\in \mathrm{SP}(T_{2^k,m})$ with the color
$$
\chi(\{x,y\})\eqdef \left\lfloor \log_{1+\delta/4}
\left(\frac{d_X(F(x),F(y))}{k\ell\lambda d_{T_{2^k,m}}(x,y)}\right)
\right\rfloor\in \{1,\ldots, r\},
$$
where $r=\lceil \log_{1+\delta/4}D \rceil$. Assuming that
\begin{equation}\label{eq:coloring assumption}
2^k\ge r^{(m+1)^2}, \end{equation}
 by  Lemma~\ref{lem:mat-ramsey}
there exists a copy $T'$ of $B_{ m}$ in $T_{2^k,m}$ such that the
colors of pairs $\{x,y\}\in \mathrm{SP}(T')$ only depend on the
levels of $x$ and $y$.

Let $P$ be  a root-leaf path in $T'$ (isometric to $P_{ m}$). The
mapping $F|_P:P\to X$ has distortion at most $D$ by~\eqref{eq:F
vert}. Assuming
\begin{equation} \label{eq:assump1}
m\ge D^{16(t\log t)/\delta},
\end{equation}
by Proposition~\ref{prop:path-boosting} there are
$\{x_i\}_{i=0}^{t}\subseteq P$ such that for some $a,b\in \N$ with
$a,a+tb\in [0,m]$, for all $i$ we have $h(x_i)=a+ib$, and for some
$\theta>0$, for all $i,j\in \{0,\ldots,t\}$,
\begin{equation}\label{eq:restriction}
\theta b|i-j|\le d_X(F(x_i),F(x_j))\le
\left(1+\frac{\delta}{4}\right)\theta b|i-j|.
\end{equation}

Define a rescaled isometry $\f:B_t\to T'$ in a top-down manner as
follows: $\f(r)=x_0$, and having defined $\f(u)\in T'$, if
$v,w$ are the  children of $u$ in $B_t$ and $v',w'$ are the children of $\f(u)$ in $T'$, the vertices $\f(v),\f(w)$ are
chosen as arbitrary descendants in $T'$ of $v',w'$ (respectively) at depth
$h(\f(u))+b$. Consider the mapping $G:B_t\to X$ given by $G=F\circ
\f=f\circ g\circ \f$. Take $x,y\in B_t$ such that $x$ is an
ancestor of $y$. Write $h(x)=i$ and $h(y)=j$. Thus $h(\f(x))=a+ib$
and $h(\f(y))=a+jb$. It follows that $\{\f(x),\f(y)\}$ is
colored by the same color as $\{x_i,x_j\}$, i.e.,
$$
\left\lfloor \log_{1+\frac{\delta}4}
\left(\frac{d_X(G(x),G(y))}{k\ell\lambda bd_{B_t}(x,y)}\right)
\right\rfloor=\chi(\{\f(x),\f(y)\})=\chi(\{x_i,x_j\})=\left\lfloor
\log_{1+\frac{\delta}4} \left(\frac{d_X(F(x_i),F(y_j))}{k\ell\lambda
bd_{B_t}(x,y)}\right) \right\rfloor.
$$
Consequently, using~\eqref{eq:restriction} we deduce that
$$
\frac{\theta b}{1+\delta/4}d_{B_t}(x,y)\le d_X(G(x),G(y))\le
\left(1+\frac{\delta}{4}\right)^2\theta b d_{B_t}(x,y).
$$
Thus $G$ is a $(1+\delta/4)^3\le 1+\delta$ vertically faithful
embedding of $B_t$ into $X$.

It remains to determine the values of the auxiliary parameters
$\ell, k$, which will lead to the desired restriction on $n$ given
in~\eqref{eq:lower n}. First of all, we want to have $\dist(g\circ
\f)\le 1+\xi$. Since $\f$ is a rescaled isometry and (for
$\ell\ge 2$) $\dist (g)\le 1+2/\ell$, we choose $\ell=\lceil
2/\xi\rceil\ge 2$. We will choose $k$ so that $4k\le n\xi$, so that
$n/(k\ell)\ge 1$. Since $m=\left\lfloor n/(k\ell)\right\rfloor$, we
have $m+1\le n\xi/k$ and $m\ge n\xi/(4k)$. Recall that $r=\lceil
\log_{1+\delta/4}D \rceil\le 2\log_{1+\delta/4}D\le 16D/\delta$.
Hence the requirement~\eqref{eq:coloring assumption} will be
satisfied if
\begin{equation}\label{eq:cubed}
2^{k^3}\ge \left(\frac{16D}{\delta}\right)^{n^2\xi^2},
\end{equation}
and the requirement~\eqref{eq:assump1} will be satisfied if
\begin{equation}\label{eq:second k}
\frac{n\xi}{4k}\ge D^{16(t\log t)/\delta}.
\end{equation}
There exists an integer $k$ satisfying both~\eqref{eq:cubed}
and~\eqref{eq:second k} provided that
\begin{equation*}\label{eq:range k}
\left(n^2\xi^2
\log_2\left(\frac{16D}{\delta}\right)\right)^{1/3}+1\le
\frac{n\xi}{4D^{16(t\log t)/\delta}},
\end{equation*}
which holds true provided the constant $c$ in~\eqref{eq:lower n} is
large enough.
\end{proof}

\section{Tree metrics do not have the dichotomy property}
\label{sec:charlie-fails}

This section is devoted to the proofs of Theorem~\ref{lem:Bn->X} and
Theorem~\ref{thm:no local rigidity}. These proofs were outlined in
Section~\ref{sec:overview}, and we will use the notation introduced
there.

\subsection{Horizontally contracted trees} \label{sec:charlie-def}

We start with the following lemma which supplies conditions on
$\{\e_n\}_{n=0}^\infty$ ensuring that the $H$-tree $(B_\infty,d_\e)$
is a metric space.

\begin{lemma} \label{lem:H-tree-metric}
Assume that $\{\e_n\}_{n=0}^\infty\subseteq (0,1]$ is non-increasing and $\{n\e_n\}_{n=0}^\infty$ is non-decreasing. Then $d_\e$ is a metric on $B_\infty$
\end{lemma}
\begin{proof}
 Take $x,y,z\in B_\infty$ and without loss of
generality assume that $h(x)\le h(y)$. We distinguish between the
cases $h(z)> h(y)$, $h(x)\le h(z)\le h(y)$ and $h(z)<h(x)$.

If $h(z)>h(y)$ then
\begin{eqnarray}\label{eq:first case}
&&\!\!\!\!\!\!\!\!\!\!\!\!\!\!\!\!\!\!d_\e(x,z)+d_\e(z,y)-d_\e(x,y)\nonumber\\&=&2[h(z)-h(y)]+2\e_{h(x)}\cdot
\left[h(\lca(x,y))-h(\lca(x,z))\right]+2\e_{h(y)}\cdot
\left[h(y)-h(\lca(z,y))\right]\nonumber\\
&\ge& 2\e_{h(x)}\cdot
\left[h(\lca(x,y))-h(\lca(x,z))\right]+2\e_{h(y)}\cdot
\left[h(y)-h(\lca(z,y))\right].
\end{eqnarray}
To show that~\eqref{eq:first case} is non-negative observe that
this is obvious if $h(\lca(x,y))\ge h(\lca(x,z))$. So assume that
$h(\lca(x,y))< h(\lca(x,z))$. In this case necessarily
$h(\lca(z,y))= h(\lca(x,y))$, so we can bound~\eqref{eq:first
case} from below as follows
\begin{eqnarray*}
&&\!\!\!\!\!\!\!\!\!\!\!\!\!\!\!\!\!\!2\e_{h(x)}\cdot
\left[h(\lca(x,y))-h(\lca(x,z))\right]+2\e_{h(y)}\cdot
\left[h(y)-h(\lca(z,y))\right]\\&\ge& 2\e_{h(x)}\cdot
\left[h(\lca(x,y))-h(x)\right]+
2\frac{h(x)}{h(y)}\e_{h(x)}\cdot\left[h(y)-h(\lca(x,y))\right]\\
&=&2\e_{h(x)}\cdot h(\lca(x,y))\left(1-\frac{h(x)}{h(y)}\right)\ge
0.
\end{eqnarray*}

If $h(z)< h(x)$ then
\begin{eqnarray}\label{eq:second case}
&&\!\!\!\!\!\!\!\!\!\!\!\!\!\!\!\!\!\!\!\!d_\e(x,z)+d_\e(z,y)=h(x)+h(y)-2h(z)+2\e_{h(z)}\cdot
\left[2h(z)-h(\lca(x,z))-h(\lca(z,y))\right]\nonumber\\&\ge&
h(y)-h(x)+2\e_{h(x)}\cdot
[h(x)-h(z)]\nonumber\nonumber+2\e_{h(x)}\cdot
\left[2h(z)-h(\lca(x,z))-h(\lca(y,z))\right]\nonumber\\&=&
h(y)-h(x)+2\e_{h(x)}\cdot
\left[h(x)+h(z)-h(\lca(x,z))-h(\lca(y,z))\right]\nonumber\\ &\ge&
h(y)-h(x)+2\e_{h(x)}\cdot \left[h(x)-h(\lca(x,y))\right]\\ &=&
d_\e(x,y)\nonumber.
\end{eqnarray}
Where in~\eqref{eq:second case} we used the fact that $h(z)\ge
h(\lca(x,z))+h(\lca(y,z))-h(\lca(x,y))$, which is true since
$h(\lca(x,y))\ge \min\left\{h(\lca(x,z)),h(\lca(y,z))\right\}$.

It remains to deal with the case $h(x)\le h(z)\le h(y)$. In this
case

\begin{eqnarray}\label{eq:third case}
&&\!\!\!\!\!\!\!\!\!\!\!\!\!\!\!\!\!\!\!\!d_\e(x,z)+d_\e(z,y)=h(y)-h(x)+2\e_{h(x)}\cdot
\left[h(x)-h(\lca(x,z))\right]+2\e_{h(z)}\cdot
\left[h(z)-h(\lca(y,z))\right]\nonumber\\&\ge&\nonumber
h(y)-h(x)+2\e_{h(x)}\cdot
\left[h(x)-h(\lca(x,z))\right]+2\frac{h(x)}{h(z)}\e_{h(x)}\cdot\left[h(z)-h(\lca(y,z))\right]
\nonumber\\&=& h(y)-h(x)+2\e_{h(x)}\cdot
\left[2h(x)-h(\lca(x,z))-\frac{h(x)}{h(z)}h(\lca(y,z))\right]\nonumber\\
&\ge& h(y)-h(x)+2\e_{h(x)}\cdot \left[h(x)-h(\lca(x,y))\right]
\\ &=& d_\e(x,y)\nonumber,
\end{eqnarray}
where~\eqref{eq:third case} is equivalent to the inequality
\begin{eqnarray}\label{eq:explain-2}
h(x)\ge h(\lca(x,z))+\frac{h(x)}{h(z)}h(\lca(y,z))-h(\lca(x,y)).
\end{eqnarray}
To prove~\eqref{eq:explain-2}, note that it is true if $h(\lca(x,y))\ge
h(\lca(x,z))$, since clearly $h(\lca(y,z))\le h(z)$. If, on the
other hand, $h(\lca(x,y))< h(\lca(x,z))$ then using the assumption
that $h(z)\ge h(x)$ it is enough to show that $ h(x)\ge
h(\lca(x,z))+h(\lca(y,z))-h(\lca(x,y))$. Necessarily
$h(\lca(x,y))=h(\lca(y,z))$, so that the required inequality follows
from the fact that $h(x)\ge h(\lca(x,z))$.
\end{proof}

\subsection{Geometry of $H$-trees}\label{sec:geometry}

\subsubsection{Classification of approximate midpoints} \label{sec:midpoints}

 From now on we will always assume that $\e=\{\e_n\}_{n=0}^\infty$
satisfies for all $n\in \mathbb N$,
$\e_n\ge \e_{n+1}>0$ and $(n+1)\e_{n+1}\ge n\e_n$. We recall the important
concept of {\em approximate midpoints} which is used frequently in
nonlinear functional analysis (see~\cite{BL} and the references
therein).

\begin{definition}[Approximate midpoints] Let $(X,d_X)$ be a metric space and $\delta\in
(0,1)$. For $x,y\in X$ the set of $\delta$-approximate midpoints
of $x$ and $z$ is defined as
\begin{equation*}
\Mid(x,z,\delta)=\left\{y\in X:\, \max\{d_X(x,y),d_X(y,z)\}\le
\frac{1+\delta}{2}\cdot d_X(x,z)\right\}.
\end{equation*}
\end{definition}

\begin{figure}[ht]
\begin{center}\includegraphics[scale=0.8]{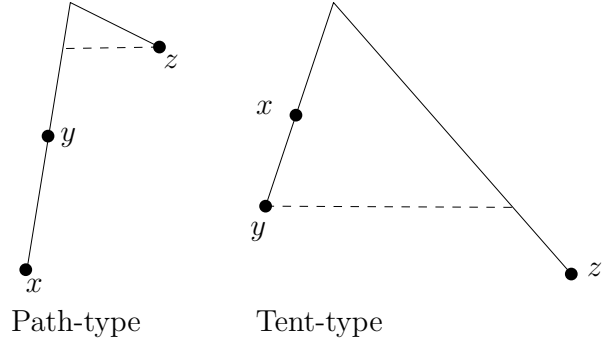}
\end{center}
 \caption{{{\em A schematic description of path-type and
 tent-type configurations. }}} \label{fig:midpoints}
\end{figure}

From now on, whenever we refer to the set $\Mid(x,z,\delta)$, the
underlying metric will always be understood to be $d_\e$. In what
follows, given $\eta>0$ we shall say that two sequences
$(u_1,\ldots,u_n)$ and $(v_1,\ldots,v_n)$ of vertices in $B_\infty$
are $\eta$-near if for every $j\in \{1,\ldots,n\}$ we have
$d_\e(u_j,v_j)\le \eta$. We shall also require the following
terminology:
\begin{definition}
An ordered triple
$(x,y,z)$ of vertices in $B_\infty$ will be called a {\em path-type
configuration} if $h(z)\le h(y)\le h(x)$, $x$ is a descendant of
$y$, and $h(\lca(z,y))<h(y)$. The triple $(x,y,z)$ will be called a
{\em tent-type configuration} if $h(y)\le h(z)$, $y$ is a descendant
of $x$, and $h(\lca(x,z))<h(x)$. These special configurations are
described in Figure~\ref{fig:midpoints}.
\end{definition}

The following useful theorem will be used extensively in the
ensuing arguments.
Its proof will be broken down into several elementary lemmas.

\begin{theorem}\label{thm:midpoints}
Assume that $\delta\in(0,\tfrac{1}{16})$, and the sequence $\e=\{\e_n\}_{n=0}^\infty$ satisfies $\e_n<\tfrac{1}{4}$ for all $n\in \N$.
Let $x,y,z\in (B_\infty, d_\e)$ be
such that $y\in \Mid(x,z,\delta)$. Then either $(x,y,z)$ or
$(z,y,x)$ is $3\delta d_\e(x,z)$-near a path-type or tent-type
configuration.
\end{theorem}

In what follows, given a vertex $v\in B_\infty$ we denote the
subtree rooted at $v$ by $T_v$.

\begin{lemma}\label{lem:monotonicity} Assume that $\e_n\le \frac12$ for all $n$. Fix $a\in B_\infty$ and let $u,v\in B_\infty$ be its
children. For every $x,z\in T_u$ such that $h(x)\ge h(z)$ consider
the function $D_{x,z}:\{a\}\cup T_v\to [0,\infty)$ defined by
$D_{x,z}(y)=d_\e(x,y)+d_\e(z,y)$. Fix an arbitrary vertex $w\in
T_v$ such that $h(w)=h(z)$. Then for every $y\in T_v$ we have
$D_{x,z}(y)\ge D_{x,z}(w)$.
\end{lemma}

\begin{proof} By the definition of $d_\e$ we have
$D_{x,z}(y)=Q(h(y))$ where
\begin{multline*}
Q(k)=\max\{h(x),k\}+\max\{k,h(z)\}-\min\{h(x),k\}-\min\{k,h(z)\}\\+
2\e_{\min\{h(x),k\}}\left[\min\{h(x),k\}-h(a)\right]+2\e_{\min\{k,h(z)\}}\left[\min\{k,h(z)\}-h(a)\right].
\end{multline*}
The required result will follow if we show that $Q$ is
non-increasing on $\{h(a),h(a)+1,\ldots,h(z)\}$ and non-decreasing
on $\{h(z),h(z)+1,\ldots\}$. If $k\in
\{h(a),h(a)+1,\ldots,h(z)-1\}$ then
\begin{multline*}
Q(k+1)-Q(k)=-2+4\e_{k+1}[k+1-h(a)]-4\e_k[k-h(a)]\\\le
-2+4\e_k[k+1-h(a)]-4\e_k[k-h(a)] =-2(1-2\e_k)\le 0.
\end{multline*}
If $k\in \{h(x),h(x)+1,\ldots\}$ then $Q(k+1)-Q(k)=2$, and if
$k\in \{h(z),\ldots,h(x)-1\}$ then
\begin{eqnarray*}
Q(k+1)-Q(k)=2\left[(k+1)\e_{k+1}-k\e_k\right]+2h(a)[\e_k-\e_{k+1}]\ge
0.
\end{eqnarray*}
This completes the proof of Lemma~\ref{lem:monotonicity}.
\end{proof}

\begin{lemma}\label{lem:nonexistant case} Assume that $\e_n<\frac 12$ for all $n\in \N$. Fix $\delta\in \left(0,\frac13\right)$ and $x,y,z\in B_\infty$ such that $h(x)\ge h(z)$, $y\in \Mid(x,z,\delta)$ and
  $h(\lca(x,z))> h(\lca(x,y))$. Then
  $$h(z)+\frac{1-3\delta}{2} d_\e(x,z)\le h(y)<h(x)\le h(y)+\frac{1+3\delta}{1-3\delta}[h(y)-h(z)].$$ Moreover, if $y'\in
  B_\infty$ is the point on the segment joining $x$ and
$\lca(x,y)$ such that $h(y')=h(y)$ then $d_\e(y,y')\le \delta
d_\e(x,z)$. Thus  $(x,y',z)$ is a path-type configuration which is
$\delta d_\e(x,z)$-near $(x,y,z)$
  \end{lemma}

\begin{proof} Write $a=\lca(x,y)$. If $u,v$ are the two
children of $a$, then without loss of generality $x,z\in T_u$ and
$y\in T_v$. Let $w\in T_v$ be such that $h(w)=h(z)$.

\parpic[r]{\includegraphics[scale=0.7]{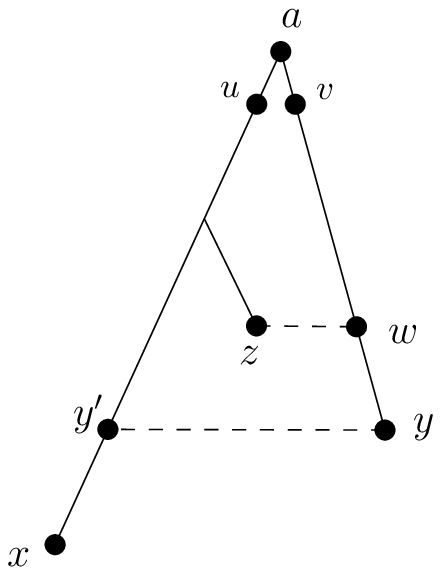}}By
Lemma~\ref{lem:monotonicity},
\begin{multline}\label{eq:use monotone}
d_\e(x,y)+d_\e(z,y)\ge d_\e(x,w)+d_\e(z,w)\\=
h(x)-h(z)+4\e_{h(z)}[h(z)-h(a)]\ge 2d_\e(x,z)-[h(x)-h(z)].
\end{multline}
On the other hand, since $y\in \Mid(x,z,\delta)$, we have that $$
d_\e(x,y)+d_\e(z,y)\le (1+\delta)d_\e(x,z).$$ Additionally, by the
definition of $d_\e$ we know that if $h(y)\le h(z)$ then
$$\frac{1+\delta}{2}d_\e(x,z)\ge d_\e(x,y)\ge h(x)-h(y)\ge
h(x)-h(z).$$ Combining these observations with~\eqref{eq:use monotone} we get that
\begin{eqnarray}\label{eq:conclude}
(1+\delta)d_\e(x,z)\ge 2d_\e(x,z)-\frac{1+\delta}{2}d_\e(x,z),
\end{eqnarray}
which is a contradiction since $\delta<\frac13$. Therefore
$h(y)>h(z)$. If $h(y)\ge h(x)$ then
$$\frac{1+\delta}{2}d_\e(x,z)\ge d_\e(z,y)\ge h(y)-h(z)\ge
h(x)-h(z),$$ so that we arrive at a contradiction as
in~\eqref{eq:conclude}. We have thus shown that $h(z)<h(y)<h(x)$.

Now, since $y\in \Mid(x,z,\delta)$,
\begin{multline*}
h(x)-h(z)+2\e_{h(y)}[h(y)-h(a)]+2\e_{h(z)}[h(z)-h(a)]=d_\e(x,y)+d_\e(z,y)
\\\le(1+\delta)d_\e(x,z)=
(1+\delta)\left(h(x)-h(z)+2\e_{h(z)}[h(z)-h(a)]\right).
\end{multline*}
Thus,  letting $y'$ be the point on the segment joining $x$ and
$a$ such that $h(y')=h(y)$, we see that
\begin{equation*}
d_\e(y,y')= 2\e_{h(y)}[h(y)-h(a)]\le
\delta\left(h(x)-h(z)+2\e_{h(z)}[h(z)-h(a)]\right)=\delta
d_\e(x,z),
\end{equation*}
Moreover
\begin{multline*}
\frac{1-\delta}{2}d_\e(x,z)\le
d_\e(y,z)=h(y)-h(z)+2\e_{h(z)}h(z)-2\e_{h(z)}h(a)\\
\le h(y)-h(z)+2\e_{h(y)}h(y)-2\e_{h(y)}h(a)\le h(y)-h(z)+\delta
d_\e(x,z).
\end{multline*}
Thus
\begin{eqnarray}\label{eq:large top interval}
h(y)-h(z)\ge \frac{1-3\delta}{2} d_\e(x,z).
\end{eqnarray}
Hence,
\begin{eqnarray*}
\frac{2}{1-3\delta}[h(y)-h(z)]\stackrel{\eqref{eq:large top
interval}}{\ge}d_\e(x,z)=[h(x)-h(y)]+[h(y)-h(z)].
\end{eqnarray*}
It follows that
\begin{eqnarray*}
h(x)-h(y)\le \frac{1+3\delta}{1-3\delta}[h(y)-h(z)].
\end{eqnarray*}
 This completes the proof of Lemma~\ref{lem:nonexistant case}.
\end{proof}

\begin{lemma}\label{lem:mid-lemma2} Assume that $\e_n<\frac14$ for all $n\in \mathbb N$. Fix $\delta\in \left(0,\frac{1}{16}\right)$ and assume that $x,y,z\in B_\infty$
are distinct vertices such that $\lca(x,y)=\lca(x,z)$, and $y\in\Mid(x,z,\delta)$.
Then either $(x,y,z)$ or $(z,y,x)$ is $3\delta d_\e(x,z)$-near a
path-type or tent-type configuration.
\end{lemma}

\begin{proof} Denote $a=\lca(x,y)$. Our assumption implies that $h(\lca(z,y))\ge h(a)$.
We perform a case analysis on the relative heights of $x,y,z$. Assume first that $h(x)\le h(y)$.

\begin{comment}
A schematic description of
the situation we are dealing with is contained in Figure~\ref{fig:second}.

\bigskip
\begin{figure}[h]
\begin{center}\includegraphics[scale=0.8]{case1}
\end{center}
 \caption{{{\em  A schematic description of the configuration of points in Lemma~\ref{lem:tent type}. }}} \label{fig:case1}
\end{figure}
\bigskip
\end{comment}

\parpic[r]{\includegraphics[scale=0.7]{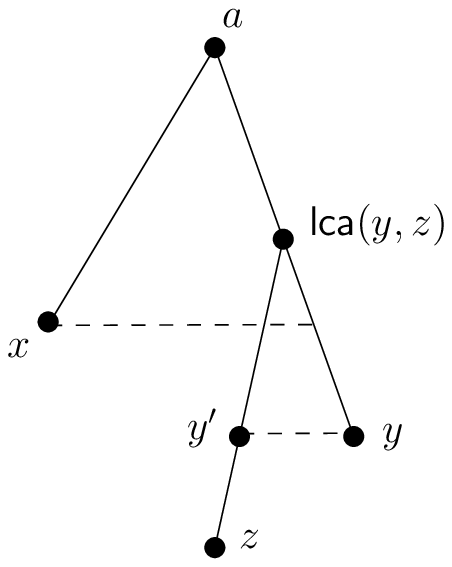}}
If $h(x)\le h(y) \le h(z)$ then
\begin{align}
(1+\delta)d_\e(x,z)&\ge \nonumber d_\e(x,y)+d_\e(y,z)\\
\nonumber &= h(y)-h(x)+2\e_{h(x)}[h(x)-h(a)]\\
\nonumber & \qquad +h(z)-h(y)+2\e_{h(y)}[h(y)-h(\lca(z,y))]\\
&= d_\e(x,z)+2\e_{h(y)}[h(y)-h(\lca(z,y))]. \label{pass to y'}
\end{align}
Let $y'$ be the point on the path from $\lca(y,z)$ to $z$ such
that $h(y')=h(y)$. Then~\eqref{pass to y'} implies that
$$
d_\e(y,y')=2\e_{h(y)}[h(y)-h(\lca(z,y))]\le \delta d_\e(x,z).
$$
Thus the triple $(z,y',x)$ is  a configuration of path-type which
is $\delta d_\e(x,z)$-near $(z,y,x)$.

\parpic[r]{\includegraphics[scale=0.7]{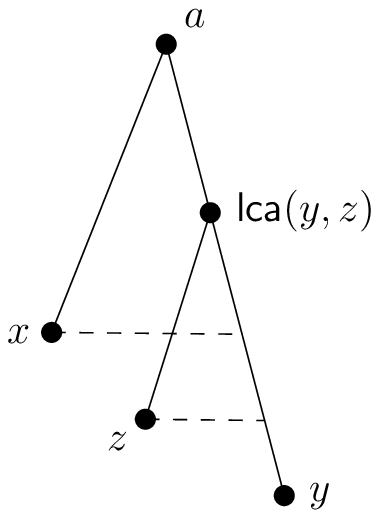}}
If $h(x)\le h(z)\le h(y)$ then since
\begin{eqnarray*}
\frac{1+\delta}{2}d_\e(x,z)\ge
d_\e(x,y)=h(y)-h(x)+2\e_{h(x)}[h(x)-h(a)]
\end{eqnarray*}
and $d_\e(x,z)=h(z)-h(x)+2\e_{h(x)}[h(x)-h(a)]$ we deduce that
\begin{eqnarray*}
-\frac{1-\delta}{2}d_\e(x,z)\ge d_\e(x,y)-d_\e(x,z)=h(y)-h(z)\ge
0.
\end{eqnarray*}
It follows that $x=z$, in contradiction to our assumption.

\parpic[r]{\includegraphics[scale=0.7]{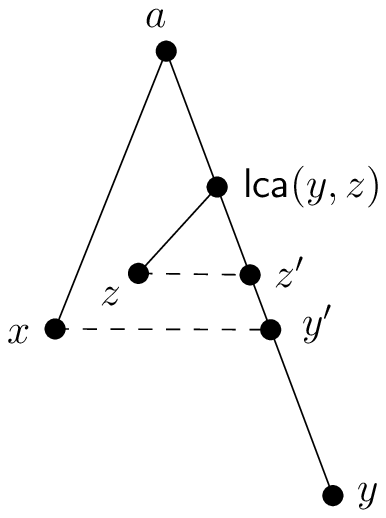}}
If $h(z)<h(x)$ then let $z'$ be the point on the segment joining $a$
and $y$ such that $h(z')=h(z)$. We thus have that
$$d_\e(z,z')=2\e_{h(z)}[h(z)-h(\lca(y,z))]=d_\e(z,y)-[h(y)-h(z)].$$
Moreover,
\begin{align*}
2\e_{h(z)}[h(\lca(z,y))-h(a)]&=d_\e(x,z)-[h(x)-h(z)]-2\e_{h(z)}[h(z)-h(\lca(z,y))]\\
&\ge \frac{2}{1+\delta}d_\e(z,y)-[h(x)-h(z)]-2\e_{h(z)}[h(z)-h(\lca(z,y))]\\
&=\frac{2}{1+\delta}\left(h(y)-h(z)+2\e_{h(z)}[h(z)-h(\lca(y,z))]\right)\\
&\phantom{\le}-\left(h(y)-h(z)+
2\e_{h(z)}[h(z)-h(\lca(y,z))]\right)+[h(y)-h(x)]\\
&=\frac{1-\delta}{1+\delta}d_\e(y,z)+[h(y)-h(x)]\displaybreak[0] \\
&\ge \frac{1-\delta}{1+\delta}\cdot
\frac{1-\delta}{2}d_\e(x,z)+[h(y)-h(x)]\\
&\ge
\left(\frac{1-\delta}{1+\delta}\right)^2d_\e(x,y)+[h(y)-h(x)]\\
&= d_\e(x,y)+[h(y)-h(x)]-\frac{4\delta}{(1+\delta)^2}d_\e(x,y)\\
&= 2[h(y)-h(x)]+2\e_{h(x)}[h(x)-h(a)]-\frac{4\delta}{(1+\delta)^2}d_\e(x,y)\\
&\ge 2[h(y)-h(x)]+ 2\e_{h(z)}h(z)-2\e_{h(z)}h(a)-\frac{2\delta}{1+\delta}d_\e(x,z).
\end{align*}
Thus
\begin{eqnarray}\label{eq:both of them}
\frac{2\delta}{1+\delta}d_\e(x,z)\ge
2[h(y)-h(x)]+2\e_{h(z)}[h(z)-h(\lca(y,z))]=2[h(y)-h(x)]+d_\e(z,z').
\end{eqnarray}
Let $y'$ be the point on the path from $a$ to $y$ such that
$h(y')=h(x)$. It follows from~\eqref{eq:both of them} that the
triple $(z',y',x)$ is a configuration of tent-type which is
$2\delta d_\e(x,z)$-near $(z,y,x)$.

This completes the proof of Lemma~\ref{lem:mid-lemma2} when
$h(x)\le h(y)$. The case $h(x)>h(y)$ is proved analogously. Here
are the details.

\begin{comment}
\bigskip
\begin{figure}[h]
\begin{center}\includegraphics[scale=0.8]{case2}
\end{center}
 \caption{{{\em  A schematic description of the configuration of points in Lemma~\ref{lem:tent type}. }}} \label{fig:case2}
\end{figure}
\bigskip
\end{comment}

\parpic[r]{\includegraphics[scale=0.7]{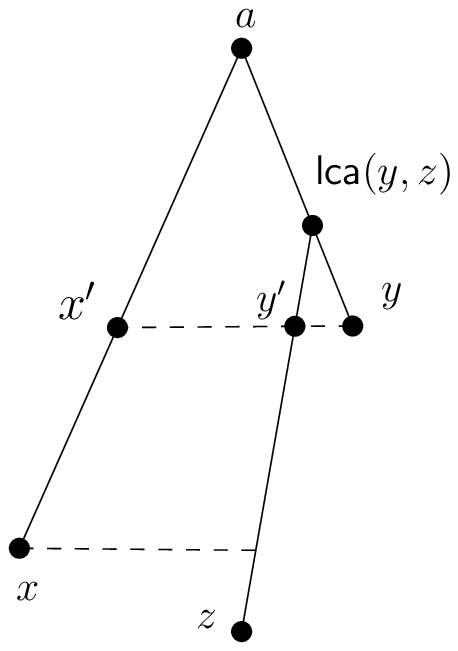}}
Assume first of all that $h(z)\ge h(x) >h(y)$. Then,
\begin{align}
d_\e(x,z)&\ge \frac{2}{1+\delta}d_\e(z,y)\nonumber
= 2d_\e(z,y)-\frac{2\delta}{1+\delta}d_\e(z,y)\nonumber\\
&= 2\left(h(z)-h(y)+2\e_{h(y)}[h(y)-h(\lca(z,y))]\right)-\frac{2\delta}{1+\delta}d_\e(z,y) \nonumber \\
&\ge 2\left(h(z)-h(y)+2\e_{h(y)}[h(y)-h(\lca(z,y))]\right)-\delta d_\e(x,z).\label{eq:reverse again}
\end{align}
On the other hand, since $h(x)>h(y)$,
\begin{eqnarray*}\label{eq:before combine}
d_\e(x,z)&=&h(z)-h(x)+2\e_{h(x)}[h(x)-h(a)]\nonumber\\
&\le& h(z)-h(x)+2\e_{h(y)}[h(x)-h(a)]\nonumber
\\
\nonumber&=&\left(h(x)-h(y)+2\e_{h(y)}[h(y)-h(a)]\right)+h(y)+h(z)-2h(x)+2\e_{h(y)}[h(x)-h(y)]\\
\nonumber &=& d_\e(x,y)+h(y)+h(z)-2h(x)+2\e_{h(y)}[h(x)-h(y)]\\
\nonumber &\le& \frac{1+\delta}{1-\delta}d_\e(y,z)+h(y)+h(z)-2h(x)+2\e_{h(y)}[h(x)-h(y)]\\
\nonumber&=& \left(h(z)-h(y)+2\e_{h(y)}[h(y)-h(\lca(z,y))]\right)+\frac{2\delta}{1-\delta}d_\e(y,z)\\\nonumber
&\phantom{\le}& +h(y)+h(z)-2h(x)+2\e_{h(y)}[h(x)-h(y)]\\
&\le& 2[h(z)-h(x)]+2\e_{h(y)}[h(x)-h(\lca(z,y))]+\tfrac{1+\delta}{1-\delta}\delta
d_\e(x,z).
\end{eqnarray*}
Combining this bound with~\eqref{eq:reverse again},
and canceling terms, gives
\begin{eqnarray}\label{eq:path appeared}
\frac{2\delta}{1-\delta}d_\e(x,z)&\ge&
2[h(x)-h(y)]-4\e_{h(y)}[h(x)-h(y)]+2\e_{h(y)}[h(x)-h(\lca(z,y))]\nonumber\\
&\ge& \nonumber
2\left(1-2\e_{h(y)}\right)[h(x)-h(y)]+2\e_{h(y)}[h(x)-h(\lca(z,y))]\\&>&[h(x)-h(y)]+2\e_{h(y)}[h(x)-h(\lca(z,y))],
\end{eqnarray}
where we used the fact that $\e_{h(y)}<\frac14$. Let $x'$ be the
point on the path from $x$ to $a$ such that $h(x')=h(y)$, and let
$y'$ be the point on the path from $a$ to $z$ such that
$h(y')=h(y)$. Then by~\eqref{eq:path appeared}
$d_\e(x,x')=h(x)-h(y)\le 3\delta d_\e(x,z)$ and
$$d_\e(y,y')=2\e_{h(y)}[h(y)-h(\lca(z,y))]\le
2\e_{h(y)}[h(x)-h(\lca(z,y))]\le 3\delta d_\e(x,z).$$ Thus the
triple $(z,y',x')$ is a configuration of path-type which is
$3\delta d_\e(x,z)$-near $(z,y,x)$.

\parpic[r]{\includegraphics[scale=0.7]{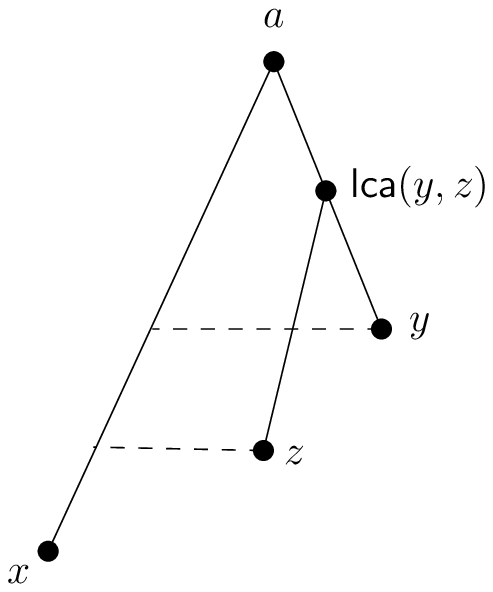}}
If $h(x)>h(z)\ge h(y)$ then
\begin{multline*}
h(z)-h(y)+2\e_{h(y)}[h(y)-h(\lca(z,y))]=d_\e(z,y) \ge
\frac{1-\delta}{1+\delta}d_\e(x,y)\\=
d_\e(x,y)-\frac{2\delta}{1+\delta}
d_\e(x,y)=h(x)-h(y)+2\e_{h(y)}[h(y)-h(a)]-\frac{2\delta}{1+\delta}
d_\e(x,y).
\end{multline*}
Canceling terms we see that
\begin{align*}
\frac{2\delta}{1+\delta} d_\e(x,y)&\ge
h(x)-h(z)+2\e_{h(y)}[h(\lca(z,y))-h(a)]\\
&=h(x)-h(z)+2\e_{h(y)}[h(z)-h(a)]-
2\e_{h(y)}[h(z)-h(\lca(z,y))]\\
&\ge h(x)-h(z)+2\e_{h(z)}[h(z)-h(a)]-
2\e_{h(y)}[h(z)-h(\lca(z,y))]\\
&= d_\e(x,z)- 2\e_{h(y)}[h(z)-h(\lca(z,y))] \displaybreak[0]\\
&\ge 2d_\e(z,y)-\frac{2\delta}{1+\delta}d_\e(x,z)-
2\e_{h(y)}[h(z)-h(\lca(z,y))]\\
&=
2\left(h(z)-h(y)+2\e_{h(y)}[h(y)-h(\lca(z,y))]\right)-\frac{2\delta}{1+\delta}d_\e(x,z)\\
&\phantom{\le} -
2\e_{h(y)}[h(z)-h(\lca(z,y))]\\
&= 2\left(1-2\e_{h(y)}\right)[h(z)-h(y)]+2\e_{h(y)}[h(z)-h(\lca(z,y))]-\frac{2\delta}{1+\delta}d_\e(x,z)\\
&\ge [h(z)-h(y)]+2\e_{h(y)}[h(z)-h(\lca(z,y))]-\frac{2\delta}{1+\delta}d_\e(x,z)\\
&= d_\e(z,y)-\frac{2\delta}{1+\delta}d_\e(x,z)\\
&\ge \left(\frac{1-\delta}{2}-\frac{2\delta}{1+\delta}\right)d_\e(x,z), %\label{eq:second noncase}
\end{align*}
which is a contradiction since $\delta<\frac{1}{16}$.

\parpic[r]{\includegraphics[scale=0.7]{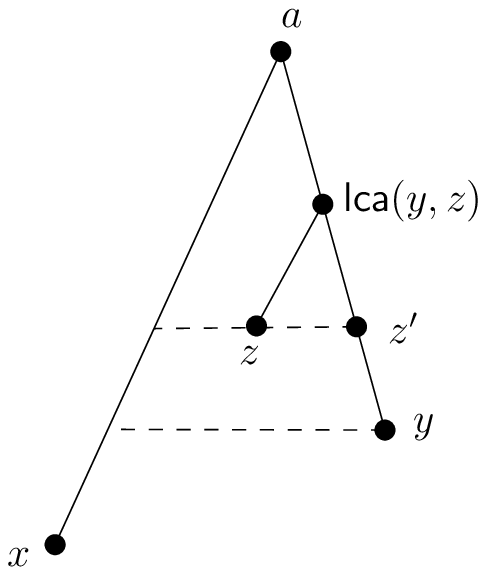}}
The only remaining case is when $h(x)>h(y)>h(z)$. In this case we
proceed as follows.
\begin{align}
 d_\e(x,z)&=h(x)-h(z)+2\e_{h(z)}[h(z)-h(a)]\nonumber\\
 &= d_\e(y,z)+[h(x)-h(y)]+2\e_{h(z)}[h(\lca(y,z))-h(a)]\nonumber\\
 &\ge d_\e(x,y)-\frac{2\delta}{1+\delta}d_\e(x,y)+[h(x)-h(y)]+2\e_{h(z)}[h(\lca(y,z))-h(a)]\nonumber
 \displaybreak[0]\\
 &\ge h(x)-h(y)+2\e_{h(y)}[h(y)-h(a)]-{\delta} d_\e(x,z)\nonumber\\
&\phantom{\le} +[h(x)-h(y)]+2\e_{h(z)}[h(\lca(y,z))-h(a)] \nonumber\\
&\ge 2[h(x)-h(y)]
+2\e_{h(z)}[h(z)-h(a)]+2\e_{h(z)}[h(\lca(y,z))-h(a)]
- {\delta}d_\e(x,z)\nonumber\\
&= 2[h(x)-h(y)]+4\e_{h(z)}[h(z)-h(a)]-2\e_{h(z)}[h(z)-h(\lca(y,z))]- {\delta} d_\e(x,z)\nonumber\\
&=
2d_\e(x,z)-2d_\e(z,y)+2\e_{h(z)}[h(z)-h(\lca(y,z))]- {\delta} d_\e(x,z)\nonumber\\
&\ge
\left(1-2\delta \right)d_\e(x,z)+2\e_{h(z)}[h(z)-h(\lca(y,z))]. \label{eq:finished}
\end{align}
Let $z'$ be the point on the path from $a$ to $y$ such that
$h(z')=h(z)$. Then
\begin{eqnarray*}
d_\e(z,z')=2\e_{h(z)}[h(z)-h(\lca(y,z))]\stackrel{\eqref{eq:finished}}{\le} 2\delta d_\e(x,z).
\end{eqnarray*}
Therefore the triple $(z',y,x)$ is of tent-type and is $2\delta
d_\e(x,z)$-near $(z,y,x)$. The proof of Lemma~\ref{lem:mid-lemma2} is
complete.
\end{proof}

\begin{proof}[Proof of Theorem~\ref{thm:midpoints}]
It remains to check that for every $x,y,z\in B_\infty$ such that $y\in \Mid(x,z,\delta)$, at least one of the triples
$(x,y,z)$ or $(z,y,x)$ satisfies the conditions of Lemma~\ref{lem:nonexistant case} or Lemma~\ref{lem:mid-lemma2}.

Indeed, if $h(\lca(x,y))=h(\lca(x,z))$ then $\lca(x,y)=\lca(x,z)$, so Lemma~\ref{lem:mid-lemma2} applies. If $h(\lca(x,y))<h(\lca(x,z))$ then $\lca(z,y)=\lca(x,z)$, so Lemma~\ref{lem:mid-lemma2} applies to the triple $(z,y,x)$. If $h(\lca(x,y))<h(\lca(x,z))$ then $\lca(x,y)=\lca(z,y)$, and so $h(\lca(x,z))>h(\lca(z,y))$. Hence Lemma~\ref{lem:nonexistant case} applies to either the triple $(x,y,z)$ or the triple $(z,y,x)$.
\end{proof}

We end this subsection with a short discussion on the distance between tent-type and path-type configurations.
It turns out that when $\e_h\ll \delta$, a $\delta$-midpoint configuration $(x,y,z)$ can be close to
a path-type configuration, and at the same time the reversed triple $(z,y,x)$ close to a tent-type configuration (or vice versa). However, it is easy to see that this is the only ``closeness" possible.

\begin{lemma} \label{lem:unique-midpoint-configuration} Fix $x,y,z\in B_\infty$ with $x\neq y$. Then the following statements are impossible:
% Let $x,y,z\in B_\infty$ be distinct points such that $y\in\Mid(x,z,\delta)$, with $\delta<1/3$. Then the following statements are impossible:
\begin{enumerate}
\item $(x,y,z)$ is $\frac15 d_\e(x,y)$-near a path-type configuration and a tent-type configuration.
\item $(x,y,z)$ is $\frac1{11}d_\e(x,y)$-near a path-type configuration and $(z,y,x)$ is  $\frac1{11}d_\e(x,y)$-near a path-type configuration.
\item $(x,y,z)$ is $\frac{1}{11} d_\e(x,y)$-near a tent-type configuration and $(z,y,x)$ is  $\frac1{11}d_\e(x,y)$-near a tent-type configuration.
\end{enumerate}
\end{lemma}
\begin{proof}
%The key point is that in both a tent-type configuration and a path-type configuration there is pair of points in ancestor/descendant relation.
%But, in all the pairs of configurations listed above, this pair of points is reversed.
%Here is the detailed argument in case 1 above. The remaining proofs are entirely analogous.
For case 1 of Lemma~\ref{lem:unique-midpoint-configuration}, assume for contradiction that $(x,y,z)$ is $\frac15 d_\e(x,y)$-near a path-type configuration $(a_1,b_1,c_1)$, and also $\frac15 d_\e(x,y)$-near a tent-type configuration $(\alpha_1,\beta_1,\gamma_1)$. By the definitions of path-type and tent-type configurations, $a_1$ is a descendant of $b_1$ and $\beta_1$ is a descendant of $\alpha_1$. Hence,
\begin{equation}\label{eq:ab}
h(a_1)-h(b_1)=d_\e(a_1,b_1)\ge d_\e(x,y)-d_\e(x,a_1)-d_\e(y,b_1)\ge \frac35 d_\e(x,y),
\end{equation}
and
\begin{equation}\label{eq:alpha beta}
h(\beta_1)-h(\alpha_1)=d_\e(\alpha_1,\beta_1)\ge d_\e(x,y)-d_\e(x,\alpha_1)-d_\e(y,\beta_1)\ge \frac35d_\e(x,y).
\end{equation}
By summing~\eqref{eq:ab} and~\eqref{eq:alpha beta} we see that,
\begin{multline*}
\frac{4}{5} d_\e(x,y)\ge d_\e(a_1,x)+d_\e(x,\alpha_1)+d_\e(b_1,y)+d_\e(y,\beta_1)\ge d_\e(a_1,\alpha_1)+d_\e(b_1,\beta_1)\\\ge h(a_1)-h(\alpha_1)+h(\beta_1)-h(b_1)
\stackrel{\eqref{eq:ab}\wedge\eqref{eq:alpha beta}}{\ge} \frac{6}{5} d_\e(x,y),
\end{multline*}
a contradiction.

For case 2 of Lemma~\ref{lem:unique-midpoint-configuration}, assume for contradiction that $(x,y,z)$ is $\frac1{11} d_\e(x,y)$-near a path-type configuration $(a_2,b_2,c_2)$, and also $(z,y,x)$ is $\frac1{11} d_\e(x,y)$-near a path-type configuration $(\alpha_2,\beta_2,\gamma_2)$. By the definitions of path-type and tent-type configurations, $a_2$ is a descendant of $b_2$ and $h(\beta_2)>h(\gamma_2)$. Hence,
\begin{equation}\label{eq:a2b2}
h(a_2)-h(b_2)=d_\e(a_2,b_2)\ge d_\e(x,y)-d_\e(x,a_2)-d_\e(y,b_2)\ge \frac9{11} d_\e(x,y),
\end{equation}
and
\begin{multline}\label{eq:alpha2beta2}
h(\beta_2)-h(\gamma_2)+2\e_{h(\gamma_2)}[h(\gamma_2)-h(\lca(\beta_2,\gamma_2))]=d_\e(\beta_2,\gamma_2)\\\ge d_\e(x,y)-d_\e(x,\gamma_2)-d_\e(y,\beta_2)\ge \frac9{11}d_\e(x,y).
\end{multline}
By summing~\eqref{eq:a2b2} and~\eqref{eq:alpha2beta2} we see that
\begin{eqnarray*}
\frac{17}{11}d_\e(x,y)&\ge&d_\e(a_2,x)+d_\e(x,\gamma_2)+d_\e(b_2,y)+d_\e(y,\beta_2)+d_\e(\beta_2,y)+d_\e(x,y)+d_\e(x,\gamma_2)\\&\ge& d_\e(a_2,\gamma_2)+d_\e(b_2,\beta_2)+d_\e(\beta_2,\gamma_2)\\&\ge& \big(h(a_2)-h(\gamma_2)\big)+\big(h(\beta_2)-h(b_2)\big)+ 2\e_{h(\gamma_2)}[h(\gamma_2)-h(\lca(\beta_2,\gamma_2))] \\&\stackrel{\eqref{eq:a2b2}\wedge \eqref{eq:alpha2beta2}}{\ge}& \frac{18}{11}d_\e(x,y),
\end{eqnarray*}
a contradiction.

For case 3 of Lemma~\ref{lem:unique-midpoint-configuration}, assume for contradiction that $(x,y,z)$ is $\frac 1{11} d_\e(x,y)$-near a tent-type configuration $(a_3,b_3,c_3)$, and also $(z,y,x)$ is $\frac1{11} d_\e(x,y)$-near a tent-type configuration $(\alpha_3,\beta_3,\gamma_3)$. Then $b_3$ is a descendant of $a_3$ and $h(\gamma_3)>h(\beta_3)$. Hence,
\begin{equation}\label{eq:a3b3}
h(b_3)-h(a_3)=d_\e(a_3,b_3)\ge d_\e(x,y)-d_\e(x,a_3)-d_\e(y,b_3)\ge \frac 9{11}d_\e(x,y),
\end{equation}
and
\begin{multline}\label{eq:alpha3beta3}
h(\gamma_3)-h(\beta_3)+2\e_{h(\beta_3)}[h(\beta_3)-h(\lca(\beta_3,\gamma_3))]=d_\e(\beta_3,\gamma_3)\\
\ge d_\e(x,y)-d_\e(x,\gamma_3)-d_\e(y,\beta_3)\ge \frac{9}{11}d_\e(x,y).
\end{multline}
Hence,
\begin{eqnarray*}
\frac{17}{11}d_\e(x,y)&\ge& d_\e(b_3,y)+d_\e(y,\beta_3)+d_\e(x,\alpha_3)+d_\e(x,\gamma_3)+d_\e(\beta_3,y)+d_\e(x,y)+d_\e(x,\gamma_3)\\&\ge& d_\e(b_3,\beta_3)+d_\e(a_3,\gamma_3)+d_\e(\beta_3,\gamma_3)\\&\ge& \big(h(b_3)-h(\beta_3)\big)+\big(h(\gamma_3)-h(a_3)\big)+2\e_{h(\beta_3)}[h(\beta_3)-h(\lca(\beta_3,\gamma_3))]
\\&\stackrel{\eqref{eq:a3b3}\wedge \eqref{eq:alpha3beta3}}{\ge}& \frac{18}{11}d_\e(x,y),
\end{eqnarray*}
a contradiction.
\end{proof}

%Note that in the notation of Theorem~\ref{thm:midpoints}, for $\delta<1/16$, $3d_\e(x,z)\le 7 d_\e(x,y)$.

\subsubsection{Classification of approximate forks}
\label{sec:forks}

We begin with three ``stitching lemmas" that roughly say that given three points $x,x',y\in(B_\infty,d_\e)$ such that
$x'$ is near $x$, there exists $y'$ near $y$ such that $d_\e(x',y')$ is close to $d_\e(x,y)$, and $y'$ relates to $x'$ in $B_\infty$
``in the same way'' that  $y$ relates $x$ in $B_\infty$.
% We check it for three cases: when $y$ is an ancestor of $x$ (Lemma~\ref{lem:ancestor}), and when
% $h(y)=h(x)$ (Lemma~\ref{lem:horizontal}).

\begin{lemma} \label{lem:ancestor}
Let $x,x',y,y'\in B_\infty$ be such that $y$ is an ancestor of $x$, and $y'$ is an ancestor of $x'$ satisfying
$h(x)-h(y)=h(x')-h(y')$. Then $d_\e(y,y')\leq d_\e(x,x')$.
\end{lemma}
\begin{proof}
Assume without loss of generality that $h(x)\geq h(x')$. So, $$d_\e(x,x')=h(x)-h(x') +2\e_{h(x')}[h(x')- h(\lca(x,x'))].$$
Note that $h(\lca(y,y')) = \min\{h(y'),h(\lca(x,x'))\}$. Hence,
\begin{align}
 d_\e(y,y') &= h(y)-h(y') + 2\e_{h(y')} [h(y') - h(\lca(y,y'))] \nonumber \\
 & = h(x)-h(x') + 2\e_{h(y')} [h(y') - \min\{h(y'), h(\lca(x,x'))\}] \nonumber \\
 & = h(x)-h(x') + 2\e_{h(y')} \max\{0, h(y') -  h(\lca(x,x'))\}. \label{eq:trans:1}
 \end{align}
If the maximum in \eqref{eq:trans:1} is $0$, then
\[ d_\e(y,y') = h(x)-h(x')\leq d_\e(x,x') .\]
If the maximum in~\eqref{eq:trans:1} equals $h(y') -  h(\lca(x,x'))$, then
\begin{multline}
d_\e(y,y')  = h(x)-h(x') + 2\e_{h(y')} [h(y')-h(\lca(x,x'))]
\\ \le h(x)-h(x') + 2\e_{h(x')} [h(x')-h(\lca(x,x'))] \label{eq:trans:2}
 = d_\e(x,x'),
\end{multline}
where in \eqref{eq:trans:2} we used the fact that the sequence $\{\e_n(n-a)\}_{n=0}^\infty$ is nondecreasing for all $a\ge 0$.
\end{proof}

\begin{lemma} \label{lem:horizontal}
Let $x,x',y\in B_\infty$ be such that $h(y)\le h(x)$. Then there exists $y'\in B_\infty$ which satisfies $
h(y')-h(x')=h(y)-h(x)$,
\begin{equation}
d_\e(y,y') \le d_\e(x,x'), \label{eq:hor-2}\end{equation}
and
\begin{equation}
d_\e(x,y)- 2d_\e(x,x') \le d_\e(y',x') \le d_\e(x,y)+ 2d_\e(x,x').  \label{eq:hor-1}
\end{equation}
\end{lemma}
\begin{proof}
Note that \eqref{eq:hor-1} follows from \eqref{eq:hor-2} by the triangle inequality.
Assume first that $h(x)\geq h(x')$.  In this case  choose $y'$ to be an ancestor of $y$ satisfying $h(y)-h(y')=h(x)-h(x')$. Then, $$d_\e(y,y')=h(y)-h(y')=h(x)-h(x')\leq d_\e(x,x').$$

We next assume that $h(x)< h(x')$.
If $h(\lca(x,x'))\neq h(\lca(x,y))$ then
 choose $y'$ to be an arbitrary descendant of $y$ such that $h(y')-h(y)=h(x')-h(x)$.
As before, we conclude that $d_\e(y,y')=h(y')-h(y)=h(x')-h(x)\leq d_\e(x,x')$.

%If $h(\lca(x,x')) > h(\lca(x,y))$, then
%we choose $y'$ to be an arbitrary descendant of $y$ such that $h(y')-h(y)=h(x')-h(x)$.
%As before $d_\e(y,y')=h(y')-h(y)=h(x')-h(x)\leq d_\e(x,x')$.

It remains to deal  with the case $h(x')>h(x)$ and $h(\lca(x,y))=h(\lca(x',x))$, which also implies that
$h(\lca(x',y))>h(\lca(x,y))$.
In this case, we choose $y'$ to be an arbitrary point on a branch containing both $\lca(x,y)$  and $x$, such that $h(y')-h(y)=h(x')-h(x)$. Then $\lca(y',y)=\lca(x,x')$, and therefore,
\begin{eqnarray*}
 d_\e(y,y')&=&h(y')-h(y)+2\e_{h(y)}[h(y)-h(\lca(y,y')]
 \\&=&  h(x')-h(x)+2\e_{h(y)}[h(y)-h(\lca(x,x'))]
\\&\le& h(x')-h(x)+2\e_{h(x)}[h(x)-h(\lca(x,x'))]\\&=&d_\e(x,x') ,
\end{eqnarray*}
proving~\eqref{eq:hor-2} in the last remaining case.
\end{proof}

\begin{lemma} \label{lem:descendant}
Let $x,x',y\in B_\infty$ be such that $y$ is a descendant of $x$. Then for any $y'\in B_\infty$ which is a descendant of
$x'$ and satisfying $h(y')-h(x')=h(y)-h(x)$, we have
\[ d_\e(y,y')\le d_\e(x,x')+2\e_{\min\{h(y'),h(y)\}}[h(y)-h(x)] \le d_\e(x,x')+2\e_{h(y)}[h(y)-h(x)+d_\e(x,x')].\]
\end{lemma}
\begin{proof}
Note that $h(\lca(y,y'))\ge h(\lca(x,x'))$. Assume first that $h(x')\ge h(x)$. Then,
\begin{align*}
d_\e(y,y') &= h(y') -h(y) +2 \e_{h(y)}[h(y)-h(\lca(y,y'))] \\
 &=h(x')-h(x) +2 \e_{h(y)}[h(x)-h(\lca(y,y'))]  +2 \e_{h(y)}[h(y)-h(x)] \\
 &\le h(x')-h(x) +2 \e_{h(y)}[h(x)-h(\lca(x,x')]  +2 \e_{h(y)}[h(y)-h(x)]
 \\
 &\le d_\e(x,x')+ 2 \e_{h(y)}[h(y)-h(x)] .
\end{align*}
When $h(x')< h(x)$, we similarly obtain  the bound:
\begin{eqnarray*}
d_\e(y,y') &=& h(y) -h(y') +2 \e_{h(y')}[h(y')-h(\lca(y,y'))] \\
&=& h(x)-h(x') +2 \e_{h(y')}[h(x')-h(\lca(y,y'))]  +2 \e_{h(y')}[h(y')-h(x')]\\
&\le& d_\e(x,x')+ 2 \e_{h(y')}[h(y)-h(x)].
\end{eqnarray*}
The last inequality in the statement of Lemma~\ref{lem:descendant} is proved by observing that
when $h(y')<h(y)$,
\[ \e_{h(y')}[h(y)-h(x)] =\e_{h(y')}[h(y')-h(x')]\le \e_{h(y)}[h(y)-h(x')] \le  \e_{h(y)}[h(y)-h(x)+d_\e(x,x')] .  \qedhere \]
\end{proof}

\newcommand{\constN}{7}

\begin{definition}\label{def:fork}
For $\delta\in (0,1)$ and $x,y,z,w\in B_\infty$, the quadruple $(x,y,z,w)$ is called a $\delta$-fork, if $$y\in
\Mid(x,z,\delta)\cap \Mid(x,w,\delta).$$
\end{definition}
$\delta$-forks in H-trees can be approximately classified using the approximate classification of midpoint configurations of Section~\ref{sec:midpoints}.
We have four types of midpoint configurations (recall Figure~\ref{fig:midpoints}):
\begin{itemize}
\item path-type; denoted (P) in what follows,
\item reverse path-type; denoted  (p)---$(x,y,z)$ is of type (p) iff $(z,y,x)$ is of type (P),
\item tent-type; denoted (T),
\item reverse tent-type; denoted  (t)---$(x,y,z)$ is of type (t) iff $(z,y,x)$ is of type (T).
\end{itemize}

Thus, there are
$\binom{5}{2}=10$ possible $\delta$-fork configurations in $(B_\infty,d_\e)$
(choose two out of the five symbols ``P",``p",``T",``t",``X", where ``X" means ``the same").
As we shall see, four of these possible configurations are impossible,
two of them have large contraction of the prongs of the forks, i.e., $d_\e(z,w)\ll d_\e(x,y)$,
which immediately implies large distortion, and the rest of the configurations are problematic
in the sense that they are not much distorted from the star $K_{1,3}$ (the metric $d$ on four points $p,q,r,s$ given by $d(p,q)=d(q,r)=d(q,s)=1$ and $d(p,s)=d(p,r)=d(r,s)=2$). The 10 possible $\delta$-fork configurations are summarized in Table~\ref{tab:types}.

\begin{table}[ht]
    \centering
        \begin{tabular}{ll}
            Midpoint configuration & Type \\ \hline
            (T$\|$T) & Type $I$\\
            (P$\|$P) & Type $II$\\
            (p$\|$T) & Type $III$\\
            (p$\|$t) & Type $IV$\\
            (p$\|$p) & prongs contracted\\
            (t$\|$t) & prongs contracted\\
            (P$\|$p) & impossible\\
            (P$\|$t) & possible only as approximate type $II$\\
            (P$\|$T) & impossible\\
            (t$\|$T) & impossible
        \end{tabular}
    \caption{{\em The ten possible fork configurations.}}
    \label{tab:types}
\end{table}

For future reference, we give names to the   four problematic configurations:
\begin{definition}\label{def:problematic}
For $\eta,\delta\in (0,1)$, a $\delta$-fork $(x,y,z,w)$ of $(B_\infty,d_\e)$ is called
\begin{itemize}
\item $\eta$-near Type $I$ (configuration (T$\|$T)) in Table~\ref{tab:types}), if both $(x,y,z)$ and $(x,y,w)$ are $\eta$-near tent-type configurations;
\item $\eta$-near Type $II$ (configuration (P$\|$P) in Table~\ref{tab:types}), if both $(x,y,z)$ and $(x,y,w)$ are $\eta$-near path-type configurations;
\item $\eta$-near Type $III$ (configuration (p$\|$T) in Table~\ref{tab:types}), if $(z,y,x)$ is $\eta$-near a path-type configuration and $(x,y,w)$
is $\eta$-near a tent-type configuration, or vice versa;
\item $\eta$-near Type $IV$ (configuration (p$\|$t) in Table~\ref{tab:types}), if $(z,y,x)$ is $\eta$-near a path-type configuration and $(w,y,x)$
is $\eta$-near a
%${\constN}\delta d_\e(x,y)$ near a
tent-type configuration, or vice versa.
\end{itemize}
\end{definition}

\begin{figure}[ht]
\begin{center}
\includegraphics[scale=0.8]{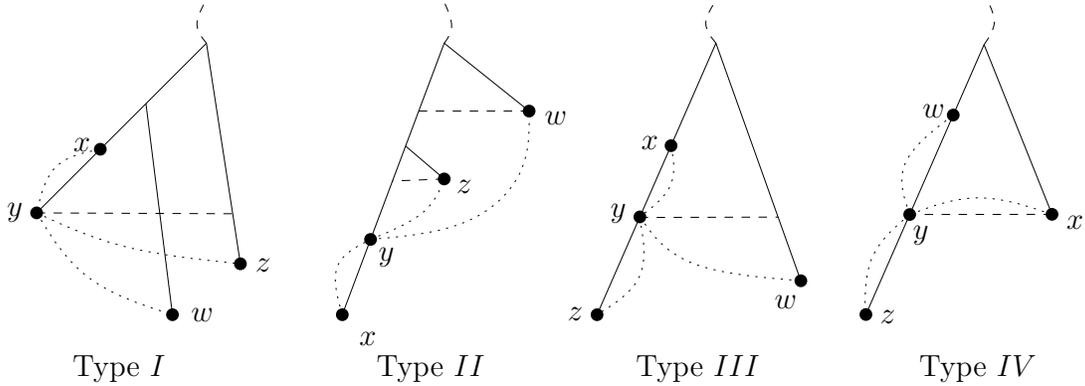}
\caption{{\em The four ``problematic'' types  of $\delta$-forks.}}
\label{fig:types}
\end{center}
\end{figure}

A schematic description of the four problematic configurations is contained in Figure~\ref{fig:types}.

The following lemma is the main result of this section.

\begin{lemma}\label{lem:fork classification0} Fix $\delta\in \left(0,\frac{1}{70}\right)$ and assume that $\e_n<\frac 14$  for all $n\in \mathbb N$.
If $(x,y,z,w)$ is a $\delta$-fork of $ (B_\infty,d_\e)$ then either it is $35\delta d_\e(x,y)$-near  one of the types $I$,  $II$,
$III$, $IV$, or we have $d_\e(z,w)\le 2(35\delta+ \e_{h_0}) d_\e(x,y),$
where $h_0=\min\{h(x),h(y),h(z),h(w)\}$.
%\begin{itemize}
%\item The $\delta$-fork $(x,y,z,w)$ is $15\delta d_\e(x,y)$ near  type $I$,  $II$,
%$III$, or  $IV$.
%\item Otherwise $d_\e(z,w)\le (83\delta+ \e_{h_0}) 2d_\e(y,z)$,
%where $h_0=\min\{h(x),h(y),h(z),h(w)\}$.
%\end{itemize}
\end{lemma}

\begin{remark}
{\em One can strengthen the statement of Lemma~\ref{lem:fork classification0} so that in the first case the fork $(x,y,z,w)$ is
$O(\delta d_\e(x,y))$ near another fork $(x',y',z',w')$ which is of (i.e. 0-near) one of the types $I$, $II$, $III$, $IV$.
This statement is more complicated to prove, and since we do not actually need it in what follows, we opted to use a weaker property which suffices for our purposes, yet simplifies (the already quite involved) proof.}
\end{remark}

\begin{figure}[ht]
\begin{center}
\includegraphics[scale=0.8]{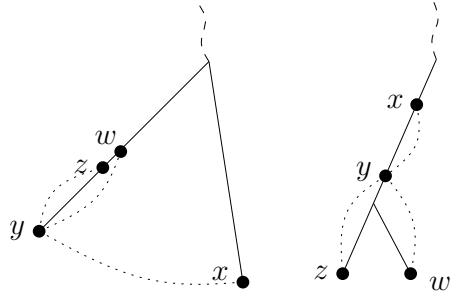}
 \caption{{\em The two configurations of
$\delta$-forks with large contraction of the prongs.}}
\label{fig:good types}
\end{center}
\end{figure}

The proof of Lemma~\ref{lem:fork classification0} proceeds by
checking that the cases marked in Table~\ref{tab:types} as
``impossible" or ``prongs contracted" are indeed so---see
Figure~\ref{fig:good types} for a schematic description of the
latter case.

We begin with the $(p\|p)$ configuration.
\begin{lemma} \label{lem:opath,opath}
Let $(x,y,z,w)$ be a $\delta$-fork of $(B_\infty,d_\e)$ and
assume that both $(z,y,x)$, and $(w,y,x)$ are $\eta d_\e(x,y)$-near path-type configurations.
Then, assuming that $\max\{\delta, \eta\}<1/8$ and $\e_n<\frac 14$  for all $n$,
we have
\[ d_\e(z,w) \le \left(9\eta+6\delta+2\e_{h(y)}\right)d_\e(x,y). \]
\end{lemma}
\begin{proof}
Let $(z',y',x')$ be a path-type configuration that is $\eta d_\e(x,y)$-near $(z,y,x)$,
and let $(w^{\prime\prime},y^{\prime\prime},x^{\prime\prime})$ be a path-type configuration that is
$\eta d_\e(x,y)$-near $(w,y,x)$.
Without loss of generality assume that $h(y^{\prime\prime})\ge h(y')$.
Let $w'$ be the descendant of $y'$ satisfying $h(w')-h(y')=h(w^{\prime\prime})-h(y^{\prime\prime})$ such that $w'$ is either an
ancestor or an arbitrary descendant of $z'$. Note that $h(w')\ge h(y)$. Indeed,
\begin{multline*}
h(w')=h(y')+h(w'')-h(y'')=h(y')+d_\e(w'',y'')\ge h(y)-|h(y)-h(y')|+d_\e(w,y)-2\eta d_\e(x,y)
\\ \ge h(y)-d_\e(y,y')+\frac{1-\delta}{1+\delta}d_\e(x,y)-2\eta d_\e(x,y)\ge h(y)+\left(\frac{1-\delta}{1+\delta}-3\eta\right)d_\e(x,y)\ge h(y).
\end{multline*}
By Lemma~\ref{lem:descendant},
 \begin{multline}\label{eq:hy}
 d_\e(w',w^{\prime\prime})\le
 d_\e(y',y^{\prime\prime})+ 2\e_{h(w')} \left[h\left(w^{\prime\prime}\right)-h\left(y^{\prime\prime}\right)\right]
 \leq d_\e(y',y^{\prime\prime})+ 2\e_{h(w')} d_\e(y^{\prime\prime},w^{\prime\prime})\\ \le 2\eta d_\e(x,y) +2\e_{h(y)}\left(2\eta+\frac{1+\delta}{1-\delta}\right)d_\e(x,y) .
 \end{multline}
%We now estimate $d_\e(w',z')$.
%Assuming that $w'$ is an ancestor of $z'$ (the opposite case is analogous),
Observe that
\begin{multline}\label{eq:primes}
\left(\frac{1-\delta}{1+\delta}-2\eta\right)d_\e(x,y)\le  d_\e(z,y)-2\eta d_\e(x,y)\le d_\e(z',y')\\\le d_\e(z,y)+2\eta d_\e(x,y) \le \left(\frac{1+\delta}{1-\delta}+2\eta\right) d_\e(x,y),\end{multline}
Since $ d_\e(w',y')=d_\e(w^{\prime\prime},y^{\prime\prime})$, we obtain similarly the bounds:
\begin{equation}\label{eq:double primes}
\left(\frac{1-\delta}{1+\delta}-2\eta\right)d_\e(x,y)\le  d_\e(w',y') \le d_\e(z,y)+2\eta d_\e(x,y) \le \left(\frac{1+\delta}{1-\delta}+2\eta\right) d_\e(x,y)
. \end{equation}
Hence
\begin{equation}\label{eq:absolute}
 d_\e(z',w')=\left|d_\e(y',z')-d_\e(y',w')\right| \stackrel{\eqref{eq:primes}\wedge\eqref{eq:double primes}}{\le} \left(\frac{4\delta}{1-\delta^2}+4\eta\right)d_\e(x,y).
 \end{equation}
So, in conclusion,
\begin{multline*} d_\e(z,w) \le
d_\e(z,z')+d_\e(w,w^{\prime\prime})+d_\e(w^{\prime\prime},w') + d_\e(z',w') \\
\stackrel{\eqref{eq:hy}\wedge\eqref{eq:absolute}}{\le} \left(8\eta+  \frac{4\delta}{1-\delta^2}+ 2\left(2\eta+\frac{1+\delta}{1-\delta}\right)\e_{h(y)}\right) d_\e(x,y)\le \left(9\eta+6\delta+2\e_{h(y)}\right)d_\e(x,y)\qedhere.
\end{multline*}
\end{proof}

We next consider the (t$\|$t) configuration.
\begin{lemma} \label{lem:otent,otent}
Let $(x,y,z,w)$ be a $\delta$-fork of $(B_\infty,d_\e)$. Assume that both $(z,y,x)$ and $(w,y,x)$ are $\eta d_\e(x,y)$-near tent-type configurations. Then, assuming that $\max\{\delta,\eta\}<1/4$, we have
\[ d_\e(z,w) \le (8\eta+5\delta) d_\e(x,y). \]
\end{lemma}
\begin{proof}
Let $(z',y',x')$ be a tent-type configuration that is $\eta d_\e(x,y)$-near $(z,y,x)$, and let $(w^{\prime\prime},y^{\prime\prime},x^{\prime\prime})$ be a tent-type configuration that is
$\eta d_\e(x,y)$-near $(w,y,x)$.
Assume without loss of generality that  $h(y^{\prime\prime})-h(w^{\prime\prime}) \ge h(y')-h(z')$. Let $\tilde w$ be a point on the path between
$w^{\prime\prime}$ and $y^{\prime\prime}$ such that $h(y^{\prime\prime})-h(\tilde w) = h(y')-h(z')$. Then,
\begin{multline}\label{eq;tilde}
 d_\e(w^{\prime\prime}, \tilde w) = h(y^{\prime\prime})-h(w^{\prime\prime})-(h(y') - h(z'))=d_\e(y'',w'')-d_\e(y',z')\\\le d_\e(y,w)-d_\e(y,z)+4\eta d_\e(x,y) \le \left(\frac{1+\delta}{1-\delta}-\frac{1-\delta}{1+\delta}+4\eta\right)d_\e(x,y).\end{multline}
 By Lemma~\ref{lem:ancestor} we have $d_\e(\tilde w,z')\le d_\e(y',y^{\prime\prime})
\le 2\eta  d_\e(x,y)$. Hence we conclude that
\[ d_\e(y,z)\le d_\e(z,z')+d_\e(\tilde w,z')+d_\e(\tilde w,w^{\prime\prime})+ d_\e(w^{\prime\prime},w) \stackrel{\eqref{eq;tilde}}{\le} \left(\frac{4\delta}{1-\delta^2}+8\eta\right) d_\e(x,y)  .  \qedhere \]
\end{proof}

\begin{lemma} \label{lem:path type =>path type}
Let $(x,y,z,w)$ be a $\delta$-fork of $B_\infty$. Assume that $(x,y,z)$ is $\eta d_\e(x,y)$-near a path-type configuration.
Assume also that $\delta<1/30$, $\eta<1/10$, and  $\e_n<1/4$ for all $n$.
Then $(x,y,w)$ is $(2\eta+21\delta) d_\e(x,y)$-near a path-type configuration,
i.e., $(x,y,z,w)$ is $(2\eta+21\delta) d_\e(x,y)$-near a type $II$ configuration.
\end{lemma}
\begin{proof}
Let $(x',y',z')$ be a path-type configuration which is $\eta d_\e(x,y)$-near $(x,y,z)$.
By Theorem~\ref{thm:midpoints}, either $(x,y,w)$ or $(w,y,x)$ must be $3\delta d_\e(x,z)\le \frac{6}{1-\delta} d_\e(x,y)\le 7\delta d_\e(x,y)$-near either a path-type configuration or a
tent-type configuration.

Suppose first that $(x,y,w)$ is $7\delta d_\e(x,y)$-near a tent-type configuration  $(x^{\prime\prime},y^{\prime\prime},w^{\prime\prime})$.
In this case,
$x^{\prime\prime}$ is an ancestor of $y^{\prime\prime}$ and
$h(y^{\prime\prime})-h(x^{\prime\prime})=d_\e(x^{\prime\prime},y^{\prime\prime})\ge (1-14\delta) d_\e(x,y)$.
At the same time, $y'$ is an ancestor of $x'$ and
$h(x')-h(y')=d_\e(x',y')\ge (1-2\eta) d_\e(x,y)$.
So,
\[
2(\eta +7\delta) d_\e(x,y)\ge d_\e(y'',y')+d_\e(x',x'') \ge h(y^{\prime\prime})-h(x^{\prime\prime}) +h(x')-h(y')
\ge 2(1-\eta-7\delta) d_\e(x,y),
\]
which is a contradiction since $\eta+7\delta<1/2$.

Next suppose that $(w,y,x)$ is $7\delta d_\e(x,y)$-near a path-type configuration $(w^{\prime\prime},y^{\prime\prime},x^{\prime\prime})$.
Then $|h(x')-h(x^{\prime\prime})|\le d_\e(x',x^{\prime\prime})\le (\eta+7\delta)d_\e(x,y)$.
So,
\begin{multline*}
 (\eta+7\delta)d_\e(x,y)  \ge d_\e(y',y^{\prime\prime})
 \ge h(y^{\prime\prime})-h(y')\\ = (h(y^{\prime\prime})-h(x^{\prime\prime})) +(h(x^{\prime\prime})-h(x')) +h(x')-h(y')
> 0 - (\eta+7\delta)d_\e(x,y) + (1-2\eta) d_\e(x,y),
\end{multline*}
which is a contradiction

Lastly, suppose that $(w,y,x)$ is $7\delta d_\e(x,y)$-near a tent-type configuration $(w^{\prime\prime},y^{\prime\prime},x^{\prime\prime})$. Note that $\left|h(y')-h(y'')\right|\le d_\e(y',y'')\le (\eta+7\delta)d_\e(x,y)$. So, $h(y')\ge h(y'')-(\eta+7\delta)d_\e(x,y)$. Also,
$$h(y'')-h(w'')=d_\e(y'',w'')\ge d_\e(y,w)-14\delta d_\e(x,y)\ge \left(\frac{1-\delta}{1+\delta}-14\delta\right)d_\e(x,y)\ge (\eta+7\delta)d_\e(x,y).$$
Consider the point $\bar w$ defined as the ancestor of $y'$ at distance $h(y'')-h(w'')-(\eta+7\delta)d_\e(x,y)$ from $y'$. Let also $w'''$ be the ancestor of $y''$ at distance $h(y'')-h(w'')-(\eta+7\delta)d_\e(x,y)$ from $y''$. By
Lemma~\ref{lem:ancestor}, we have $d_\e(\bar w,w''')\le d_\e(y',y'')\le (\eta+7\delta)d_\e(x,y)$. Therefore,
\[ d_\e(\bar w,w)\le d_\e(\bar w,w''')+d_\e(w''',w'')+d_\e(w'',w) \le (2\eta+21\delta)d_\e(x,y). \]
Hence $(x,y,w)$ is $(2\eta+21\delta)d_\e(x,y)$-near the path-type configuration $(x',y',\bar w)$.
\end{proof}

\begin{lemma} \label{lem:tent => not opposite tent}
Let $(x,y,z,w)$ be a $\delta$-fork of $B_\infty$. Assume that $(x,y,z)$ is $\eta d_\e(x,y)$-near a tent-type configuration. Assume also that $\eta<1/10$ and  $\e_n<1/4$ for all $n$.
Then  $(w,y,x)$
cannot be $\eta d_\e(x,y)$-near a tent-type configuration.
\end{lemma}
\begin{proof}
Let $(x',y',z')$ be a tent type configuration that is $\eta d_\e(x,y)$-near $(x,y,z)$. Suppose for contradiction
that there exists a tent type configuration $(w^{\prime\prime},y^{\prime\prime},x^{\prime\prime})$ that is $\eta d_\e(x,y)$-near $(w,y,x)$. Note that $h(y'')\ge h(y')-d_\e(y',y'')\ge h(y')-2\eta d_\e(x,y)$ and $h(y')-h(x')\ge (1-2\eta)d_\e(x,y)>2\eta d_\e(x,y)$.  Let $x^*$ be the ancestor of $y'$ at distance $h(y')-h(x')-2\eta d_\e(x,y)$ from $y'$, and let $\tilde x$ be the ancestor of $y''$ at distance $h(y')-h(x')-2\eta d_\e(x,y)$ from $y''$. An application of Lemma~\ref{lem:ancestor} yields the estimate $d_\e(\tilde x,x^*)\le d_\e(y',y'')\le 2\eta d_\e(x,y)$. But, since $h(x^{\prime\prime})\ge h(y^{\prime\prime})$, we also know that $d_\e(\tilde x,x^{\prime\prime})\ge h(y^{\prime\prime})-h(\tilde x)=d_\e(y',x')-2\eta d_\e(x,y)$. Hence,
\begin{equation*}
2\eta d_\e(x,y)\ge d_\e(\tilde x,x^*)\ge d_\e(\tilde x,x'')-d_\e(x^*,x')-d_\e(x',x'')\ge d_\e(x',y')-6\eta d_\e(x,y)\ge (1-8\eta)d_\e(x,y),
\end{equation*}
which is a contradiction, since $\eta<1/10$.
\end{proof}

\begin{proof}[Proof of Lemma~\ref{lem:fork classification0}]
Since $(x,y,z,w)$ is a $\delta$-fork, by Theorem~\ref{thm:midpoints}, both $(x,y,z)$ and $(x,y,w)$ are
$7\delta d_\e(x,y)$-near a tent-type configuration, a path-type configuration, or the corresponding reverse configurations. We have 10 possible combinations
of these pairs, as appearing in Table~\ref{tab:types}.
By applying Lemmas~\ref{lem:path type =>path type} and~\ref{lem:tent => not opposite tent} with $\eta=7\delta$,
we rule out three
of these configurations, and a fourth configuration is possible but only as $35\delta d_\e(x,y)$-near a type $II$ configuration.

We are left with six possible configurations.
By applying Lemmas~\ref{lem:opath,opath} and~\ref{lem:otent,otent} with $\eta=7\delta$ we conclude that in two of those
configurations we have $d_\e(w,z)\le (69\delta+2\e_{h_0})d_\e(x,y)$, and the rest are configurations that are $7\delta d_\e(x,y)$-near one of the types $I$--$IV$.
\end{proof}

\subsubsection{Classification of approximate 3-paths}
\label{sec:3-path}

We start with the following natural notion:
\begin{definition}\label{def:3 path}
For $x_0,x_1,x_2,x_3\in B_\infty$ the quadruple $(x_0,x_1,x_2,x_3)$ is called a $(1+\delta)$-approximate $P_3$ if there exists $L>0$ such that for every $0\le i\le j\le 3$ we have
\[ (j-i)L\le d_\e(x_i,x_j)\le (1+\delta) (j-i)L. \]
Note that in this case $x_1\in \Mid(x_0,x_2,\delta)$ and
$x_2\in \Mid(x_1,x_3,\delta)$.
\end{definition}
As in the case of $\delta$-forks, there are 10 possible concatenations of two midpoints configurations (path-type or tent-type):
P-P, P-p, P-T, P-t, p-P, p-T, p-t, T-T, T-t, t-T (the midpoint configurations
p-p, P-p, t-p, T-p, p-P, t-P, T-P, t-t, T-t, t-T are respectively such concatenations with the order of $x_0,x_1,x_2,x_3$ reversed).
We will rule out some of these possibilities, and obtain some stronger properties for the rest. See Table~\ref{tab:3-paths-types}.

\begin{table}[ht]
    \centering
        \begin{tabular}{lll}
            Midpoint configuration & Reverse configuration & Type \\ \hline
            (P-P) & (p-p) &type $A$\\
            (P-p) & (P-p) &impossible\\
            (P-T) & (t-p) &impossible \\
            (P-t) & (T-p) &type $B$\\
            (p-P) & (p-P) &impossible\\
            (p-T) & (t-P) &type $C$\\
            (p-t) & (T-P) &impossible\\
            (T-T) & (t-t) &possible only as type $C$\\
            (T-t) & (T-t) &impossible\\
            (t-T) & (t-T) &impossible
        \end{tabular}
    \caption{{\em The possible configurations of 3 paths.}}
    \label{tab:3-paths-types}
\end{table}

As in the case of $\delta$-forks, it will be beneficial to give names to three special types approximate $3$-paths:

\begin{definition}\label{def:ABC}
For $x_0,x_1,x_2,x_3\in B_\infty$ and $\eta>0$, a quadruple $(x_0,x_1,x_2,x_3)$ is called:
\begin{itemize}
\item $\eta$-near a type $A$ configuration if both $(x_0,x_1,x_2)$ and $(x_1,x_2,x_3)$ are $\eta$-near
path-type configurations,
\item $\eta$-near a type $B$ configuration  if  $(x_0,x_1,x_2)$ is $\eta$-near a path-type configuration, and
 $(x_3,x_2,x_1)$ is $\eta$-near tent-type configuration,
\item $\eta$ near type $C$ configuration  if  $(x_2,x_1,x_0)$ is $\eta$-near a path-type configuration, and
 $(x_1,x_2,x_3)$ is $\eta$-near a tent-type configuration.
\end{itemize}
See also Figure~\ref{fig:P3-types}.
\begin{figure}[ht]
\begin{center}
\includegraphics[scale=0.8]{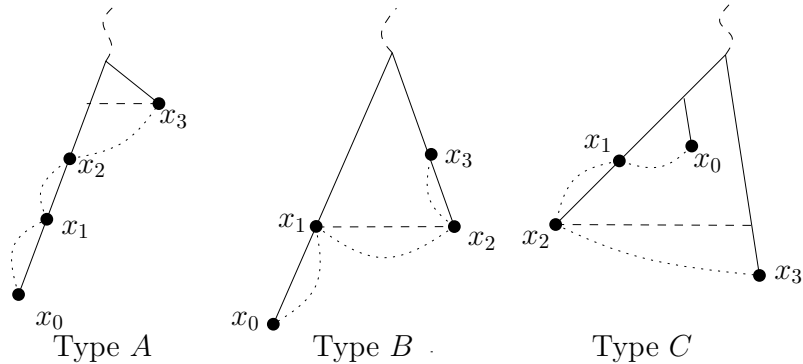}
\caption{{\em The three possible types of approximate 3-paths.}}
\label{fig:P3-types}
\end{center}
\end{figure}
\end{definition}

The following lemma is the main result of this subsection.
\begin{lemma} \label{lem:3-path classification}
Assume that $\e_n<\frac14$ for all $n$ and fix $\delta<1/200$. Assume that $(x_0,x_1,x_2,x_3)$ is a $(1+\delta)$-approximate $P_3$. Then either $(x_0,x_1,x_2,x_3)$ or $(x_3,x_2,x_1,x_0)$ is $35\delta d_\e(x_0,x_1)$-near a configuration of type $A$, $B$ or $C$.
\end{lemma}

The proof of Lemma~\ref{lem:3-path classification} is again a case analysis that examines all 10 possible ways (up to symmetry)
to concatenate two midpoint configurations. The proof is divided into a few lemmas according to the cases, and is completed at the end of this subsection.

\begin{lemma} \label{lem:3-path-1} Assume that $\e_n<\frac14$ for all $n$ and that $(x_0,x_1,x_2,x_3)$ is a $(1+\delta)$-approximate $P_3$ such that
$(x_0,x_1,x_2)$ is $\eta d_\e(x_0,x_1)$-near a path-type configuration. If $\max\{\delta,\eta\}<1/200$
then either $(x_1,x_2,x_3)$ is $7\delta d_\e(x_0,x_1)$-near a path-type configuration (type $A$), or $(x_3,x_2,x_1)$
is $7\delta d_\e(x_0,x_1)$-near a tent-type configuration (type $B$).
\end{lemma}
\begin{proof}
Due to Theorem~\ref{thm:midpoints} we only need to rule out the possibility that $(x_3,x_2,x_1)$ is $7\delta d_\e(x_1,x_2)$-near
a path-type configuration, or that $(x_1,x_2,x_3)$ is $7\delta d_\e(x_1,x_2)$-near a tent-type configuration.
Let $(x'_0,x'_1,x'_2)$ be a path-type configuration that is $\eta d_\e(x_0,x_1)$-near $(x_0,x_1,x_2)$.

Suppose first that $(x_3,x_2,x_1)$ is $7\delta d_\e(x_1,x_2)$-near the path-type configuration $(x_3^{\prime\prime},x^{\prime\prime}_2,x^{\prime\prime}_1)$. Since $h(x'_1)\ge h(x'_2)$ and $h(x^{\prime\prime}_2)\ge h(x^{\prime\prime}_1)$ we have,
\begin{multline}\label{eq:12'}
|h(x_1)-h(x_2)|\le |h(x_1)-h(x_1')|+\left(h(x_1')-h(x_2')\right)+|h(x_2')-h(x_2)|\\ \le d_\e(x_1,x_1')+\left(h(x_1')-h(x_2')\right)+d_\e(x_2',x_2)\le 2\eta d_\e(x_0,x_1)+\left(h(x_1')-h(x_2')\right),
\end{multline}
and similarly,
\begin{multline}\label{eq:12''}
|h(x_1)-h(x_2)|\le |h(x_1)-h(x_1'')|+\left(h(x_2'')-h(x_1'')\right)+|h(x_2'')-h(x_2)|\\ \le d_\e(x_1,x_1'')+\left(h(x_2'')-h(x_1'')\right)+d_\e(x_2',x_2)
\le 14\delta d_\e(x_0,x_1)+\left(h(x_2'')-h(x_1'')\right).
\end{multline}
By summing~\eqref{eq:12'} and~\eqref{eq:12''} we obtain the bound
$$
2|h(x_1)-h(x_2)|\le \left(2\eta+14\delta\right)d_\e(x_0,x_1)+d_\e(x_1',x_1'')+d_\e(x_2',x_2'')\le \left(4\eta+28\delta\right)d_\e(x_0,x_1).
$$
Thus
\begin{equation}\label{eq:height bound12}
|h(x_1) - h(x_2)| \le (2\eta+14\delta)d_\e(x_0,x_1).
\end{equation}

Since $x'_0$ is a descendant of $x'_1$,
\begin{multline}\label{eq:01eps}
| h(x_0)- h(x_1)  - d_\e(x_0,x_1)| \le |h(x'_0) - h(x'_1) -d_\e(x_0,x_1)| + 2\eta d_\e(x_0,x_1) \\
= |d_\e(x'_0,x'_1) -d_\e(x_0,x_1)|+ 2\delta d_\e(x_0,x_1) \le 4\eta d_\e(x_0,x_1) .
\end{multline}
Similarly, since $x^{\prime\prime}_3$ is a descendant of $x^{\prime\prime}_2$,
\begin{multline}\label{eq:23eps}
|h(x_3)- h(x_2)  - d_\e(x_0,x_1)| \le |h(x^{\prime\prime}_3) - h(x^{\prime\prime}_2) -d_\e(x_0,x_1)|
+ 14\delta d_\e(x_1,x_2) \\
= |d_\e(x^{\prime\prime}_3,x^{\prime\prime}_2) -d_\e(x_0,x_1)|+ 14\delta d_\e(x_0,x_1)
\le 28\delta d_\e(x_0,x_1) .
\end{multline}
Hence,
\begin{eqnarray}\label{eq:49-7}
&&\!\!\!\!\!\!\!\!\!\!\!\!\!\!\!\!\nonumber |h(x_3'')-h(x_0')|\le |h(x_3'')-h(x_3)|+ |h(x_3)- h(x_2)  - d_\e(x_0,x_1)|+|h(x_2)-h(x_1)|\\&&\quad\quad\quad\quad\quad\quad
+ | h(x_0)- h(x_1)  - d_\e(x_0,x_1)|+|h(x_0)-h(x_0')|\nonumber\\
&&\stackrel{\eqref{eq:height bound12}\wedge\eqref{eq:01eps}\wedge\eqref{eq:23eps}}{\le} d_\e(x_3'',x_3)+28\delta d_\e(x_0,x_1)+ (2\eta+14\delta)d_\e(x_0,x_1)+4\eta d_\e(x_0,x_1)+d_\e(x_0,x_0')\nonumber\\
&&\quad\quad\, \le (49\delta+7\eta)d_\e(x_0,x_1).
\end{eqnarray}

We record for future reference the following consequence of~\eqref{eq:height bound12} and~\eqref{eq:49-7}:
\begin{eqnarray}\label{eq:mins}
&&\nonumber\!\!\!\!\!\!\!\!\!\!\!\!\!\!\!\!\min\{h(x'_0),h(x_3'')\}-\min\{h(x'_1),h(x_2'')\}\le \max\left\{h(x_0')-h(x_1'),h(x_0')-h(x_2'')\right\}\\
&\stackrel{\eqref{eq:49-7}}{\le}&\nonumber\max\left\{d_\e(x_0',x_1'),h(x_3'')-h(x_2'')+(49\delta+7\eta)d_\e(x_0,x_1)\right\}\\
&\le& \max\left\{(1+2\eta)d_\e(x_0,x_1),d_\e(x_3'',x_2'')+(49\delta+7\eta)d_\e(x_0,x_1)\right\}\nonumber\\
&\le& (1+64\delta+7\eta)d_\e(x_0,x_1).
\end{eqnarray}

We next claim that
\begin{equation}\label{eq:lca12}
\lca(x_0',x_3'')=\lca(x_1',x_2'').
\end{equation}
Indeed, since $x_1'$ is an ancestor of $x_0'$ and $x_2''$ is an ancestor of $x_3''$, if  $\lca(x_0',x_3'')\neq\lca(x_1',x_2'')$ then either $x_1'$ is a descendant of $x_2''$, or $x_2''$ is a descendant of $x_1'$. If $x_1'$ is a descendant of $x_2''$ then
\begin{multline*}
(\eta+7\delta)d_\e(x_0,x_1)\ge d_\e(x_1',x_1'')\ge d_\e(x_2'',x_1')\ge d_\e(x_2,x_1)-(\eta+7\delta)d_\e(x_0,x_1)\\\ge \frac{1}{1+\delta}d_\e(x_0,x_1)-(\eta+7\delta)d_\e(x_0,x_1),
\end{multline*}
which is a contradiction since $\delta,\eta<1/200$. Similarly, if $x_2''$ is a descendant of $x_1'$ then
\begin{equation*}
(\eta+7\delta)d_\e(x_0,x_1)\ge d_\e(x_2',x_2'')\ge d_\e(x_2'',x_1')\ge \frac{1}{1+\delta}d_\e(x_0,x_1)-(\eta+7\delta)d_\e(x_0,x_1),
\end{equation*}
arriving once more at a contradiction. This proves~\eqref{eq:lca12}.

Now,
\begin{align}
\nonumber\frac{3}{1+\delta}d_\e(x_0,x_1)&\le d_\e(x_0,x_3)\\\nonumber &\le d_\e(x^{\prime\prime}_3,x'_0) +(\eta+7\delta)d_\e(x_0,x_1)\\ \nonumber
&\stackrel{\eqref{eq:49-7}}{\le} 2\e_{\min\{h(x_0'),h(x_3'')\}} \left[\min\{h(x'_0),h(x_3'')\}- h(\lca(x'_0,x^{\prime\prime}_3))\right] +(8\eta+56\delta) d_\e(x_0,x_1)\\ \nonumber
& \stackrel{\eqref{eq:lca12}}{=} 2\e_{\min\{h(x_0'),h(x_3'')\}} \left[\min\{h(x'_1),h(x_2'')\}- h(\lca(x'_1,x^{\prime\prime}_2))\right]+(8\eta+56\delta) d_\e(x_0,x_1) \\\nonumber&\quad\quad+2\e_{\min\{h(x_0'),h(x_3'')\}} \left[\min\{h(x'_0),h(x_3'')\}-\min\{h(x'_1),h(x_2'')\}\right]\\\nonumber
&\stackrel{\eqref{eq:mins}}{\le} 2\e_{\min\{h(x_1'),h(x_2'')\}} \left[\min\{h(x'_1),h(x_2'')\}- h(\lca(x'_1,x^{\prime\prime}_2))\right]\\&\quad\quad+\left(8\eta+56\delta+\frac{1+64\delta+7\eta}{2}\right)d_\e(x_0,x_1)
\label{eq:eps<1/4}\\
& \le \nonumber  d_\e(x'_1,x^{\prime\prime}_2) +\left(\frac12+88\delta+12\eta\right)d_\e(x_0,x_1)\\
&\le \left(\frac32+96\delta+13\eta\right)d_\e(x_0,x_1)\label{eq:1/200},
\end{align}
where in~\eqref{eq:eps<1/4} we used $\min\{h(x_0'),h(x_3'')\}\ge \min\{h(x_1'),h(x_2'')\}$ and $\e_{\min\{h(x_0'),h(x_3'')\}}<1/4$. Since   $\max\{\eta,\delta\}<1/200$, the bound~\eqref{eq:1/200} is a contradiction.

Next suppose that $(x_1,x_2,x_3)$ is $7\delta d_\e(x_0,x_1)$-near a tent-type configuration  $(x^{\prime\prime}_1,x^{\prime\prime}_2,x^{\prime\prime}_3)$. Since $h(x_2')\le h(x_1')$ and $h(x_2'')\ge h(x_1'')$, we have
\begin{eqnarray}\label{eq:h12}
|h(x_1)-h(x_2)|&\le&\nonumber |h(x_1)-h(x_1')|+\left(h(x_1')-h(x_2')\right)+|h(x_2')-h(x_2)|\\ \nonumber&\le&
d_\e(x_1,x_1')+\left(h(x_1')-h(x_2')\right)+\left(h(x_2'')-h(x_1'')\right)+d_\e(x_2',x_2)\\
&\le&\nonumber d_\e(x_1,x_1')+d_\e(x_1',x_1'')+d_\e(x_2',x_2'')+d_\e(x_2',x_2)\\
&\le& (4\eta+14\delta)d_\e(x_0,x_1).
\end{eqnarray}
On the other hand, $x^{\prime\prime}_1$ is an ancestor of $x^{\prime\prime}_2$, and therefore we have
\begin{multline}\label{eq:for next lemma}
\left(\frac{1}{1+\delta}-14\delta\right)d_\e(x_0,x_1)\le d_\e(x_1'',x_2'')=h(x_2'')-h(x_1'')\\\le |h(x_1)-h(x_2)|+14\delta d_\e(x_0,x_1)\stackrel{\eqref{eq:h12}}{\le}(4\eta+28\delta)d_\e(x_0,x_1),
\end{multline}
which is a contradiction since $\max\{\eta,\delta\}<1/200$.
\end{proof}

\begin{lemma} \label{lem:3-path-2}
Assume that $\e_n<\frac14$ for all $n$ and that $(x_0,x_1,x_2,x_3)$ is a $(1+\delta)$-approximate $P_3$ such that
$(x_2,x_1,x_0)$ is $\eta d_\e(x_0,x_1)$-near a path-type configuration. If $\max\{\delta,\eta\}<1/200$
then either $(x_3,x_2,x_1)$ is $7\delta d_\e(x_0,x_1)$-near a path-type configuration (reverse type $A$), or $(x_1,x_2,x_3)$
is $7\delta d_\e(x_0,x_1)$-near a tent-type configuration (type $C$).
\end{lemma}
\begin{proof}
Let $(x'_2,x'_1,x'_0)$ be in path-type configuration that is $\eta d_\e(x_0,x_1)$-near $(x_2,x_1,x_0)$. First, assume for contradiction that $(x_3,x_2,x_1)$ is $7\delta d_\e(x_1,x_2)$-near a  tent-type configuration $(x^{\prime\prime}_3,x^{\prime\prime}_2,x^{\prime\prime}_1)$.
Then $h(x^{\prime\prime}_1) \ge h(x^{\prime\prime}_2)$, where as
$h(x'_2)-h(x'_1)=d_\e(x'_2,x'_1)$. Arguing as in~\eqref{eq:h12}, it follows that $|h(x_1)-h(x_2)|\le (2\eta+28\delta)d_\e(x_0,x_1)$, and we arrive at a contradiction by arguing similarly to~\eqref{eq:for next lemma}.

Next, assume for contradiction that $(x_1,x_2,x_3)$ is $7\delta d_\e(x_1,x_2)$-near a path-type configuration $(x^{\prime\prime}_1,x^{\prime\prime}_2,x^{\prime\prime}_3)$. Then
$h(x^{\prime\prime}_1)- h(x^{\prime\prime}_2)=d_\e(x^{\prime\prime}_1, x^{\prime\prime}_2)$, whereas
$h(x^{\prime}_2)- h(x^{\prime}_1)=d_\e(x^{\prime}_1, x^{\prime}_2)$. By summing these two identities, we arrive at a contradiction as follows:
\begin{multline*}
\left(\frac{2}{1+\delta}-2\eta-14\delta\right)d_\e(x_0,x_1)\le d_\e(x^{\prime}_1, x^{\prime}_2)+d_\e(x^{\prime\prime}_1, x^{\prime\prime}_2)= \left(h(x_2')-h(x_2'')\right)+\left(h(x_1'')-h(x_1')\right)\\\le d_\e(x_2',x_2'')+d_\e(x_1',x_1'')\le (2\eta+14\delta)d_\e(x_0,x_1).\qedhere
\end{multline*}
\end{proof}

\begin{lemma} \label{lem:3-path-3}
Assume that $\e_n<\frac14$ for all $n$ and that $(x_0,x_1,x_2,x_3)$
is a $(1+\delta)$-approximate $P_3$ such that $(x_0,x_1,x_2)$ is
$\eta d_\e(x_0,x_1)$-near a tent-type configuration. If
$\max\{\delta,\eta\}<1/200$ then either $(x_2,x_1,x_0)$ is
$(14\delta+3\eta)d_\e(x_0,x_1)$-near a path-type configuration and $(x_1,x_2,x_3)$ is $7\delta d_\e(x_1,x_2)$-near a tent-type
configuration (type
$C$), or
 $(x_3,x_2,x_1)$ is $7\delta d_\e(x_0,x_1)$-near a path-type configuration (reverse type $B$).
\end{lemma}
\begin{proof}
Let $(x'_0,x'_1,x'_2)$ be a tent-type configuration that is $\eta
d_\e(x_0,x_1)$-near $(x_0,x_1,x_2)$. First, suppose that
$(x_1,x_2,x_3)$ is $7\delta d_\e(x_1,x_2)$-near a tent-type
configuration
$(x^{\prime\prime}_1,x^{\prime\prime}_2,x^{\prime\prime}_3)$. Note
that $|h(x_1')-h(x_1'')|\le d_\e(x_1',x_1'')\le
(\eta+7\delta)d_\e(x_0,x_1)$. So, let $x_0''$ be an ancestor of
$x_1''$ at distance $h(x_1')-h(x_0')-(\eta+7\delta)d_\e(x_0,x_1)\in
[0,h(x_1'')]$ from $x_1''$, and let $x_0^*$ be an ancestor of $x_1'$
at distance $h(x_1')-h(x_0')-(\eta+7\delta)d_\e(x_0,x_1)$ from
$x_1'$. Then $h(x_1')-h(x_0^*)=h(x_1'')-h(x_0'')$ and
$d_\e(x_0^*,x_0')\le (\eta+7\delta)d_\e(x_0,x_1)$. By
Lemma~\ref{lem:ancestor},
\begin{multline*}
d_\e(x_0,x_0'')-(2\eta+7\delta)d_\e(x_0,x_1\le
d_\e(x_0,x_0'')-d_\e(x_0^*,x_0')-d_\e(x_0',x_0)\le
d_\e(x_0^*,x_0'')\\\le d_\e(x_1',x_1'')\le
(\eta+7\delta)d_\e(x_0,x_1).
\end{multline*}
Hence $(x^{\prime\prime}_2,x^{\prime\prime}_1,x^{\prime\prime}_0)$
is a path-type configuration that is $(14\delta+3\eta)d_\e(x_0,x_1)
$-near $(x_2,x_1,x_0)$.

Next assume for contradiction that $(x_3,x_2,x_1)$ is $7\delta
d_\e(x_1,x_2)$ \-near a tent-type configuration
$(x^{\prime\prime}_3,x^{\prime\prime}_2,x^{\prime\prime}_1)$. Then
\begin{equation}\label{eq:height23''}
\left(1-15\delta\right)d_\e(x_0,x_1)\le
\left(\frac{1}{1+\delta}-14\delta\right)d_\e(x_0,x_1)\le
h(x_2'')-h(x_3'')\le (1+15\delta)d_\e(x_0,x_1),
\end{equation}
and
\begin{equation}\label{eq:height10'}
\left(1-\delta-2\eta\right)d_\e(x_0,x_1)\le\left(\frac{1}{1+\delta}-2\eta\right)d_\e(x_0,x_1)\le
h(x_1')-h(x_0')\le (1+\delta+2\eta)d_\e(x_0,x_1).
\end{equation}
So, let $x_3^{\#\#}$ be an ancestor of $x_2''$ at distance
$h(x_2'')-h(x_3'')-\left(16\delta+2\eta\right)d_\e(x_0,x_1)\in
[0,h(x_2'')]$ from $x_2''$, and let $x_0^\#$ be an ancestor of
$x_1'$ at distance
$h(x_2'')-h(x_3'')-\left(16\delta+2\eta\right)d_\e(x_0,x_1)\in
[0,h(x_1')]$ from $x_1'$. Then
\begin{equation}\label{eq:sharp}
d_\e(x_3'',x_3^{\#\#})\le \left(16\delta+2\eta\right)d_\e(x_0,x_1),
\end{equation}
and
\begin{multline}\label{eq:sharp2}
d_\e(x_0',x_0^\#)=\left|h(x_0')-h(x_0^\#)\right|=
\left|h(x_0')-\left(h(x_1')-h(x_2'')+h(x_3'')+\left(16\delta+2\eta\right)d_\e(x_0,x_1)\right)\right|
\\
\stackrel{\eqref{eq:height23''}\wedge\eqref{eq:height10'}}{\le}
2(16\delta+2\eta)d_\e(x_0,x_1).
\end{multline}
Moreover, $h(x_1)-h(x_0^\#)=h(x_2'')-h(x_3^{\#\#})$, so by
Lemma~\ref{lem:ancestor} we have
\begin{multline*}
\left(\frac{3}{1+\delta}-55\delta-7\eta\right)d_\e(x_0,x_1)\le
d_\e(x_0,x_3)-(55\delta+7\eta)d_\e(x_0,x_1)\stackrel{\eqref{eq:sharp}\wedge\eqref{eq:sharp2}}{\le}
d_\e(x_0^\#,x_3^{\#\#})\\\le d_\e(x_1',x_2'')\le
(1+8\delta+\eta)d_\e(x_0,x_1),
\end{multline*}
which is a contradiction since $\max\{\delta,\eta\}<1/200$.

Lastly, assume for contradiction that $(x_1,x_2,x_3)$ is $7\delta
d_\e(x_1,x_2)$-near a path-type configuration
$(x^{\prime\prime}_1,x^{\prime\prime}_2,x^{\prime\prime}_3)$. Then
since $h(x_1')\le h(x_2')$ we have
\begin{multline*}
\left(\frac{1}{1+\delta}-14\delta\right)d_\e(x_0,x_1)\le
d_\e(x_1'',x_2'')=h(x_1'')-h(x_2'')\le
\left(h(x_1'')-h(x_2'')\right)+\left(h(x_2')-h(x_1')\right)\\
\le d_\e(x_1'',x_1')+d_\e(x_2'',x_2')\le
(14\delta+2\eta)d_\e(x_0,x_1),
\end{multline*}
a contradiction.
\end{proof}

\begin{lemma} \label{lem:3-path-4}
Assume that $\e_n<\frac14$ for all $n$ and that $(x_0,x_1,x_2,x_3)$
is a $(1+\delta)$-approximate $P_3$ such that
$(x_2,x_1,x_0)$ is $\eta d_\e(x_0,x_1)$-near a tent-type configuration. If
$\max\{\delta,\eta\}<1/200$ then $(x_1,x_2,x_3)$ cannot be $7\delta d_\e(x_0,x_1)$ near a tent-type configuration.
\end{lemma}
\begin{proof}
Let $(x'_2,x'_1,x'_0)$ be a tent-type configuration that is $\eta d_\e(x_0,x_1)$-near $(x_2,x_1,x_0)$. Suppose for contradiction that $(x_1,x_2,x_3)$ is $7\delta d_\e(x_0,x_1)$-near a tent-type configuration $(x^{\prime\prime}_1,x^{\prime\prime}_2,x^{\prime\prime}_3)$.
Then $h(x'_1)-h(x'_2)=d_\e(x'_1,x'_2)$, whereas
$h(x^{\prime\prime}_2)-h(x^{\prime\prime}_1)= d_\e(x^{\prime\prime}_1,x^{\prime\prime}_2)$.
Taking the sum of these two inequalities we conclude that
\[ d_\e(x^{\prime\prime}_1,x^{\prime\prime}_2)+d_\e(x'_1,x'_2) \le d_\e(x_1',x_1'')+d_\e(x_2',x_2'')\le (2\eta+14\delta)d_\e(x_0,x_1). \]
At the same time, $\left(\frac{2}{1+\delta}-2\eta-14\delta\right)d_\e(x_0,x_1)\le d_\e(x^{\prime\prime}_1,x^{\prime\prime}_2)+d_\e(x'_1,x'_2)$, which leads to the desired contradiction.
\end{proof}

%Thus the resulting possible configuration of 3-paths are as in Table~\ref{tab:3-paths-types},
%and Figure~\ref{fig:3-paths-types}.

\begin{proof}[Proof of Lemma~\ref{lem:3-path classification}]
Since $(x_0,x_1,x_2,x_3)$ is a $(1+\delta)$-approximate $P_3$, we have $x_1\in\Mid(x_0,x_2,\delta)$, and
$x_2\in\Mid(x_1,x_3,\delta)$. Since the assumptions of  Theorem~\ref{thm:midpoints} hold, we can apply with $\eta=7\delta$
Lemmas~\ref{lem:3-path-1},~\ref{lem:3-path-2},~\ref{lem:3-path-3},~\ref{lem:3-path-4}, and
conclude that either $(x_0,x_1,x_2,x_3)$ or $(x_3,x_2,x_1,x_0)$ must be $35\delta d_\e(x_0,x_1)$-near a configuration of type $A$, $B$ or $C$.
\end{proof}

\subsection{Nonembeddability of vertically faithful $B_4$}
\label{sec:inembed-B4}

In what follows we need some standard notation on trees. As before, $B_n$ is the complete binary tree of height $n$;
the root of $B_n$ is denoted by $r$.
Denote by $I(B_n)$ the set of internal vertices of $B_n$, i.e., vertices of $B_n$ which are not the root or a leaf.
For a vertex $v$ in $\{r\}\cup I(B_n)$ we denote by $v_0$ and
$v_1$ its children. For $\alpha\in\{0,1\}^*$ (the set of finite sequences of '0' and '1') and $a\in\{0,1\}$ we
denote by $v_{\alpha a}=(v_\alpha)_a$.
%When $u$ is a vertex in a rooted tree, we denote by $h(u)$ the depth
%of $u$, i.e.,  the distance of $u$ from the root $r$.

The aim of the current section is to prove the following lemma.

\begin{lemma} \label{lem:B5}
Fix $0<\delta<1/400$ and let $f:B_4 \to (B_\infty,d_\e)$ be a $(1+\delta)$-vertically faithful embedding.
Then the distortion of $f$ satisfies
 $$
 \dist(f)\geq \frac{1}{500\delta+ \e_{h_0}},
$$
where   $h_0=\min_{x\in B_4} h(f(x))$.
\end{lemma}

\begin{comment}
We fix a $1+\delta$ vertically faithful embedding of $B_4$, $f:B_4 \to (b_\infty, d_\e)$, and denote $h_0=\min_{u\in B_4} h(f(u))$.
We also set $L=d_\e(f(r),f(r_0))$.
For $y\in I(B_4)$ we denote by $F(y)=(f(x),f(y),\{f(z),f(w)\})$,
where $x$ is the parent of $y$, and $z,w$ are the children of $y$ in $B_4$.
That is the $F(y)$ is the image of the 3-leaf star around $y$ in $B_4$.
The notation is chosen so as to signify that the roles of $z$ and $w$ are interchangeable.
\end{comment}

The proof of Lemma~\ref{lem:B5} is by a contradiction. By
Lemma~\ref{lem:fork classification0}, assuming the distortion of $f$
is small, all the $\delta$-forks in the $(1+\delta)$-vertically
faithful embedding must be of types $I$--$IV$. By exploring the
constrains implied by Lemma~\ref{lem:3-path classification} on how
those $\delta$-forks can be ``stitched" together, we reach the
conclusion that they are sufficiently severe to force any vertically
faithful embedding of $B_4$ to have a large contraction, and
therefore high distortion.

Fix $f:B_4\to (B_\infty,d_\e)$. For $u\in I(B_4)$ we denote by $\F(u)$ the fork in which $u$ is the center point,  i.e., if $v$ be the parent of $u$ in $B_4$, then
$$
\F(u)\eqdef(f(v),f(u),f(u_0),f(u_1)).
$$
We shall assume from now on that $f$ satisfies the assumptions of
Lemma~\ref{lem:B5}, i.e., that it satisfies~\eqref{eq:def vertical}
with $D=1+\delta$ for some $\delta<1/400$ and $\lambda>0$.

\begin{lemma} \label{lem:I,III->}
Fix $u\in B_4$ with $h(u)\in\{1,2\}$. If the fork $\F(u)$ is
$37\delta \lambda$-near a type $I$ or type $III$ configuration, then
there exists $w\in I(B_4)$ satisfying
\begin{equation}\label{eq:find w}
d_\e(f(w_0),f(w_1))\le (170\delta +\e_{h_0})\cdot 2\lambda.
\end{equation}
\end{lemma}
\begin{proof}
Let $v$ be the parent of $u$. Hence, $(f(v),f(u),f(u_{0}),f(u_{1}))$
is $35\delta(1+\delta) \lambda$-near a type $I$ or a type $III$ configuration.
Assume first that $(f(v),f(u),f(u_{0}),f(u_{1}))$ is $37\delta
\lambda$-near a type $I$ configuration. If both $(f(u_0),f(u),f(v))$
and $(f(u_1),f(u),f(v))$ were $37\delta\lambda$-near a path type
configuration then by Lemma~\ref{lem:opath,opath} (with
$\eta=37\delta$) we would have
\begin{equation}\label{eq:w is u}
d_\e(f(u_0),f(u_1))\le (339\delta+2\e_{h_0})(1+\delta)\lambda\le
(170\delta +\e_{h_0})\cdot 2\lambda,
\end{equation}
proving~\eqref{eq:find w} with $w=u$. The same conclusion holds when
$(f(v),f(u),f(u_{0}),f(u_{1}))$ is $37\delta \lambda$-near a type
$III$ configuration: in this case without loss of generality
$(f(v),f(u),f(u_{0}))$ is $37\delta \lambda$-near a tent-type
configuration and $(f(u_1),f(u),f(v))$ is $37\delta\lambda$-near a
path-type configuration. Using Lemma~\ref{lem:opath,opath} as above
we would arrive at the conclusion~\eqref{eq:w is u} if
$(f(u_0),f(u),f(v))$ were $37\delta \lambda$-near a path-type type
configuration. Thus, in both the type $I$ and type $III$ cases of
Lemma~\ref{lem:I,III->} we may assume that $(f(v),f(u),f(u_{0}))$ is
$37\delta \lambda$-near a tent-type configuration, and that, by
Lemma~\ref{lem:unique-midpoint-configuration},
$(f(v),f(u),f(u_{0}))$ is not $37\delta\lambda$-near a path-type
configuration, and $(f(u_0),f(u),f(v))$ is not
$37\delta\lambda$-near a path-type configuration or a tent-type
configuration.

 By
Lemma~\ref{lem:3-path classification} (and
Table~\ref{tab:3-paths-types}) $(f(u_{0c}),f(u_0),f(u),f(v))$ must
be $35\delta(1+\delta) \lambda$-near a type $B$ configuration for both
$c\in\{0,1\}$. This means that $(f(u_{0c}),f(u_0),f(u))$ are both
$35\delta(1+\delta) \lambda$-near a path-type configuration, and so by
Lemma~\ref{lem:opath,opath} (with $\eta=35\delta(1+\delta)$) we deduce that
$d_\e(f(u_{00}),f(u_{01}))\le(170\delta+ \e_{h_0})\cdot 2\lambda$.
\end{proof}

\begin{lemma}\label{lem:II->}
Fix $u\in B_4$ with $h(u)\in\{1,2\}$. If $\F(u)$ is $37\delta
\lambda$-near a type $II$ configuration then for both $b\in\{0,1\}$
either $\F(u_{b})$ is $99\delta \lambda$-near a type $II$
configuration, or $d_\e(f(u_{b0}),f(u_{b1}))\le 400\delta \lambda$.
\end{lemma}
\begin{proof}
Let $v$ be the parent of $u$. For both $c\in\{0,1\}$ we know that $(f(v),f(u),f(u_0),f(u_{0c}))$
is a $(1+\delta)$-approximate $P_3$, and therefore by Lemma~\ref{lem:3-path classification} either $(f(v),f(u),f(u_0),f(u_{0c}))$ or $(f(u_{0c}),f(u_0),f(u),f(v))$
is $35\delta(1+\delta)\lambda $-near a configuration of type $A$, $B$ or $C$. Note that since $(f(v),f(u),f(u_{0}))$ is  assumed to be $37\delta\lambda$-near a path-type configuration, we rule out the possibility that $(f(v),f(u),f(u_0),f(u_{0c}))$ is $35\delta(1+\delta)\lambda $-near a configuration of type $C$, since otherwise both $(f(v),f(u),f(u_{0}))$ and $(f(u_0),f(u),f(v))$ would be $37\delta\lambda $-near path-type configurations, contradicting Lemma~\ref{lem:unique-midpoint-configuration}. For the same reason we rule out the possibility that $(f(u_{0c}),f(u_0),f(u),f(v))$ is $35\delta(1+\delta)\lambda $-near a configuration of type $A$ or type $B$. An inspection of the three remaining possibilities shows that either $(f(u),f(u_0),f(u_{0c}))$ is $37\delta\lambda$-near a path-type configuration, or $(f(u_{0c}),f(u_0),f(u))$ is $37\delta\lambda$-near a tent-type configuration.

Now,

\begin{itemize}
\item If for both $c\in\{0,1\}$ we have that $(f(u),f(u_0),f(u_{0c}))$ is $37\delta\lambda$-near a path-type configuration, then $\F(u_0)$ is $37\delta\lambda$-near a
type $II$ configuration.
\item If for both $c\in\{0,1\}$ we have that $(f(u_{0c}),f(u_0),f(u))$ are $37\delta\lambda$-near a  tent-type configuration,
then by Lemma~\ref{lem:otent,otent} we have $d_\e(f(u_{01}),f(u_{00}))\le 400\delta$.
\item By Lemma~\ref{lem:path type =>path type}, the only way that
$(f(u),f(u_0),f(u_{00}))$ could be $37\delta\lambda$-near a path type configuration while at the same time
$(f(u_{01}),f(u_0),f(u))$ is $37\delta\lambda$-near a tent-type configuration (or vice versa), is that
$\F(u_0)$ is $99\delta \lambda$-near a type $II$ configuration. \qedhere
\end{itemize}
\end{proof}

\begin{lemma} \label{lem:IV->}
Fix $u\in B_4$ with $h(u)\in\{1,2\}$.
If $\F(u)$ is $35\delta(1+\delta) \lambda$-near a type $IV$ configuration, then there exists $b\in\{0,1\}$ such that
$\F(u_b)$ is $37\delta \lambda$-near a type $II$ configuration.
\end{lemma}
\begin{proof}
Let $v$ be the parent of $u$.
Without loss of generality $(f(u_0),f(u),f(v))$ is $35\delta(1+\delta) \lambda$-near a tent-type configuration.
By Lemma~\ref{lem:3-path classification} (using Lemma~\ref{lem:unique-midpoint-configuration} to rule out the remaining possibilities),
this means that for both $c\in\{0,1\}$ the quadruple $(f(u_{0c}),f(u_0),f(u),f(v))$ is
$35\delta(1+\delta)^2 \lambda$-near a type $C$ configuration, and
therefore $\F(u_0)$ is $35\delta(1+\delta)^2 \lambda$ near a type $II$ configuration.
\end{proof}

\begin{lemma} \label{lem:II,II}
Fix $u\in B_4$ with $h(u)\in\{0,1,2\}$.
If $\F(u_{0})$ and $\F(u_{1})$ are both $99\delta \lambda$-near type a  $II$ configuration then
\( d_\e(f(u_{0}),f(u_{1})) \le 1000\delta\lambda. \)
\end{lemma}
\begin{proof}
By our assumptions, $(f(u),f(u_0),f(u_{00}))$ is $99\delta \lambda$-near a path type configuration $(u',u'_0,u'_{00})$ and
$(f(u),f(u_1),f(u_{10}))$ is $99\delta \lambda$-near a path-type configuration $(u^{\prime\prime},u^{\prime\prime}_1,u^{\prime\prime}_{10})$. We may assume without loss of generality that $h(u'')-h(u_1')\le h(u')-h(u_0')$. We may therefore consider the ancestor $u_1^*$ of $u'$ such that $h(u')-h(u_1^*)=h(u^{\prime\prime})- h(u^{\prime\prime}_1)$, implying in particular that $h(u_1^*)\ge h(u_0')$ (recall that $u'_0$ is an ancestor of $u'$, and $u^{\prime\prime}_1$ is ancestor of $u^{\prime\prime}$). By Lemma~\ref{lem:ancestor} we have
\begin{equation}\label{eq:smooth use of lemma}
d_\e(u^*_1,u^{\prime\prime}_1)\le d_\e(u',u^{\prime\prime})\le 198\delta\lambda.
\end{equation}
 Hence,
\begin{equation}\label{eq:400}
h(u')-h(u_1^*)=d_\e(u',u_1^*)\stackrel{\eqref{eq:smooth use of lemma}}{\ge} d_\e(u',u_1'')-198\delta\lambda \ge d_\e(f(u),f(u_1))-394\delta\lambda\ge (1-394\delta)\lambda.
\end{equation}
But, we also know that
\begin{equation}\label{eq:201}
h(u')-h(u_0')=d_\e(u',u_0')\le d_\e(f(u),f(u_0))+198\delta\lambda\le (1+200\delta)\lambda.
\end{equation}
It follows from~\eqref{eq:400} and~\eqref{eq:201} that $d_\e(u_0',u_1^*)=h(u_1^*)-h(u_0')\le 601\delta\lambda$. Therefore,
\begin{equation*}
 d_\e(f(u_1),f(u_0)) \le d_\e(f(u_0),u'_0)+ d_\e(u'_0,u_1^*)+ d_\e(u^*_1,u^{\prime\prime}_1)+
 d_\e(u^{\prime\prime}_1, f(u_1))=1000\delta\lambda. \qedhere
\end{equation*}
\end{proof}

\begin{proof}[Proof of Lemma~\ref{lem:B5}]
We may assume that for all $u\in I(B_4)$ the fork $\F(u)$ is $35\delta(1+\delta) \lambda$-near a configuration of type $I$, $II$, $III$, or $IV$. Indeed, otherwise the proof is complete by Lemma~\ref{lem:fork classification0}.
If $\F(r_0)$ is $35\delta(1+\delta)\lambda$-near a type $I$ or type $III$ configuration, then by Lemma~\ref{lem:I,III->}
the proof is complete.
%\mnote{STOPPED HERE}
If $F(r_0)$ is $35\delta(1+\delta) \lambda$-near a type $IV$ configuration then by Lemma~\ref{lem:IV->}
there exists $b\in\{0,1\}$ such that $F(r_{0b})$ is $37\delta \lambda$-near  a type $II$ configuration. It therefore remains to deal with the case in which for some $u\in\{r_0,r_{0b}\}$ the fork $\F(u)$ is $37\delta \lambda$-near a type $II$ configuration.
Applying Lemma~\ref{lem:II->}, either we are done, or both $\F(u_0)$ and $\F(u_1)$ are $99\delta \lambda$-near a type $II$
configuration, but then by Lemma~\ref{lem:II,II} the proof of Lemma~\ref{lem:B5} is complete.
\end{proof}

\subsection{Nonembeddability of binary trees}
\label{sec:no-dich-Bn}

We are now in position to complete the proof of Theorem~\ref{lem:Bn->X}.

\begin{proof}[Proof of Theorem~\ref{lem:Bn->X}]
Write $\e_n=1/s(n)$, and $\e=\{\e_n\}_{n=0}^\infty$. Thus
$\{\e_n\}_{n=0}^\infty$ is non-increasing, $\{n\e_n\}_{n=0}^\infty$
is non-decreasing, and $\e_n\le 1/4$. We can therefore choose the
metric space $(X,d_X)=(B_\infty,d_\e)$. The identity embedding of
$B_n$ into the top $n$-levels of $B_\infty$ shows that $c_X(B_n)\le
s(n)$. It remains to prove the lower bound on $c_X(B_n)$. To this
end take an arbitrary injection $f:B_n\to X$ satisfying $\dist(f)\le
s(n)$, and we will now prove that
\begin{equation}
\label{eq:repeat dist lower} \dist(f)\ge s\left(\left\lfloor
\frac{n}{40 s(n)}\right\rfloor \right)\left(1-\frac{Cs(n)\log s(n)}
{\log n}\right).
\end{equation}

By adjusting the constant $C$ in~\eqref{eq:repeat dist lower}, we
may assume below that $n$ is large enough, say, $n\ge 100$. Write
$h_0=\lfloor n/(40s(n))\rfloor$ and define $X_{>h_0}= \{x\in
B_\infty:\ h(x)>h_0\}$. We claim that there exists a complete binary
subtree $T\subseteq B_n$ of height at least $\lceil n/3\rceil $,
such that we have $f(T)\subseteq X_{>h_0}$. Indeed, let
$h_{\min}=\min\{h(x):\; x\in f(B_n)\}$ and $h_{\max}=\max\{h(x):\;
x\in f(B_n)\}$. If $h_{\min}>h_0$ then
 $f(B_n)\subseteq X_{>h_0}$, and we can take $T=B_n$. So assume
that $h_{\min}<h_0$. Since $f$ is an injection it must satisfy
$h_{\max}\geq n$. Hence $\|f\|_{\Lip}\geq \frac{h_{\max}-h_{\min}}{2n}\geq \frac{n-h_0}{2n}\ge
\frac14$. Since $\dist(f)\leq s(n)$ we conclude that
$\|f^{-1}\|_{\Lip} \leq 4 s(n)$. It follows that, since  $\diam(X\setminus X_{>h_0})\leq 2h_0$, we have
 \( \diam \bigl (f^{-1}(X\setminus X_{>h_0}) \bigr )\leq 8 h_0 s(n)\le n/5. \) If the top $\lceil n/3\rceil $
 levels of $B_n$ are mapped into $X_{>h_0}$ then we are done, so
assume that there exists $u\in f^{-1}(X\setminus X_{>h_0})$ of depth at most $\le \lceil
n/3\rceil$. In this case  $f^{-1}(X\setminus X_{>h_0})$ must be
contained in the first $\lceil n/3\rceil+n/5<2n/3-1$ levels of
$B_n$, so we can take $T$ to be any subtree of $B_n$ contained in
the last $\lceil n/3\rceil$ levels of $B_n$.

Fix $\delta\in (0,1)$. By Theorem~\ref{lem:dich-vertical-Bn} (with
$t=4$, $D=s(n)$ and $\xi=\delta$), there exists a universal constant
$\kappa>0$ such that if $n\ge s(n)^{\kappa/\delta}$ then there
exists a mapping $\phi:B_4\to B_n$ with $\dist(\phi)\le 1+\delta$
such that $f\circ \phi$ is a $(1+\delta)$-vertically faithful
embedding of $B_4$ into $X_{>h_0}$. Choosing
$\delta=\kappa\frac{\log s(n)}{\log n}$, by increasing $C$
in~\eqref{eq:repeat dist lower} if necessary, we may assume that
$\delta<1/400$. Lemma~\ref{lem:B5} then implies
$$
(1+\delta)\dist(f)\ge\dist(f\circ\phi)\ge
\frac{1}{500\delta+\e_{h_0}} =\frac{1}{500\kappa\frac{\log
s(n)}{\log n}+\frac{1}{s\left(\lfloor
n/(40s(n))\rfloor\right)}}.
$$
%implying~\eqref{eq:repeat dist lower}.

The deduction of~\eqref{eq:subsequence} from~\eqref{eq:sharp c_X} is a simple exercise: if $s(n)=o(\log n/\log\log n)$ then we have $(s(n)\log s(n))/\log n=o(1)$. The desired claim will then follow once we check that \begin{equation}\label{eq:limsup}
\limsup_{n\to\infty} \frac{s\left(\lfloor
n/(40s(n))\rfloor\right)}{s(n)}=1.
\end{equation}
Indeed, if~\eqref{eq:limsup} failed then there would exist $\e_0\in (0,1)$ and $n_0\in \N$ such that for all $n\ge n_0$,
\begin{equation}\label{eq:to iterate}
s\left(\lfloor
n/\log n\rfloor\right)\le s\left(\lfloor
n/(40s(n))\rfloor\right)\le (1-\e_0)s(n).
\end{equation}
Iterating~\eqref{eq:to iterate}, it would follow that $s(n_j)\ge n_j^{\Omega(1)}$ for some subsequence $\{n_j\}_{j=1}^\infty$, a contradiction.
\end{proof}

\begin{proof}[Proof of Theorem~\ref{thm:no local rigidity}] The proof is identical to the above argument: all one has to notice is that when $s(n)=D$ for all $n\in \N$ the resulting metric $d_{\e}$ on $B_\infty$ is $D$-equivalent to the original shortest path metric on $B_\infty$. In this case, if $c_X(B_n)\le D-\e$ then the bound~\eqref{eq:repeat dist lower} implies that $n\le D^{CD^2/\e}$.
\end{proof}

%Note that Theorem~\ref{thm:Btree-non-dich} is an immediate corollary of Lemma~\ref{lem:Bn->X}.

%\subsection{Proofs of the approximate classification results}

%\subsubsection{Some preparatory stitching lemmas}

%\subsubsection{Proof of Lemma~\ref{thm:midpoints}}

%\subsubsection{Proof of Lemma~\ref{lem:fork classification0}}

%\subsubsection{Proof of Lemma~\ref{lem:3-path classification}}

\section{Discussion and open problems}\label{sec:open}

A very interesting question that arises naturally from
Theorem~\ref{thm:convexity-coincide} and is also a part of the Ribe
program, is finding a metric characterization of $q$-smoothness. A
Banach space $(X,\|\cdot\|_X)$ is called $q$-smooth if it admits an
equivalent norm ${\tb} \cdot {\tb}$ such that there is a constant
$S>0$ satisfying:
$$
{\tb} x{\tb} =1\ \wedge \   y\in X\implies \frac{{\tb} x+y{\tb}
+{\tb} x-y{\tb} }{2}\le 1+S{\tb} y{\tb} ^q.
$$
A Banach space $X$ is $p$-convex if and only if its dual space $X^*$
is $q$-smooth, where $\frac{1}{p}+\frac{1}{q}=1$~\cite{Lind63}. It
is known that a Banach space $X$ is $p$-convex for some $p<\infty$
(i.e., superreflexive) if and only if it is $q$-smooth for some
$q>1$ (this follows from~\cite{James72,Pisier-martingales}). Hence
Bourgain's metric characterization of superreflexivity can be viewed
as a statement about uniform smoothness as well. However, we still
lack a metric characterization of the more useful notion of
$q$-smoothness. Trees are natural candidates for finite metric
obstructions to $p$-convexity, but it is unclear what would be the
possible finite metric witnesses to the ``non-$q$-smoothness" of a
metric space.

$H$-trees are geometric objects that are quite simple
combinatorially, yet as we have seen, they have interesting
bi-Lipschitz properties. It would therefore be of interest to
investigate the geometry of $H$-trees for its own right. In
particular, what is the $L_1$ distortion of an $H$-tree? How close
can an $H$-tree be to a metric of negative type?

%\cite{charlie-soda,JS09,Gro83,Enf78,JS72,OS94,MP76,Pis73,Jam64,Men09}
%\cite{ckn}
%\cite{Bau07}
%\cite{Dvo60}
\bibliographystyle{abbrvurl}
\bibliography{dich}

\end{document}